%&amstex          
\input amstex\documentstyle{amsppt}  
\pagewidth{12.5cm}\pageheight{19cm}\magnification\magstep1
\topmatter
\title Non-unipotent representations and categorical centres\endtitle
\author G. Lusztig\endauthor
\address{Department of Mathematics, M.I.T., Cambridge, MA 02139}\endaddress
\thanks{Supported by NSF grant DMS-1566618.}\endthanks
\endtopmatter   
\document
\define\dz{\dot z}
\define\da{\dagger}
\define\btco{\bar{\ti{\co}}}
\define\tco{\ti{\co}}
\define\cdo{\cdot}

\define\Irr{\text{\rm Irr}}

\define\Bpq{\Bumpeq}

\define\du{\dot u}
\define\dw{\dot w}

\define\dy{\dot y}
\define\dZ{\dot Z}

     \define\bco{\bar{\co}}

\define\mpb{\medpagebreak}

\define\bh{\bar h}

\define\frl{\forall}

\define\si{\sim}
\define\wt{\widetilde}
\define\sqc{\sqcup}

\define\qua{\quad}

\define\bG{\bar G}

\define\bX{\bar X}
\define\bZ{\bar Z}
\define\lb{\linebreak}

\define\bin{\binom}
\define\op{\oplus}
   
\define\part{\partial}
\define\emp{\emptyset}

\define\ra{\rangle}
\define\n{\notin}
\define\iy{\infty}
\define\m{\mapsto}
\define\do{\dots}
\define\la{\langle}
\define\bsl{\backslash}

\define\lra{\leftrightarrow}

\define\sub{\subset}    
\define\bxt{\boxtimes}
\define\T{\times}
\define\ti{\tilde}
\define\nl{\newline}
\redefine\i{^{-1}}

\define\un{\underline}
\define\ov{\overline}
\define\ot{\otimes}
\define\bbq{\bar{\QQ}_l}

\define\Ad{\text{\rm Ad}}
\define\Hom{\text{\rm Hom}}

\define\tr{\text{\rm tr}}

\define\supp{\text{\rm supp}}
\define\card{\text{\rm card}}

\define\di{\diamond}
\redefine\spa{\spadesuit}

\define\a{\alpha}
\redefine\b{\beta}
\redefine\c{\chi}
\define\g{\gamma}
\redefine\d{\delta}
\define\e{\epsilon}
\define\et{\eta}

\redefine\o{\omega}
\define\p{\pi}
\define\ph{\phi}

\define\r{\rho}
\define\s{\sigma}
\redefine\t{\tau}

\define\k{\kappa}
\redefine\l{\lambda}
\define\z{\zeta}
\define\x{\xi}

\define\vt{\vartheta}

\redefine\G{\Gamma}
\redefine\D{\Delta}

\define\Si{\Sigma}

\define\Ph{\Phi}
\define\Ps{\Psi}

\define\boc{\bold c}

\define\ee{\bold e}

\define\kk{\bold k}
\define\mm{\bold m}

\define\pp{\bold p}

\define\rr{\bold r}

\redefine\tt{\bold t}

\define\ww{\bold w}

\redefine\AA{\bold A}
\define\BB{\bold B}
\define\CC{\bold C}
\define\DD{\bold D}

\define\FF{\bold F}

\define\HH{\bold H}

\define\JJ{\bold J}

\define\LL{\bold L}

\define\NN{\bold N}

\define\QQ{\bold Q}

\define\TT{\bold T}
\define\UU{\bold U}

\define\ZZ{\bold Z}
\define\XX{\bold X}

\define\ca{\Cal A}
\define\cb{\Cal B}
\define\cc{\Cal C}
\define\cd{\Cal D}
\define\ce{\Cal E}
\define\cf{\Cal F}

\define\ch{\Cal H}
\define\ci{\Cal I}

\define\cl{\Cal L}
\define\cm{\Cal M}

\define\co{\Cal O}
\define\cp{\Cal P}

\define\car{\Cal R}
\define\cs{\Cal S}

\define\cu{\Cal U}
\define\cv{\Cal V}
\define\cw{\Cal W}
\define\cz{\Cal Z}
\define\cx{\Cal X}
\define\cy{\Cal Y}

\define\fb{\frak b}
\define\fc{\frak c}

\define\fh{\frak h}

\define\fl{\frak l}

\define\fo{\frak o}

\define\fs{\frak s}

\define\fA{\frak A}

\define\fD{\frak D}

\define\fL{\frak L}

\define\fR{\frak R}
\define\fS{\frak S}

\define\fZ{\frak Z}

\define\tb{\ti b}
\define\tc{\ti c}

\define\tf{\ti f}
\define\tg{\ti g}
\define\tih{\ti h}
\define\tj{\ti j}
\define\tk{\ti k}

\define\tit{\ti t}

\define\ty{\ti y}

\define\tG{\ti G}

\define\tK{\ti K}
\define\tL{\ti L}

\define\tP{\ti P}

\define\tV{\ti V}

\define\tX{\ti X}

\define\sha{\sharp}

\define\bc{\bar c}
\define\bp{\bar p}

\define\bul{\bullet}

\define\che{\check}
\define\cha{\che{\a}}

\define\cir{\circ}

\define\tcb{\ti{\cb}}

\define\BBD{BBD}
\define\BFO{BFO}
\define\DL{DL}
\define\ENO{ENO}
\define\KL{KL}
\define\ORA{L1}
\define\ICM{L2}
\define\CELLSII{L3}
\define\CSI{L4}
\define\CSII{L5}
\define\CSIII{L6}
\define\TENS{L7}
\define\HEC{L8}
\define\CDGVI{L9}
\define\CDGVII{L10}
\define\CDGVIII{L11}
\define\CDGIX{L12}
\define\CDGX{L13}
\define\CONV{L14}
\define\URE{L15}
\define\MONO{L16}
\define\MUG{Mu}
\define\YO{Yo}

\head Introduction\endhead
\subhead \endsubhead   
\subhead 0.1\endsubhead
Let $\kk$ be an algebraic closure of the finite field with $p$ elements. 
Let $G$ be a connected reductive
group over $\kk$. We denote by $F_q$ the subfield of $\kk$ with exactly $q$ elements; here $q$ is a 
power of $p$. Let $F:G@>>>G$ be the Frobenius map for an $F_q$-rational structure on $G$.
We fix a prime number $l$ different from $p$. Let $\Irr(G^F)$ be the set of isomorphism classes of 
irreducible representations (over $\bbq$) of the finite group $G^F=\{g\in G;F(g)=g\}=G(F_q)$. 
In \cite{\ICM} I gave a parametrization of $\Irr(G^F)$ in terms of the
group of type dual to that of $G$. (For ``most'' representations in $\Irr(G^F)$ this has benn
already done in \cite{\DL}.) For the part of $\Irr(G^F)$ consisting of unipotent representations in a 
fixed two-sided cell of $W$ (with $G$ assumed to be $F_q$-split) the parametrization was in terms of a 
set $M(\G)$ where $\G$ is a certain finite group associated to the two-sided cell and $M(\G)$ is the 
set of simple objects (up to isomorphism) of the category $Vec_\G(\G)$
of $\G$-equivariant vector bundles on $\G$ 
(here $\G$ acts on $\G$ by conjugation). In the early 1990's, Drinfeld pointed out to me that the 
category $Vec_\G(\G)$ can be interpreted as the categorical centre of the monoidal category of finite 
dimensional representations of $\G$. (The notion of categorical centre of a monoidal category is
due to Joyal, Street, Majid and Drinfeld.) This suggested that one should be able the
reformulate the parametrization of $\Irr(G^F)$ in terms of categorical centres of suitable
monoidal categories associated with $G$. 
This is achieved in the present paper, except that we must allow certain
twisted categorical centres instead of usual categorical centres.
Note that in our approach the representation theory of $G(F_q)$ cannot be separated 
from the theory of character sheaves on $G$ which appears as the limit of the first theory when
$q$ tends to $1$; in particular we also obtain the parametrization of character sheaves on $G$ 
in terms of categorical centres (no twisting needed in this case).

Earlier results of this type were known in the following cases: 

(i) the case \cite{\BBD} of character sheaves on $G$ (with centre assumed to be connected and with
$\kk$ replaced by $\CC$);

(ii) the case \cite{\CONV} of unipotent character sheaves on $G$;

(iii) the case \cite{\URE} of unipotent representations of $G^F$;

(iv) the case \cite{\MONO} of not necessarily unipotent character sheaves on $G$.
\nl
The papers \cite{\URE},\cite{\MONO} were generalizations of \cite{\CONV} in different directions; the 
present paper is a common generalization of \cite{\URE},\cite{\MONO}; the methods used in 
(ii),(iii),(iv) and the present paper are quite different from those used in (i) which relied
on techniques not available in positive characteristic.

Let $\BB$ be a Borel subgroup of $G$ and let $\TT$ be a maximal torus of $\BB$. 
In this subsection we assume that $F(\BB)=\BB$, $F(\TT)=\TT$. Let $W$ be the Weyl group of $G$ 
with respect to $\TT$. Let $\fs$ be an indexing set for the isomorphism classes of Kummer local 
systems (over $\bbq$); note that $W$ acts naturally on $\fs$.

Let $\ca=\ZZ[v,v\i]$ where $v$ is an indeterminate. A key role in this paper is played by an
$\ca$-algebra $\HH$ (without $1$ in general) which has $\ca$-basis $\{T_w1_\l;w\in W,\l\in\fs\}$
and multiplication defined in 1.5 (see also \cite{\CDGVI, 31.2}). This is a monodromic version of the 
usual Hecke algebra of $W$, closely related to an algebra defined in \cite{\YO}; it contains the
usual Hecke algebra as a subalgebra.
Now $\HH$ has a canonical basis, two-sided cells and an asymptotic version $H^\iy$ 
(introduced in \cite{\CDGVII},\cite{\MONO}) which generalize the analogous notions for the usual
Hecke algebra, see \cite{\KL}, \cite{\CELLSII}; the two-sided cells form a partition of
$W\T\fs$ and we have $H^\iy=\op_\boc H^\iy_\boc$ as rings ($\boc$ runs over the two-sided cells and
each $H^\iy_\boc$ is a ring with $1$). For any $\boc$, $H^\iy_\boc$ admits a category version
(for which $H^\iy$ is the Grothendieck group) which is a semisimple monoidal category $\cc^\boc$ with 
finitely many simple objects (up to isomorphism) indexed by the elements of $\boc$, see \S5. (In the 
case where $\boc\sub W\T\{1\}$, this reduces to the monoidal category defined is \cite{\TENS}.) Now 
$\cc^\boc$ has a well defined categorical centre which is again a semisimple abelian category. 
Note that $F$ acts naturally on $\fs$ and on $W$ hence on $W\T\fs$; this induces an action
of $F$ on the set of two-sided cells. If $\boc$ is a two-sided cell such that $F(\boc)=\boc$ then $F$ 
defines an equivalence of categories $\cc^\boc@>>>\cc^\boc$ and one can define the notion of 
$F$-centre of $\cc^\boc$ (see 5.5) which is a twisted version of the usual centre; it is a
semisimple abelian category. We denote by $[\boc]$ the set of isomorphism classes of simple objects 
of this category (a finite set).

Our main result is that $\Irr(G^F)$ is in natural bijection with $\sqc_\boc[\boc]$
(disjoint union over all $F$-stable two-sided cells $\boc$). (See Theorem 7.3.)
In the case where $\boc\sub W\T\{1\}$, this reduces to the main result in \cite{\URE}.

The fact that the asymptotic Hecke algebra $\HH^\iy$ plays a role in the classification is perhaps
not surprising since its non-monodromic versions appeared implicitly in the arguments of \cite{\ORA}, 
through the traces of their canonical basis elements in their various simple modules (the algebras 
themselves were not defined at the time where \cite{\ORA} was written).

Many arguments in this paper follow very closely the arguments in \cite{\MONO}; we generalize
them by taking into account also the arguments in \cite{\URE}.
We have written the proofs in such a way that they apply at the same time in the case of
character sheaves on a connected component of a possibly disconnected algebraic group 
with identity component $G$. In this case, the classification involves 
twisted categorical centers, unlike that for the character sheaves on $G$.

We plan to show elsewhere that the parametrization of $\Irr(G^F)$ given in \cite{\ICM}
can be deduced from the main result of this paper. 

\subhead 0.2\endsubhead
{\it Notation.} 
Let $\NN^*=\{n\in\ZZ-p\ZZ;n\ge1\}$. Let $T$ be a torus over $\kk$. For $n\in\NN^*$ let 
$T_n=\{t\in T;t^n=1\}$; we have $\sha(T_n)=n^{\dim T}$. For $n,n'$ in $\NN^*$ such that 
$n'/n\in\ZZ$ we have a surjective homomorphism $N^{n'}_n:T_{n'}@>>>T_n$, $t\m t^{n'/n}$.
Hence we can form the projective limit $T^\iy$ of the groups $T_n$ with $n\in\NN^*$ (a profinite 
abelian group). Then for any $n\in\NN^*$, $T_n$ is naturally a quotient of $T^\iy$.

\mpb

All algebraic varieties are over $\kk$. We denote by $\pp$ the algebraic variety consisting of a 
single point. 
For an algebraic variety $X$ we write $\cd(X)$ for the bounded derived category of constructible
$\bbq$-sheaves on $X$. Let $\cm(X)$ be the subcategory of 
$\cd(X)$ consisting of perverse sheaves on $X$. 
For $K\in\cd(X)$ and $i\in\ZZ$ let $\ch^iK$ be the $i$-th cohomology sheaf of $K$ and let $K^i$ be
the $i$-th perverse cohomology sheaf of $K$. Let $\fD(K)$ be the Verdier dual of $K$.
For any constructible sheaf $\ce$ on $X$ let $\ce_x$ be the stalk of $\ce$ at $x\in X$.
If $X$ has a fixed $F_q$-structure $X_0$, we denote by 
$\cd_m(X)$ what in \cite{\BBD, 5.1.5} is denoted by $\cd_m^b(X_0,\bbq)$;
let $\cm_m(X)$ be the corresponding category of mixed perverse sheaves.
In this paper we often encounter maps of algebraic varieties which are not morphisms but 
only quasi-morphisms (as in \cite{\URE, 0.3}). For such maps the usual operations with
derived categories are defined as in \cite{\URE, 0.3}.

Note that if $K\in\cd_m(X)$ then $K$ can be viewed as an object of $\cd(X)$ denoted again by $K$.
If $K\in\cm_m(X)$ and $h\in\ZZ$, we denote by $gr_h(K)$ the subquotient of
pure weight $h$ of the weight filtration of $K$.
If $K\in\cd_m(X)$ and $i\in\ZZ$ we write $K\la i\ra=K[i](i/2)$ where $[i]$ is a shift and $(i/2)$
is a Tate twist; we write $K^{\{i\}}=gr_i(K^i)(i/2)$. 
If $K$ is a perverse sheaf on $X$ and $A$ is a simple perverse sheaf on $X$ we write $(A:K)$ for
the multiplicity of $A$ in a Jordan-H\"older series of $K$.
If $C\in\cd_m(X)$ and $\{C_i;i\in I\}$ is a family of objects of $\cd_m(X)$ then the relation
$C\Bpq\{C_i;i\in I\}$ is as in \cite{\MONO, 0.2}.

Let $\bar{}:\ca@>>>\ca$ be the ring
homomorphism such that $\ov{v^m}=v^{-m}$ for any $m\in\ZZ$. If $f\in\QQ[v,v\i]$ and $j\in\ZZ$ we 
write $(j;f)$ for the coefficient of $v^j$ in $f$.

Let $\cb$ be the variety of Borel subgroups of $G$. For any $B\in\cb$ 
let $U_B$ be the unipotent radical of $B$. 
In this paper we fix a Borel subgroup $\BB$ of $G$ and a maximal torus $\TT$ of $\BB$.
Let $\UU=U_\BB$. Let $\nu=\dim\UU=\dim\cb$, $\r=\dim\TT$, $\D=\dim G=2\nu+\r$. 

For any algebraic variety $X$ let $\fL=\fL_X=\a_!\bbq\in\cd(X)$ where $\a:X\T\TT@>>>X$ is the obvious 
projection. When $X$ and $T$ are defined over $\FF_q$, $\fL$ is naturally an object of $\cd_m(X)$.

Unless otherwise specified, all vector spaces are over $\bbq$;  in particular, all 
representations of finite 
 groups are assumed to be in (finite dimensional) $\bbq$-vector spaces. 

\head Contents\endhead

1. The monodromic Hecke algebra and its asymptotic version.

2. The group $\tG$.

3. Sheaves on $\tcb^2$.

4. Sheaves on $Z_s$.

5. The monoidal category $\cc^\boc\tcb^2$.

6. Truncated induction, truncated restriction, truncated convolution.

7. Equivalence of $\cc^\boc\tG_s$ with the $\ee^s$-centre of $\cc^\boc\tcb^2$.

\head 1. The monodromic Hecke algebra and its asymptotic version\endhead
\subhead 1.1\endsubhead
Let $N\TT$ be the normalizer of $\TT$ in $G$, let $W=N\TT/\TT$ be the Weyl group and let 
$\k:N\TT@>>>W$ be the obvious homomorphism. For $w\in W$ we set $G_w=\UU\k\i(w)\UU$ so that 
$G=\sqc_{w\in W}G_w$; let $\co_w=\{(x\BB x\i,y\BB y\i);x\in G,y\in G,x\i y\in G_w\}$ so that
$\cb\T\cb=\sqc_{w\in W}\co_w$. For $w\in W$ let $\bG_w$ be the closure of $G_w$ in $G$; we have 
$\bG_w=\cup_{y\le w}G_y$ for a well defined partial order $\le$ on $W$. Let $\bco_w$ be the 
closure of $\co_w$ in $\cb\T\cb$. Now $W$ is a (finite) Coxeter group with length function 
$w\m|w|=\dim\co_w-\nu$ and with set of generators $S=\{\s\in W;|\s|=1\}$; it acts on $\TT$ by 
$w:t\m w(t)=\o t\o\i$ where $\o\in\k\i(w)$.

\subhead 1.2\endsubhead
Let $R\sub\Hom(\TT,\kk^*)$ be the set of roots of $G$ with respect to $\TT$.
Now $W$ acts on $R$ by $w:\a\m w(\a)$ where $(w(\a))(t)=\a(w\i(t))$ for $t\in\TT$.
Let $R^+$ be the set of $\a\in R$ such that the corresponding root subgroup is
contained in $G$. For $\a:\TT@>>>\kk^*$ we denote by $\cha:\kk^*@>>>\TT$ the 
corresponding coroot and by
$\s_\a$ the corresponding reflection in $W$. For any $\s\in S$ let $\UU_\s$ be the unique root 
subgroup of $\UU$ with respect to $\TT$ such that $\UU^-_\s:=\o\UU_\s\o\i\not\sub\UU$ for some/any
$\o\in\k\i(\s)$. Let $\a_\s:\TT@>>>\kk^*$ be the root corresponding to $\UU_\s$; then the coroot 
$\cha_\s:\kk^*@>>>\TT$ is well defined. 

For any $\s\in S$ we fix an element $\x_\s\in\UU_\s-\{1\}$; there is a unique
$\x'_\s\in\UU^-_\s-\{1\}$ such that $\x_\s\x'_\s\x_\s=\x'_\s\x_\s\x'_\s\in\k\i(\s)\sub N\TT$;
the two sides of the last equality are denoted by $\dot\s$.
We have $\k(\dot\s)=\s$ and $\dot\s^2=\cha_\s(-1)$. 
For any $w\in W$ we define $\dw\in N\TT$ by $\dw=\dot\s_1\dot\s_2\do\dot\s_r$ where
$w=\s_1\s_2\do\s_r$ with $r=|w|,\s_j\in S$; note that, by a result of Tits, $\dw$ 
is well defined. Let $N_0\TT$ be the 
subgroup of $N\TT$ generated by $\{\dot\s;\s\in S\}$. This is a finite subgroup of $N\TT$ 
containing $\dw$ for any $w\in W$. Let $\k_0:N_0\TT@>>>W$ be the restriction of 
$\k:N\TT@>>>W$.

\subhead 1.3\endsubhead
For $n\in\NN^*$ let $\fs_n=\Hom(\TT_n,\bbq^*)$; we have $\sha(\fs_n)=n^\r$.
For $n,n'$ in $\NN^*$ such that $n'/n\in\ZZ$, the surjective homomorphism 
$N^{n'}_n:\TT_{n'}@>>>\TT_n$, $t\m t^{n'/n}$ induces an imbedding $\fs_n\sub\fs_{n'}$, 
$\l\m\l N^{n'}_n$. Hence we can form the union $\fs_\iy=\cup_{n\in\NN^*}\fs_n$ (a countable 
abelian group). Then for any $n\in\NN^*$, $\fs_n$ is a subgroup of $\fs_\iy$. Note also that 
$\fs_\iy$ is the group of homomorphisms $\TT^\iy@>>>\bbq^*$ which factor through $\TT_n$ for some 
$n\in\NN^*$. For any $\l\in\fs_\iy$ there is a well defined local system $L_\l$ on $\TT$ such that
for some/any $n\in\NN^*$ for which $\l\in\fs_n$, $L_\l$ is equivariant for the $\TT$-action 
$t_1:t\m t_1^nt$ on 
$\TT$ and the natural $\TT_n$ action the stalk of $L_\l$ at $1$ is through the character $\l$. For
$\l,\l'\in\fs_\iy$ we have canonically $L_\l\ot L_{\l'}=L_{\l\l'}$; for $\l\in\fs_\iy$ we have 
canonically $L_\l^*=L_{\l\i}$; here $()^*$ denotes the dual local system.

The $W$-action on $\TT$ restricts to a $W$-action on $\TT_n$ for any $n\in\NN^*$. This induces a
$W$-action on $\TT^\iy$, a $W$-action on $\fs_n$ for any $n\in\NN^*$; for $\l\in\fs_n$, $w\in W$ 
and $t\in\TT_n$ we have $(w(\l))(t)=\l(w\i(t))$. There is a unique $W$-action of $\fs_\iy$ which 
for any $n\in\NN^*$ restricts to the $W$-action on $\fs_n$ just described.
We set $I=W\T\fs_\iy$; for 
$w\in W,\l\in\fs_\iy$ we write $w\cdo\l$ instead of $(w,\l)$. 

\subhead 1.4 \endsubhead  
If $\a\in R$, the coroot $\cha:\kk^*@>>>\TT$ restricts to
a homomorphism $\kk^*_n@>>>\TT_n$ for any $n\in\NN^*$ and by passage to
projective limits, this induces a homomorphism $\cha^\iy:\kk^\iy@>>>\TT^\iy$
(notation of 0.2). Let $\l\in\fs_\iy$. We say that $\a\in R_\l$ if the 
composition $\kk^\iy@>\cha^\iy>>\TT^\iy@>\l>>\bbq^*$ is identically $1$ or
equivalently if $\cha^*L_\l\cong\bbq$ as local systems on $\kk^*$. Note that
for $w\in W$ we have $w(R_\l)=R_{w(\l)}$. Let $R_\l^+=R_\l\cap R^+$, $R_\l^-=R_\l-R_\l^+$.
Let $W_\l$ be the subgroup of $W$ generated by $\{\s_\a;\a\in R_\l\}$. We have $W_\l=W_{\l\i}$. 
Let $W'_\l=\{w\in W;w(\l)=\l\}$. We have $W_\l\sub W'_\l$. As in \cite{\CSI, 5.3}, there is a 
unique Coxeter group structure on $W_\l$ with length function $W_\l@>>>\NN$, 
$w\m|w|_\l=\sha\{\a\in R_\l^+;w(\a)\in R^-_\l\}$; note that, if $w\in W_\l$ and 
$w=\s_1\s_2\do\s_r$ is any reduced expression of $w$ in $W$, then 
$$|w|_\l=\card\{i\in[1,r];\s_r\do\s_{i+1}\s_i\s_{i+1}\do\s_r\in W_\l\}.$$

\subhead 1.5 \endsubhead
For $n\in\NN^*$ we set $I_n=\{w\cdo\l\in I;\l\in\fs_n\}$. 
As in \cite{\CDGVI, 31.2}, let $\HH_n$ be 
the associative $\ca$-algebra with generators $T_w (w\in W)$, $1_\l (\l\in\fs_n)$ and 
relations:

$1_\l1_{\l'}=\d_{\l,\l'}1_\l$ for $\l,\l'\in\fs_n$;

$T_wT_{w'}=T_{ww'}$ if $w,w'\in W$ and $|ww'|=|w|+|w'|$;

$T_w1_\l=1_{w(\l)}T_w$ for $w\in W,\l\in\fs_n$;

$T_\s^2=v^2T_1+(v^2-1)\sum_{\l\in\fs_n;\s\in W_\l}T_\s1_\l$ for $\s\in S$;

$T_1=\sum_{\l\in\fs_n}1_\l$.
\nl
The algebra $\HH_n$ is closely related to the algebra introduced by Yokonuma \cite{\YO}. (It 
specializes to it under $v=\sqrt{q},n=q-1$ where $q$ is a power of a prime; this is shown in
\cite{\CDGVII, \S35}.) Note that $T_1$ is the unit element of $\HH_n$. In \cite{\CDGVII, 31.2} it is
shown that $\{T_w1_\l;w\cdo\l\in I_n\}$ is an $\ca$-basis of $\HH_n$. (In \cite{\MONO,1.7} we 
write $\HH$ instead of $\HH_n$, but here we shall not do so.)

Now, for $\s\in S$, $T_\s$ is invertible in $\HH_n$; indeed, we have
$$T_\s\i=v^{-2}T_\s+(1-v^{-2})(\sum_{\l\in\fs_n;\s\in W_\l}1_\l).$$
It follows that $T_w$ is invertible in $\HH_n$ for any $w\in W$. As shown in \cite{\CDGVI, 31.3}, 
there is a unique ring homomorphism $\HH_n@>>>\HH_n$, $h\m\bh$ such that $\ov{T_w}=T_{w\i}\i$ for 
any $w\in W$ and $\ov{f1_\l}=\bar f1_\l$ for any $f\in\ca$, $\l\in\fs_n$. It is an involution 
called the {\it bar involution}.

If $n,n'\in\NN^*$ and $n'/n\in\ZZ$, then $I_n\sub I_{n'}$ and the $\ca$-linear map 
$j_{n,n'}:\HH_n@>>>\HH_{n'}$ given by $T_w1_\l\m T_w1_\l$ for $w\cdo\l\in I_n$ is an $\ca$-algebra 
imbedding which does not necessarily preserve the unit element. Let $\HH$ be the union of all 
$\HH_n$ for various $n\in\NN^*$ according to the imbeddings $j_{n,n'}$above. Then $\HH$ is an 
$\ca$-algebra without $1$ in general; it has an $\ca$-basis 
$\{T_w1_\l=1_{w(\l)}T_w;w\cdo\l\in I\}$. If $n\in\NN^*$, then $\HH_n$ is the $\ca$-submodule of 
$\HH$ with basis $\{T_w1_\l;w\cdo\l\in I_n\}$; it is an $\ca$-subalgebra of $\HH$. The algebra 
$\HH_n$ has been studied in \cite{\CDGVII} and \cite{\MONO, 1.7}. We shall often refer to 
{\it loc.cit.} for properties of $\HH$ which in {\it loc.cit.} are stated for $\HH_n$ with $n$ 
fixed and which imply immediately the corresponding properties of $\HH$.

We show that, if $n,n'\in\NN^*$ and $n'/n\in\ZZ$, then $j_{n,n'}:\HH_n@>>>\HH_{n'}$ is compatible
with the bar-involution on $\HH_n$ and $\HH_{n'}$. 
It is enough to show that $j_{n,n'}(\ov\x)=\ov{j_{n,n'}(\x)}$ for $\x=1_\l,\l\in\fs_n$ or
$\x=T_\s$, $\s\in S$. The case where $\x=1_\l,\l\in\fs_n$ is immediate. For $\s\in S$ we have 
$j_{n,n'}(T_\s)=T_\s\sum_{\l\in\fs_n}1_\l$, hence
$$\align&j_{n,n'}(\ov{T_\s})=j_{n,n'}(v^{-2}T_\s+(1-v^{-2})(\sum_{\l\in\fs_n;\s\in W_\l}1_\l))\\&
=v^{-2}T_\s\sum_{\l\in\fs_n}1_\l+(1-v^{-2})(\sum_{\l\in\fs_n;\s\in W_\l}1_\l)
=T_\s\i\sum_{\l\in\fs_n}1_\l=\ov{j_{n,n'}(T_\s)},\endalign$$
as desired.
It follows that there is a unique ring homomorphism $\HH@>>>\HH$, $h\m\bar h$, whose restriction 
to $\HH_n$ (for any $n\in\NN^*$) is the bar involution. This has square $1$ and is again called 
the bar involution. 

The $\ca$-linear map $\HH@>>>\HH$, $h\m\tih$ given by $T_w1_\l\m T_w1_{\l\i}$ for
$w\cdo\l\in I$ is an algebra involution. The $\ca$-linear map $\HH@>>>\HH$, $h\m h^\flat$, given 
by $T_w1_\l\m1_\l T_{w\i}$ is an involutive algebra antiautomorphism. (See \cite{\CDGVII, 32.19}.)

\subhead 1.6\endsubhead    
As in \cite{\CDGVII, 34.4}, for any $w\cdo\l\in I$ there is a unique element $c_{w\cdo\l}\in\HH$ 
such that
$$c_{w\cdo\l}=\sum_{y\in W}p_{y\cdo\l,w\cdo\l}v^{-|y|}T_y1_\l$$
where $p_{y\cdo\l,w\cdo\l}\in v\i\ZZ[v\i]$ if $y\ne w$, 
$p_{w\cdo\l,w\cdo\l}=1$ and $\ov{c_{w\cdo\l}}=c_{w\cdo\l}$. 
For $\l\in\fs_\iy$, $y',w'$ in $W_\l$ let $P^\l_{y',w'}$ be the polynomial defined in \cite{\KL} 
in terms of the Coxeter group $W_\l$; let 
$$p^\l_{y',w'}=v^{-|w'|_\l+|y'|_\l}P^\l_{y',w'}(v^2)\in\ZZ[v\i].$$
Let $w\cdo\l\in I$. From \cite{\ORA, 1.9(i)} we see that $wW_\l$ contains a
unique element $z$ such that $|z|$ is minimum; we write $z=\min(wW_\l)$; we 
have $w=zw'$ with $w'\in W_\l$. We have
$$c_{w\cdo\l}=\sum_{y'\in W_\l}p^\l_{y',w'}v^{-|zy'|}T_{zy'}1_\l.\tag a$$
See \cite{\MONO, 1.8(a)}. From (a) we see that
$$p_{y\cdo\l,zw'\cdo\l}=p^\l_{y',w'}(v^2)\text{ if }y=zy',y'\in W_\l,$$
$$p_{y\cdo\l,zw'\cdo\l}=0\text{ if }y\n zW_\l.$$
In particular we have $p_{y\cdo\l,w\cdo\l}\in\NN[v\i]$.
From \cite{\MONO, 1.8} for $w\cdo\l\in I$ we have
$$\wt{c_{w\cdo\l}}=c_{w\cdo\l\i}, c_{w\cdo\l}^\flat=c_{w\i\cdo w(\l)}.$$

\subhead 1.7\endsubhead
Now $\HH$ can be regarded as a two-sided ideal in an $\ca$-algebra $\HH'$ with $1$ as follows.

Let $[\fs_\iy]$ be the set of formal $\ca$-linear combinations $\sum_{\l\in\fs_\iy}c_\l1_\l$ with 
$c_\l\in\ca$; this is an $\ca$-module in an obvious way. We regard $[\fs_\iy]$ as a (commutative) 
$\ca$-algebra with multiplication
$$(\sum_{\l\in\fs_\iy}c_\l1_\l)(\sum_{\l\in\fs_\iy}c'_\l1_\l)=\sum_{\l\in\fs_\iy}c_\l c'_\l1_\l.$$
This algebra has a unit element $1=\sum_{\l\in\fs_\iy}1_\l$. 

Let $\HH'$ be the $\ca$-algebra with generators $T_w (w\in W)$ and $\ph\in[\fs_\iy]$ and relations:

$T_wT_{w'}=T_{ww'}$ if $w,w'\in W$ and $|ww'|=|w|+|w'|$;

$T_\s^2=v^2T_1+(v^2-1)T_\s(\sum_{\l\in\fs_\iy;\s\in W_\l}1_\l)$ for $\s\in S$;

$T_w\ph=\ph'T_w$ for $\ph=\sum_{\l\in\fs_\iy}c_\l1_\l,\ph'=\sum_{\l\in\fs_\iy}c_{w\i(\l)}1_\l$ in 
$[\fs_\iy]$, $w\in W$;

the map $[\fs_\iy]@>>>\HH'$, $\x\m\x$ respects the algebra structures.
\nl
It follows that $\HH'$ is a free left $[\fs_\iy]$-module with basis $\{T_w;w\in W\}$ and a
right free $[\fs_\iy]$-module with basis $\{T_w;w\in W\}$. Note that the algebra $\HH'$ has a unit
element $\sum_{\l\in\fs_\iy}1_\l$.
Now $\HH$ can be identified with the two-sided ideal of $\HH'$ which as an $\ca$-submodule is free
with basis $\{T_w1_\l=1_{w(\l)}T_w;w\cdo\l\in I\}$.

\subhead 1.8\endsubhead
Let $W\bsl\fs_\iy$ be the set of $W$-orbits on $\fs_\iy$.
For any $\fo\in W\bsl\fs_\iy$ we set $I_\fo=\{w\cdo\l\in I;\l\in\fo\}$. This is a finite set.
 We have
$I=\sqc_\fo I_\fo$, $\HH=\op_\fo\HH_\fo$ where $\HH_\fo$ is the $\ca$-submodule of $\HH$ spanned 
by $\{T_w1_\l=1_{w(\l)}T_w;w\cdo\l\in I_\fo\}$ (thus, $H_\fo$ is a free $\ca$-module of finite
rank). If $\fo,\fo'$ are distinct in $W\bsl\fs_\iy$,
then clearly $\HH_\fo\HH_{\fo'}=0$. Thus, each $\HH_\fo$ is a subalgebra of $\HH$; unlike $\HH$, it
has a unit element $\sum_{\l\in\fo}1_\l$. It is stable under $h\m\bar h$ and under $h\m h^\flat$.
Moreover, $h\m\tih$ is an isomorphism of $\HH_\fo$ onto $\HH_{\fo\i}$. For any $w\cdo\l\in I_\fo$ 
we have $c_{w\cdo\l}\in\HH_\fo$; moreover, $\{c_{w\cdo\l};w\cdo\l\in I_\fo\}$ is an $\ca$-basis
of $\HH_\fo$. 

\subhead 1.9\endsubhead
For $i,i'$ in $I$ we write $c_ic_{i'}=\sum_{j\in I}h_{i,i',j}c_j$ (product in $\HH$) where 
$h_{i,i',j}\in\ca$. 
Let $j\underset\text{left}\to\preceq i$ (resp. $j\preceq i$) be the preorder on $I$ generated by 
the relations $h_{i',i,j}\ne0\text{ for some }i'\in I$, resp. by the relations
$$h_{i,i',j}\ne0\text{ or }h_{i',i,j}\ne0\text{ for some }i'\in I.$$
We say that $i\underset\text{left}\to\si j$ (resp. $i\si j$) if 
$i\underset\text{left}\to\preceq j$ and
$j\underset\text{left}\to\preceq i$ (resp. $i\preceq j$ and $j\preceq i$). This is an 
equivalence relation on $I$; the equivalence classes are called left cells (resp. two-sided 
cells). Note that any two-sided cell is a union of left cells. Since for
$\fo\in W\bsl\fs_\iy$, $\HH_\fo$ is closed under left and right multiplication by elements in 
$\HH$, we see that

{\it $h_{i,i',j}\ne0,i\in I_\fo$ implies $i',j\in I_\fo$; $h_{i,i',j}\ne0,i'\in I_\fo$ implies 
$i,j\in I_\fo$.}
\nl
It follows that $j\preceq i, i\in I_\fo$ implies $j\in I_\fo$. In particular, $j\si i, i\in I_\fo$
implies $j\in I_\fo$. Thus any two-sided cell is contained in $I_\fo$ for a unique $\fo$.

For $i=w\cdo\l\in I$ we set 
$$i^!=w\i\cdo w(\l)\in I.$$
Note that $i\m i^!$ is an involution of $I$ preserving $I_\fo$ for any $\fo$.

If $\boc$ is a two-sided cell and $i\in I$, we write $i\preceq\boc$ (resp. $\boc\preceq i$) if 
$i\preceq i'$ (resp. $i'\preceq i$) for some $i'\in\boc$; we write $i\prec\boc$ (resp. 
$\boc\prec i$) if $i\preceq\boc$ (resp. $\boc\preceq i$) and $i\n\boc$. If $\boc,\boc'$ are 
two-sided cells, we write $\boc\preceq\boc'$ (resp. $\boc\prec\boc'$) if $i\preceq i'$ (resp. 
$i\preceq i'$ and $i\not\si i'$) for some $i\in\boc,i'\in\boc'$.

Let $j\in I$. We can find an integer $m\ge0$ such that $h_{i,i',j}\in v^{-m}\ZZ[v]$ for all 
$i,i'$; let $a(j)$ be the smallest such $m$. For $i,i',j$ in $I$ there is a well defined integer 
$h^*_{i,i',j}$ such that
$$h_{i,i',j^!}=h^*_{i,i',j}v^{-a(j^!)}+\text{ higher powers of }v.$$ 
Note that

{\it $h^*_{i,i',j}\ne0,i\in I_\fo$ implies $i',j\in I_\fo$; $h^*_{i,i',j}\ne0,i'\in I_\fo$ implies 
$i,j\in I_\fo$.}
\nl
Let $\DD$ be the set of all $w\cdo\l\in I$ where $w$ is a distinguished involution of the Coxeter 
group $W_\l$, see \cite{\CELLSII}. We have $\DD=\sqc_\fo(\DD\cap\fo)$.

By \cite{\MONO, 1.11}, the following properties hold:

Q1. If $j\in\DD$ and $i,i'\in I$ satisfy $h^*_{i,i',j}\ne0$ then $i'=i^*$.

Q2. If $i\in I$, there exists a unique $j\in\DD$ such that $h^*_{i^!,i,j}\ne0$.

Q3. If $i'\preceq i$ then $a(i')\ge a(i)$. Hence if $i'\si i$ then 
$a(i')=a(i)$.

Q4. If $j\in\DD$, $i\in I$ and $h^*_{i^!,i,j}\ne0$ then $h^*_{i^!,i,j}=1$.

Q5. For any $i,j,k$ in $I$ we have $h^*_{i,j,k}=h^*_{j,k,i}$.

Q6. Let $i,j,k$ in $I$ be such that $h^*_{i,j,k}\ne0$. Then $i\underset\text{left}\to\si j^!$, 
$j\underset\text{left}\to\si k^!$, $k\underset\text{left}\to\si i^!$.

Q7. If $i'\underset\text{left}\to\preceq i$ and $a(i')=a(i)$ then $i'\underset\text{left}\to\si i$.

Q8. If $i'\preceq i$ and $a(i')=a(i)$ then $i'\si i$.

Q9. Any left cell $\G$ of $I$ contains a unique element of $j\in\DD$. We have $h^*_{i^!,i,j}=1$ 
for all $i\in\G$.

Q10. For any $i\in I$ we have $i\si i^!$.

\mpb

Note that $h^*_{i,j,k}\in\NN$ for all $i,j,k$ in $I$, see \cite{\MONO, 1.11}.

\mpb

Let $\HH^\iy$ be the free abelian group with basis $\{t_i;i\in I\}$. We define a $\ZZ$-bilinear 
multiplication $\fA^\iy\T\fA^\iy@>>>\fA^\iy$ by 
$$t_it_{i'}=\sum_{j\in I}h^*_{i,i',j^!}t_j.$$
For any $\fo\in W\bsl\fs_\iy$ let $\HH^\iy_\fo$ be the free abelian subgroup of $\HH^\iy$ with
basis $\{t_i;i\in I_\fo\}$. We have $\HH^\iy=\op_\fo\HH^\iy_\fo$; moreover, if $\fo,\fo'$ are 
distinct in $W\bsl\fs_\iy$, then $\HH^\iy_\fo\HH^\iy_{\fo'}=0$. 
Thus each $\HH^\iy_\fo$ is a subalgebra of $\HH$; unlike $\HH^\iy$, $\HH^\iy_\fo$ has a unit 
element 
$\sum_{i\in\DD\cap\fo}t_i$. The $\ZZ$-linear map $\HH^\iy@>>>\HH^\iy$, $h\m h^\flat$ defined by 
$t_i^\flat=t_{i^!}$ for all $i\in I$ is a ring antiautomorphism preserving each $\HH^\iy_\fo$.
We define an $\ca$-linear map $\psi:\HH@>>>\ca\ot\HH^\iy$ by 
$$\psi(c_i)=\sum_{i'\in I,j\in\DD;i'\si j}h_{i,j,i'}t_{i'}\text{ for all }i\in I.$$
(This last sum is finite. We have $i\in I_\fo$ for some $\fo$. If $h_{i,j,i'}\ne0$ then we have
$i'\in\fo,j\in\fo$. Thus $i',j$ run through a finite set.) By \cite{\MONO, 1.9, 1.11(vi)}, $\psi$ 
is a homomorphism of $\ca$-algebras. For any $\fo$, $\psi$ restricts to a homomorphism of 
$\ca$-algebras $\psi_\fo:\HH_\fo@>>>\ca\ot\HH^\iy_\fo$ which takes $1$ to $1$. 

We set $\HH^v=\QQ(v)\ot_\ca\HH$, $\JJ=\QQ\ot\HH^\iy$; for any $\fo$ we set 
$\HH_\fo^v=\QQ(v)\ot_\ca\HH_\fo$, $\JJ_\fo=\QQ\ot_\ca\HH^\iy_\fo$.
For any $\fo$, $\psi$ induces  an algebra isomorphism 
$\psi^v_\fo:\HH^v_\fo@>\si>>\bbq(v)\ot\JJ_\fo$; 
hence $\psi$ induces an algebra isomorphism $\psi^v:\HH^v@>\si>>\bbq(v)\ot\JJ$.

\mpb

We define a group homomorphism $\tt:\HH^\iy@>>>\ZZ$ by $\tt(t_i)=1$ if $i\in\DD$,
$\tt(t_i)=0$ if $i\in I-\DD$. As in \cite{\MONO, 1.9(a)}, the following can be deduced from
Q1,Q2,Q4.

(a) {\it For $i,j\in I$ we have $\tt(t_it_j)=1$ if $j=i^!$ and $\tt(t_it_j)=0$ if $j\ne i^!$.}

\subhead 1.10\endsubhead
For $n\in\NN^*$ we set $\HH^1_n=\bbq\ot_\ca\HH_n$; this is a $\bbq$-algebra with $1$.
It is the algebra with generators $T_w (w\in W)$, $1_\l (\l\in\fs_n)$ and relations:

$1_\l1_{\l'}=\d_{\l,\l'}1_\l$ for $\l,\l'\in\fs_n$;

$T_wT_{w'}=T_{ww'}$ for $w,w'\in W$;

$T_w1_\l=1_{w(\l)}T_w$ for $w\in W,\l\in\fs_n$;

$T_1=\sum_{\l\in\fs_n}1_\l$.
\nl
It has a basis $\{T_w1_\l;w\cdo\l\in I_n\}$. Let $\HH^1=\bbq\ot_\ca\HH$. This is a $\bbq$-algebra 
without $1$ in general. As a vector space it has basis $\{T_w1_\l,w\cdo\l\in I\}$. It contains 
naturally $\HH_n^1$ as a subalgebra for any $n\in\NN^*$. For any $\fo\in W\bsl\fs_\iy$ we set 
$\HH^1_\fo=\bbq\ot_\ca\HH_\fo$; this is a $\bbq$-algebra with $1$. It has a basis $\{T_w1_\l;w\cdo\l\in I_\fo\}$. We have $\HH^1=\op_\fo\HH^1_\fo$.
Now $\psi$ in 1.9 induces an algebra isomorphism $\psi^1:\HH^1@>\si>>\JJ$; for any $\fo$,
$\psi_\fo$ in 1.9 induces an algebra isomorphism $\psi^1_\fo:\HH^1_\fo@>\si>>\JJ_\fo$ taking $1$
to $1$. 

\subhead 1.11\endsubhead
Let $n\in\NN^*$. Consider the group algebra $\bbq[W\TT_n]$ where $W\TT_n$ is the semidirect 
product of $W$ and $\TT_n$ with $\TT_n$ normal and $W$ acting on $\TT_n$ by $w:t\m w(t)$. Now
$wt\m\sum_{\l\in\fs_n}\l(t)T_w1_\l$ defines a $\bbq$-linear isomorphism 
$u_n:\bbq[W\TT_n]@>\si>>\HH^1_n$ which is in fact an algebra isomorphism taking $1$ to $1$.

Now let $n,n'\in\NN^*$ be such that $n'/n\in\ZZ$.
We define a $\bbq$-linear imbedding $h_{n,n'}:\bbq[W\TT_n]@>>>\bbq[W\TT_{n'}]$ by
$$h_{n,n'}(wt)=(n/n')^\r\sum_{t'\in\TT_{n'};t'{}^{n'/n}=t}wt'.$$
We show that $h_{n,n'}$ is compatible with multiplication, that is, for $w,w'$ in $W$ and $t,t'$ in
$\TT_n$ we have
$$\align&((n/n')^\r\sum_{\tit\in\TT_{n'};\tit{}^{n'/n}=t}w\tit)
((n/n')^\r\sum_{\tit'\in\TT_{n'};\tit'{}^{n'/n}=t'}w'\tit')\\&
=(n/n')^\r\sum_{\tit''\in\TT_{n'};\tit''{}^{n'/n}=w'{}\i(t)t'}ww'\tit'',\endalign$$
or equivalently
$$((n/n')^\r\sum_{\tit,\tit'\in\TT_{n'};\tit{}^{n'/n}=t,\tit'{}^{n'/n}=t'}w'{}\i(\tit)\tit'
\sum_{\tit''\in\TT_{n'};\tit''{}^{n'/n}=w'{}\i(t)t'}\tit'',$$
which is easily verified.

Let $j^1_{n,n'}:\HH^1_n@>\si>>\HH^1_{n'}$ be the specialization of $j_{n,n'}$ (see 1.5) at $v=1$.
We have $u_{n'}h_{n,n'}= j_{n,n'}u_n$; equivalently for $w\in W,t\in\TT_n$, we have
$$(n/n')^\r\sum_{t'\in\TT_{n'};t'{}^{n'/n}=t}\sum_{\l\in\fs_{n'}}\l(t')T_w1_\l
=\sum_{\l\in\fs_n}\l(t)T_w1_\l.$$
(It is enough to show that for any $\l\in\fs_{n'}$,
$$(n/n')^\r\sum_{t'\in\TT_{n'};t'{}^{n'/n}=t}\l(t')=\l(t).$$
is equal to $\l(t)$ if $\l\in\fs_n$ and to $0$ if $\l\n\fs_n$.
This is immediate: we use that the kernel of the surjective homomorphism 
$\TT_{n'}@>>>\TT_n$, $t'\m t'{}^{n'/n}$ has exactly $(n'/n)^\r$ elements.)

We can form the union $\cup_{n\in\NN^*}\bbq[W\TT_n]$ over all imbeddings $h_{n,n'}$ as above. 
This union has an algebra structure whose restriction to $\bbq[W\TT_n]$ (for any $n\in\NN^*$) is
the algebra structure of $\bbq[W\TT_n]$. Moreover, there is a unique isomorphism of
algebras $\cup_{n\in\NN^*}\bbq[W\TT_n]@>\si>>\HH^1$ whose restriction to $\bbq[W\TT_n]$ (for any 
$n\in\NN^*$) is $u_n:\bbq[W\TT_n]@>\si>>\HH^1_n$.

\subhead 1.12\endsubhead
For $\fo\in W\bsl \fs_\iy$, $\HH^1_\fo$ is a semisimple $\bbq$-algebra. Let $\Irr(H^1_\fo)$ be a 
set of representatives for the isomorphism classes of simple $\HH^1_\fo$-modules. 

\subhead 1.13\endsubhead
We have $\HH^\iy=\op_{\boc}\HH^\iy_\boc$, $\JJ=\op_{\boc}\JJ_\boc$, where $\boc$ runs over the 
two-sided cells in $I$, $\HH^\iy_\boc$ is the $\ca$-submodule of $\HH^\iy$ with basis 
$\{t_i;i\in\boc\}$ and $\JJ_\boc$ is the $\bbq$-subspace of $\JJ$ with basis $\{t_i;i\in\boc\}$.  
Each $\HH^\iy_\boc$ is an $\ca$-subalgebra of $\HH^\iy$ with unit $\sum_{i\in\DD_\boc}t_i$ where
$\DD_\boc=\DD\cap\boc$. Each $\JJ_\boc$ is a $\bbq$-subalgebra of $\JJ$ with the same unit as
$\HH^\iy_\boc$. Moreover if $\boc,\boc'$ are distinct two-sided cells in $I$ we have
$\JJ_\boc\JJ_{\boc'}=0$.
Recall from 1.9 that any two-sided cell in $I$ is contained in $I_\fo$ for a unique 
$\fo\in W\bsl\fs_\iy$. It follows that for any $\fo\in W\bsl\fs_\iy$ we have
$\JJ_\fo=\op_{\boc\sub I_\fo}\JJ_\boc$.
Hence, if $E\in\Irr(H^1_\fo)$ then there is a unique two-sided cell $\boc_E$ such that
$\JJ_{\boc}$ acts as zero on $E^\iy$ for any $\boc\sub I_\fo$ with $\boc\ne\boc_E$. Thus $E^\iy$
can be viewed as a simple $\JJ_{\boc_E}$-module. We define $a_E\in\NN$ to be the constant
value of the restriction of $a:I@>>>\NN$ to $\boc_E$.

\subhead 1.14\endsubhead
If $\boc$ is a two-sided cell of $I$ then its image $\wt\boc$ under $I@>>>I$, 
$w\cdo\l\m w\cdo\l\i$ is a two-sided cell of $I$. (See \cite{\MONO, 1.14}.) As noted in 1.9, we
have $\boc\sub I_\fo$ for a unique $\fo$; from the definitions we have $\wt\boc\sub I_{\fo\i}$. 
Moreover, the value of the $a$-function on $\wt\boc$ is equal to the value of the $a$-function on 
$\boc$. From Q3,Q10 in 1.9, we see that $a(i^!)=a(i)$ for $i\in I$.

\subhead 1.15\endsubhead
For $i,i'$ in $I$ we show:

(a) {\it If $i\underset\text{left}\to\si i'$, then for some $u\in I$, $t_{i'}$ appears  with 
$\ne0$ coefficient in $t_ut_i$.}

(b) {\it If $i^!\underset\text{left}\to\si i'{}^!$, then for some $u\in I$, $t_{i'}$ appears with 
$\ne0$ coefficient in $t_it_u$.}

(c) {\it If $i\si i'$, then for some $u,u'$ in $I$, $t_{i'}$ appears with nonzero coefficient in 
$t_ut_it_u'$.}

(d) {\it If $i\si i'$, then $t_it_jt_{i'}\ne0$ for some $j\in I$.}
\nl
The proof is along the lines of that of \cite{\HEC, 18.4}. Let $J^+=\sum_{k\in I}\NN t_k$. We will
use repeatedly that $J^+J^+\sub J^+$. 

Let $i,i'$ be as in (a). Let $d,d'\in\DD$ be such that $h^*_{i^!,i,d}\ne0$ and
$h^*_{i'{}^!,i',d'}\ne0$. Then $i\underset\text{left}\to\si d$, $i'\underset\text{left}\to\si d'$.
Hence $d\underset\text{left}\to\si d'$. By Q9 in 1.9 we have $d=d'$ and $h^*_{i^!,i,d}=1$,
$h^*_{i'{}^!,i',d}=1$. Hence $t_{i^!}t_i=t_d+J^+$, $t_{i'{}^!}t_{i'}=t_d+J^+$, $t_dt_d=t_d$; it 
follows that $t_{i^!}t_it_{i'{}^!}t_{i'}\in t_dt_d+J^+=t_d+J^+$. In particular,
$t_it_{i'{}^!}\ne0$. Thus, $h^*_{i,i'{}^!,u}\ne0$ for some $u\in I$. Using Q5 in 1.9 we deduce
that $h^*_{u,i,i'{}^*}\ne0$ hence $t_{i'}$ appears with $\ne0$ coefficient in $t_ut_i$. This proves
(a). Now (b) follows from (a) using the antiautomorphism of $\HH^\iy$ such that $t_u\m t_{u^!}$
for all $u\in I$.

Let $i_1,i_2,i_3$ in $I$ be such that $i_1\si i_2\si i_3$. If the conclusion of (c) holds for
$(i,i')=(i_1,i_2)$ and for $(i,i')=(i_2,i_3)$ then clearly it holds for $(i,i')=(i_1,i_3)$. 
Applying this repeatedly, we see that it is enough to prove (c) in the case where $i,i'$ satisfy 
either $i\underset\text{left}\to\si i'$ or 
$i^!\underset\text{left}\to\si i'{}^!$. In these cases the desired result follws from (a),(b).

Let $i,i'$ be as in (d). Then $i\si i'{}^!$. By (c), we have
$t_{u'}t_it_u\in at_{i'{}^!}+J^+$ for some $u,u'\in I$ and some $a\in\ZZ_{>0}$. Hence 
$t_{u'}t_it_ut_{i'}\in at_{i'{}^!}t_{i'}+J^+$. Since $t_{i'{}^!}t_{i'}$ has some coefficient $1$ 
and the other coefficienrs are $\ge0$, it follows that $t_{u'}t_it_ut_{i'}\ne0$. Thus, 
$t_it_ut_{i'}\ne0$. This proves (d).

\head 2. The group $\tG$\endhead
\subhead 2.1\endsubhead
In this paper (except in 2.2) we fix a group $\tG$ containing $G$ as a subgroup, 
such that $\tG/G$ is cyclic of order $\mm\le\iy$ with a fixed generator. For $s\in\ZZ$ let 
$\tG_s$ be the inverse image of the $s$-th power of this generator under the obvious map 
$\tG@>>>\tG/G$. For $\g\in\tG$, the map $G@>>>G$, $g\m\g g\g\i$ is denoted by $\Ad(\g)$. 

We shall always assume that we are in one of the two cases below (later referred to as case A
and case B).

(A) We have $\mm=\iy$ and one of the following two equivalent conditions
are satisfied ($q$ denotes a fixed power of $p$):

(i) for some $\g\in\tG_1$, $\Ad(\g):G@>>>G$ is the Frobenius map for an $F_q$-rational 
rational structure on $G$;

(ii) for any $s>0$ and any $\g\in\tG_s$, $\Ad(\g):G@>>>G$ is the Frobenius map for an 
$F_{q^s}$-rational rational structure on $G$.
\nl
(B) $\mm<\iy$ and $\tG$ is an algebraic group with identity component $G$.

\mpb

We show the equivalence of (i),(ii) in case A. Clearly, if (ii) holds then (i) holds. 
Conversely, assume that (i) holds for $\g\in\tG_1$.
If $\g'\in\tG_s$ with $s>0$, then we have $\g'=g_1\g^s$ where $g_1\in G$. By Lang's theorem applied
to $\Ad(\g^s):G@>>>G$, which is the Frobenius map for an $F_{q^s}$-rational structure on $G$, we
have $g_1=g_2\i\Ad(\g^s)(g_2)$ for some $g_2\in G$ hence 
$\g'=g_2\i\Ad(\g^s)(g_2)\g^s=g_2\i\g^sg_2$ and $\Ad(\g')=\Ad(g_2)\i\Ad(\g^s)\Ad(g_2)$. Since 
$\Ad(g_2):G@>>>G$ is an isomorphism of algebraic varieties, it follows that
$\Ad(\g'):G@>>>G$ is the Frobenius map for an $F_{q^s}$-rational structure on $G$. Thus (ii)
holds.

Let $s\in\ZZ$. In case B, $\tG_s$ is naturally an algebraic variety.
In case A, we view $\tG_s$ as an algebraic variety using the bijection $g\m g\g$ where
$\g$ is fixed in $\tG_s$; this algebraic structure on $\tG_s$ is independent of the choice 
of $\g$. For $s=0$ this gives the usual structure of algebraic variety of $G$. 
For $s\in\ZZ,s'\in\ZZ$, the multiplication $\tG_s\T\tG_{s'}@>>>\tG_{s+s'}$ is obviously a
morphism of algebraic varieties in case B, but is only a quasi-morphism in the sense of 
\cite{\URE, 0.3} in case A. Similarly, for $s\in\ZZ$, $\tG_s@>>>\tG_{-s}$, $\g\m\g\i$ is a 
morphism of algebraic varieties in case B, but is only a quasi-morphism in case A.

\mpb

Note that in case A with $s\ne0$, the conjugation action of $G$ on $\tG_s$ is transitive. 
(If $s>0$, this follows from as above using Lang's theorem, while if $s<0$ this follows 
using the bijection $\tG_s@>>>\tG_{-s}$, $\g\m\g\i$, which commutes with the $G$-actions.) 
Moreover in this case for any $\g\in\tG_s$, the stabilizer of $\g$ for this $G$-action is 
finite. (This stabilizer is the fixed point set of $\Ad(\g):G@>>>G$ which is a Frobenius map
relative to an $F_{q^s}$-structure if $s>0$ or the inverse of a Frobenius map if $s<0$.) 

We show:

(a) {\it If $\g\in\tG_s$ and $B\in\cb$ then $\Ad(\g)(B)\in\cb$, $\Ad(\g)(U_B)=U_{\Ad(\g)B}$ and 
$\Ad(\g):\cb@>>>\cb$ is a bijection.}
\nl
In case A with $s=0$ and in case B, (a) is obvious. In case A with $s>0$, (a) follows from 
(ii); in case A with $s<0$, (a) follows from (ii) applied to $\g\i$. 

\subhead 2.2\endsubhead
Here are some examples in case A.

(i) Let $F:G@>>>G$ be the Frobenius map for an $F_q$-rational structure on $G$.
Let $\tG=G\T\ZZ$ regarded as a group with multiplication $(g,s)(g',s')=(gF^s(g'),s+s')$. Define 
a homomorphism $\tG@>>>\ZZ$ by $(g,s)\m s$. Its kernel $\{(g,s)\in\tG;s=0\}$ can be identified 
with $G$. Note that $\tG$ and $\tG@>>>\ZZ$ are as in case A; we have $(1,1)\in\tG_1$ and
$\Ad(1,1):G@>>>G$ is just $F:G@>>>G$. Moreover, any $\tG$ and $\tG@>>>\ZZ$ as in case A is 
obtained by the procedure above. 

\mpb

(ii) In the case where $G$ is adjoint we define $\tG_s$ for $s\in\ZZ_{<0}$ to be the set of
Frobenius maps $G@>>>G$ with respect to various split $F_{q^s}$-rational structures on $G$; we
define $\tG_s$ for $s\in\ZZ_{<0}$ to be the set of maps $G@>>>G$ whose inverse is in $\tG_{-s}$
and we set $\tG_0=G$. Then $\tG=\sqc_{s\in\ZZ}\tG_s$ is as in case A. (This case has been 
considered
in \cite{\URE}.)

(iii) Let $V$ be a finite dimensional $\kk$-vector space. For any $s\in\ZZ$ let $\wt{GL(V)}_s$ be 
the set of all group isomorphisms $T:V@>>>V$ such that $T(zx)=z^{q^s}T(x)$ for all 
$z\in\kk,x\in V$; in particular we have $\wt{GL(V)}_0=GL(V)$. Then 
$\wt{GL(V)}:=\sqc_{s\in\ZZ}\wt{GL(V)}_s$ is a group under composition of maps; it is of the form 
$\tG$ (as in case A) where $G=GL(V)$.

\mpb

(iv) Let $V$ be a finite dimensional $\kk$-vector space with a nondegenerate symplectic form
$(,):V\T V@>>>\kk$. For any $s\in\ZZ$ let $\wt{Sp(V)}_s$ be the set of all $T\in\wt{GL(V)}_s$ such
that $(T(x),T(x'))=(x,x')^{q^s}$ for all $x,x'$ in $V$; in particular we have 
$\wt{Sp(V)}_0=Sp(V)$. Then $\wt{Sp(V)}:=\sqc_{s\in\ZZ}\wt{Sp(V)}_s$ is a group under composition 
of maps; it is of the form $\tG$ (as in case A) where $G=Sp(V)$.

\subhead 2.3\endsubhead
{\it In the rest of this paper we fix $\t\in\tG_1$ such that $\t\BB\t\i=\BB,\t\TT\t\i=\TT$.
and such that for any $\s\in S$, $\Ad(\t)$ carries $\x_\s\in\UU_\s-\{1\}$ to
$\x_{\s'}\in\UU_{\s'}-\{1\}$ for some $\s'\in S$.}
\nl
Note that such $\t$ exists.

We define a group homomorphism $\ee:\tG@>>>\tG$ by $\ee(\g)=\t\g\t\i$. We have
$\ee(\tG_s)=\tG_s$ for all $s\in\ZZ$, $\ee(\TT)=\TT$, $\ee(\BB)=\BB$ (hence $\ee(\UU)=\UU$),
$\ee(N\TT)=N\TT$; thus $\ee$ induces an automorphism of $W$ denoted again by $\ee$ which
preserves the Coxeter group structure.
If $B\in\cb$ then $\ee(B)\in\cb$ and $B\m\ee(B)$, $\cb@>>>\cb$ is an automorphism in case B
and is the Frobenius map for an $\FF_q$-rational structure on $\cb$ in case A. We define 
$\ee:\cb\T\cb@>>>\cb\T\cb$ by $\ee(B,B')=(\ee(B),\ee(B'))$. For $w\in W$ we have 
$\ee(G_w)=G_{\ee(w)}$ and $\ee(\co_w)=\co_{\ee(w)}$.

The set $\{\dot\s;\s\in S\}$ of $N\TT$ is stable under $\ee:N\TT@>>>N\TT$. For $w\in W$ we 
have $(\ee(w))\dot{}=\ee(\dw)$. Hence $N_0\TT$ is stable under $\ee:N\TT@>>>N\TT$.

\mpb

Now for $n\in\NN^*$, $\ee:\TT@>>>\TT$ restricts to an isomorphism $\ee:\TT_n@>>>\TT_n$ and this
induces an isomorphism $\ee:\fs_n@>>>\fs_n$ by $\l\m\ee(\l)$ where $(\ee(\l))(t)=\l(\ee\i(t))$ for
$t\in\TT_n$. Let $\ee:\fs_\iy@>>>\fs_\iy$ be the isomorphism whose restriction to $\fs_n$ is 
$\ee:\fs_n@>>>\fs_n$ as above for any $n\in\NN^*$.

\mpb

We shall fix a Frobenius map $\Ps:G@>>>G$ relative to some sufficiently large finite subfield
$F_{q'}$ of $\kk$ such that $\BB,\TT$ are $\Ps$-stable, $\Ps$ acts on $t$ by $t\m t^{q'}$ (hence it acts as the identity on $W$) and such that $\Ps\ee=\ee\Ps:G@>>>G$ and $\Ps(\o)=\o$ for any 
$\o\in N_0\TT$; in case B we also require that $\Ps(\t^\mm)=\t^\mm$.

For any $s\in\ZZ$ we define an $F_{q'}$-rational structure on $\tG_s$ with Frobenius map
$\Ps:\tG_s@>>>\tG_s$ by the requirement that $\Ps(g\t^s)=\Ps(g)\t^s$ for any $g\in G$;
in case B, this rational structure depends only on $\tG_s$ not on $s$.

Now for any $n\in\NN^*$ we have $\Ps(\TT_n)=\TT_n$;  hence we can define 
$\Ps:\fs_n@>\si>>\fs_n$ by
$(\Ps\l)(t)=\l(\Ps\i(t))$ for $t\in\TT_n$, $\l\in\fs_n$. There is a unique bijection 
$\Ps:\fs_\iy@>>>\fs_\iy$ whose restriction to $\fs_n$ is as above for any $n\in\NN^*$.
Now $\Ps$ induces $F_{q'}$-rational structures on various varieties that will appear in the sequel.
 When we consider $\cd_m()$ or $\cm_m()$ for such varieties, we will refer to these specific 
$F_{q'}$-structures. 

\subhead 2.4\endsubhead
We define a bijection $\ee:I@>>>I$ by $\ee(w\cdo\l)=\ee(w)\cdo\ee(\l)$. The $\ca$-linear map 
$\ee:\HH@>>>\HH$ defined by $\ee(T_w1_\l)=T_{\ee(w)}1_{\ee(\l)}$ for $w\cdo\l\in I$ is an algebra 
isomorphism commuting with $\bar{}:\HH@>>>\HH$. It follows that $\ee(c_i)=c_{\ee(i)}$ 
for all $i\in I$ and that $\ee:I@>>>I$ maps any left (resp. two-sided) cell of $I$ 
onto a left (resp. two-sided) cell of $I$. It also maps any $W$-orbit in $\fs_\iy$ onto a
$W$-orbit in $\fs_\iy$.

Let $\fo\in\fs_\iy$ and $s\in\ZZ$ be such that $\ee^s(\fo)=\fo$.
The $\ca$-linear map $\ee^s:\HH@>>>\HH$ restricts to an $\ca$-algebra isomorphism 
$\ee^s:\HH_\fo@>>>\HH_\fo$; this gives rise by extension of scalars 
to a $\bbq$-algebra isomorphism $\ee^s:\HH^1_\fo@>>>\HH^1_\fo$ and to a $\bbq(v)$-algebra 
isomorphism $\ee:\HH^v_\fo@>>>\HH^v_\fo$; moreover the $\bbq$-linear map 
$\ee^s:\JJ_\fo@>>>\JJ_\fo$ given by $t_i\m t_{\ee^s(i)}$ for $i\in I_\fo$ is an algebra
isomorphism and $\psi^v_\fo:\HH^v_\fo@>\si>>\bbq(v)\ot\JJ_\fo$, 
$\psi^1_\fo:\HH^1_\fo@>\si>>\JJ_\fo$ are compatible with the action of $\ee^s$.

Let $\Irr_s(\HH^1_\fo)$ be the set of all $E\in\Irr(\HH^1_\fo)$ with the following property:
there exists a linear isomorphism $\ee_s:E@>>>E$ such that for any $w\cdo\l\in I_\fo$ and 
any $e\in E$ we have 
$$\ee_s((T_w1_\l)(e)))=(T_{\ee^s(w)}1_{\ee^s(\l)})(\ee_s(e)).$$
(Such $\ee_s$ is clearly unique up to a nonzero scalar, if it exists.) We assume that for any 
$E\in\Irr_s(\HH^1_\fo)$, an $\ee_s$ as above has been chosen; we can assume that $\ee_s$ has 
finite order (since $\ee^s:I_\fo@>>>I_\fo$ has finite order); moreover, when $s=0$ we have 
$\Irr_s(\HH^1_\fo)=\Irr(\HH^1_\fo)$ and for any $E$ in this set we can take $\ee_s=1$. If 
$E\in\Irr(H^1_\fo)$ we can view $E$ as a simple $\JJ_\fo$-module via $\psi^1_\fo$; we denote this 
$\JJ_\fo$-module by $E^\iy$. Moreover we can view $\bbq(v)\ot E^\iy$ as a simple 
$\HH^v_\fo$-module via $\psi^v_\fo$; we denote this $\HH^v_\fo$-module by $E^v$. If in addition we
have $E\in\Irr_s(H^1_\fo)$, then $\ee_s$ can be viewed as a $\bbq$-linear isomorphism 
$E^\iy@>>>E^\iy$ (denoted again by $\ee_s$) and as a $\bbq(v)$-linear isomorphism $E^v@>>>E^v$ 
(denoted again by $\ee_s$).

Note that for any $\x\in\JJ_\fo$, $e\in E^\iy$  we have $\ee_s(\x(e))=\ee^s(\x)(\ee_s(e))$; for 
any $\x'\in\HH_\fo$, $e'\in E^v$ we have $\ee_s(\x'(e'))=\ee^s(\x')(\ee_s(e'))$. 

\subhead 2.5\endsubhead
For $s\in\ZZ$ let 
$$I^s=\{w\cdo\l\in I;w(\l)=\ee^{-s}(\l)\}.$$
For any two-sided cell $\boc$ of $I$ we set
$$\boc^s=I^s\cap\boc.$$
We show:

(a) {\it If $\ee^s(\boc)=\boc$ and $i\in\boc$, $j\in I$ satisfy $t_{i^!}t_jt_{\ee^s(i)}\ne0$, then
$j\in\boc^s$.}

(b) {\it If $\ee^s(\boc)=\boc$, then $\boc^s\ne\emp$.}
\nl
We prove (a). Let $i=w\cdo\l$, $j=z\cdo\l'$. From our assumption  we have \lb
$t_{z\cdo\l'}t_{\ee^s(w)\cdo\ee^s(\l)}\ne0$ (which implies $\l'=\ee^s(w(\l))$) and
$t_{w\i\cdo w(\l)}t_{z\cdo\l'}\ne0$ (which implies $w(\l)=z(\l')$). We deduce that
$z(\l')=\ee^{-s}(\l')$ so that $j\in I^s$. 
Since $t_{i^!}t_j\ne0$ and $i^!\in\boc$ we must have 
$j\in\boc$. Thus we have $j\in I^s\cap\boc$ and (a) is proved.

We prove (b). Let $i\in\boc$. By assumption we have $\ee^s(i)\in\boc$; by Q10 in 1.9 we have
$i^!\in\boc$. Using 1.15(d) with $i,i'$ replaced by $i^!,\ee^s(i)$ we see that for some
$j=z\cdo\l'\in I$  we have $t_{i^!}t_jt_{\ee^s(i)}\ne0$. Using (a) we deduce that $j\in\boc^s$
and (b) is proved.

\head 3. Sheaves on $\tcb^2$\endhead
\subhead 3.1\endsubhead
Let $\tcb=G/\UU$. We have $\tcb^2=\sqc_{w\in W}\tco_w$ where
$$\tco_w=\{(x\UU,y\UU)\in\tcb^2;x\i y\in G_w\}.$$ 
The closure of $\tco_w$ in $\tcb^2$ is $\btco_w=\cup_{y\in W;y\le w}\tco_y$.
For $w\in W$ and $\o\in\k_0\i(w)$ we define $G_w@>>>\TT$ by $g\m g_\o$ where $g\in\UU\o g_\o\UU$, 
$g_\o\in\TT$. We define $j^\o:\tco_w@>>>\TT$ by $j^\o(x\UU,y\UU)=(x\i y)_\o$. For $\l\in\fs_\iy$ 
we set $L_\l^\o=(j^\o)^*L_\l$, a local system on $\tco_w$. Let $L_\l^{\o\sha}$ be its
extension to an intersection cohomology complex on $\btco_w$ viewed as a complex on 
$\tcb^2$, equal to $0$ on $\tcb^2-\btco_w$. We shall view $L_\l^\o$ as a 
constructible sheaf on $\tcb^2$ which is $0$ on $\tcb^2-\tco_w$. 
Let $\LL_\l^\o=L_\l^{\o\sha}\la|w|+\nu+2\r\ra$, a simple perverse sheaf on $\tcb^2$. 

(a) {\it In the remainder of this section we fix a two-sided cell $\boc$ of $I$ and we set 
$a=a(i)$ for some/any $i\in\boc$. We define $\fo\in W\bsl\fs_\iy$ by $\boc\sub I_\fo$. We denote 
by $n$ the smallest integer in $\NN^*$ such that $\fo\sub\fs_n$. We shall assume that $\Ps$ in
2.3 acts as $1$ on the finite subset $\{t\in\TT;t^n\in\TT\cap N_0\TT\}$ of $\TT$.}
\nl
In particular, $\Ps(t)=t$ for any $t\in\TT_n$ (hence $\Ps(\l)=\l$ for any $\l\in\fs_n$). 

Now, if $w\in W,\o\in\k_0\i(w),\l\in\fs_n$, then $L_\l^\o|_{\tco_w}$, $L_\l^{\o\sha}$ and 
$\LL_\l^\o$ can be regarded naturally as objects in the mixed derived category of pure weight 
zero. Moreover, $L_\l^\o|_{\tco_w}$ (resp. $L_\l^{\o\sha}$, $\LL_\l^\o$) is (noncanonically) 
isomorphic to $L_\l^{\dw}|_{\tco_w}$ (resp. $L_\l^{\dw\sha}$, $\LL_\l^{\dw}$) in the mixed derived
category. (It is enough to show that if $t,t'\in\TT, t^n=t'=\dw\o\i$ and $h_{t'}:\TT@>>>\TT$ is 
translation by $t'$, then $t$ defines an isomorphism $h_{t'}^*L_\l@>>>L_\l$; see 
\cite{\MONO, 1.5}.)

\mpb

We define $\ti\fh:\tcb^2@>>>\tcb^2$ by $(x\UU,y\UU)\m(y\UU,x\UU)$.

We define an action of $G\T\TT^2$ on $\tcb^2$ (resp. on $\TT$) by
$$(g,t_1,t_2):(x\UU,y\UU)\m(gxt_1^n\UU,gyt_2^n\UU)$$
(resp. by $(g,t_1,t_2):t\m w\i(t_1)^{-n}tt_2^n$). For any $w\in W$, the $G\T\TT^2$-action leaves 
stable $\tco_w$ and its restriction to $\tco_w$ is transitive; moreover, $j^\o$ is compatible with
actions of $G\T\TT^2$ on $\tco_w$ and $\TT$. 

If $\l\in\fs_n$ then $L_\l$ is a $G\T\TT^2$-equivariant local system on $\TT$ hence $L_w^\l$ is a 
$G\T\TT^2$-equivariant local system on $\tco_w$. By \cite{\MONO, 2.1}, the following holds.

(c) {\it For fixed $w\in W,\o\in\k_0\i(w)$, the local systems $L_\l^\o$ with $\l\in\fs_n$ form a 
set of representatives for the isomorphism classes of irreducible $G\T\TT^2$-equivariant local 
systems on $\tco_w$.} 

\subhead 3.2\endsubhead
We define $p_{01}:\tcb^3@>>>\tcb^2$, $p_{12}:\tcb^3@>>>\tcb^2$, $p_{02}:\tcb^3@>>>\tcb^2$ by
$$\align&p_{01}(x\UU,y\UU,z\UU)=(x\UU,y\UU), p_{12}(x\UU,y\UU,z\UU)=(y\UU,z\UU),\\&
p_{02}(x\UU,y\UU,z\UU)=(x\UU,z\UU).\endalign$$ 
For any $L\in\cd(\tcb^2)$, $L'\in\cd(\tcb^2)$, we set
$$L\cir L'=p_{02!}(p_{01}^*L\ot p_{12}^*L')\in\cd(\tcb^2).$$
This defines a monoidal structure on $\cd(\tcb^2)$. Thus, if ${}^iL\in\cd(\tcb)$ for $i=1,\do,k$, then
${}^1L\cir{}^2L\cir\do\cir{}^kL\in\cd(\tcb)$ is well defined. Note that, if 
$L\in\cd_m(\tcb^2)$, $L'\in\cd_m(\tcb^2)$ then $L\cir L'$ is naturally in $\cd_m(\tcb^2)$.

\subhead 3.3\endsubhead
Now assume that $w,w'\in W$, $\o\in\k\i_0(w),\o'\in\k\i_0(w')$, $\l,\l'\in\fs_\iy$. From
\cite{\MONO, 2.3} we see that:

(a) {\it if $w'(\l')\ne\l$, then $L_\l^\o\cir L_{\l'}^{\o'}=0$.}

\subhead 3.4\endsubhead
Now assume that $w,w'\in W$, $\o\in\k_0\i(w),\o'\in\k_0\i(w')$, $\l,\l'\in\fs_\iy$.  
Let $\Xi$ be the set of all $(x\UU,y\UU,z\UU)\in\tcb^3$ such that $x\i y\in\UU\o t\UU$, 
$y\i z\in\UU\o't'\UU$ for some $t,t'$ in $\TT$ (which are in fact uniquely determined). 
Define $c:\Xi@>>>\TT\T\TT$ by $c(x\UU,y\UU,z\UU)=(t,t')$ where $x\i y\in\UU\o t\UU$,
$y\i z\in\UU\o't'\UU$. Define $p'_{02}:\Xi@>>>\tcb^2$ by $(x\UU,y\UU,z\UU)\m(x\UU,z\UU)$.
From the definitions we see that
$$L_\l^\o\cir L_{\l'}^{\o'}=p'_{02!}(c^*(L_\l\bxt L_{\l'})).\tag a$$
We show: 

(b) {\it If $w'(\l')=\l$ and $|ww'|=|w|+|w'|$, then we have canonically 
$L_\l^\o\cir L_{\l'}^{\o'}=L_{\l'}^{\o\o'}\ot\fL$, with $\fL$ as in 0.2.}
\nl
Let $Y=\{(x\UU,z\UU,t,t')\in\tcb\T\tcb\T\TT\T\TT;x\i z\in\UU\o t\UU\o't'\UU\}$. We define 
$\Xi@>>>Y$ by $(x\UU,y\UU,z\UU)\m(x\UU,z\UU,t,t')$ where $t,t'$ in $\TT$ are given by 
$x\i y\in\UU\o t\UU$, $y\i z\in\UU\o't'\UU$. This is an isomorphism since $|ww'|=|w|+|w'|$. We 
identify $\Xi=Y$ through this isomorphism. Then $c:\Xi@>>>\TT\T\TT$ becomes $c:Y@>>>\TT\T\TT$, 
$(x\UU,z\UU,t,t')\m(t,t')$. We define $h:\TT\T\TT@>>>\TT$ by $h(t,t')=w'{}\i(t)t'$. We have 
$$Y=\{(x\UU,z\UU,t,t')\in\tcb\T\tcb\T\TT\T\TT;x\i z\in\UU\o\o'h(t,t')\UU\}.$$
Define $j:Y@>>>\tco_{ww'}$ by $(x\UU,z\UU,t,t')\m(x\UU,z\UU)$. 
Let $j'=j^{\o\o'}:\tco_{ww'}@>>>\TT$. Using (a) and the cartesian diagram
$$\CD
Y@>c>>\TT\T\TT  \\
@VjVV                  @VhVV  \\
\tco_{ww'}@>j'>>\TT\\
\endCD$$
we see that
$$L_\l^\o\cir L_{\l'}^{\o'}=j_!c^*(L_\l\bxt L_{\l'})=j'{}^*h_!(L_\l\bxt L_{\l'})$$.
Since $L_{\l'}^{\o\o'}\ot\fL=j'{}^*(L_{\l'}\ot\fL)$, we see that to prove (b) it is enough to show
that $h_!(L_\l\bxt L_{\l'})=L_{\l'}\ot\fL$ (assuming that $w'(\l')=\l$). This is proved as in the 
last paragraph of \cite{\MONO, 2.4}.

\subhead 3.5\endsubhead
Let $\s\in S$ and let $\o\in\k_0\i(\s)$, $\l'\in\fs_\iy$. Define $\d_\o:\UU_\s-\{1\}@>>>\TT$ by 
$\x\m t_\x\i$ where $t_\x\in\TT$ is given by $\o\i\x\i\o\in\UU\o\i t_\x\UU$; let 
$\ce=\d_\o^*L_{\l'}^*$. Let $\d':\UU_\s-\{1\}@>>>\pp$ be the obvious map. From the definitions we 
see that:

(a) {\it $\d'_!\ce=0$ if $\s\n W_{\l'}$; $\d'_!\ce\Bpq\{\bbq\la-2\ra,\bbq[-1]\}$ if 
$\s\in W_{\l'}$.}
\nl
Consider the diagram $\TT@<\tk<<\TT\T(\UU_\s-\{1\})@>\tih>>\TT$ where $\tk:(t,\x)\m t_\x\i$ and 
$\tih:(t,\x)\m tt_\x\i$. We show:

(b) {\it Let $\l'\in\fs_\iy$. If $\s\n W_{\l'}$, then $\tih_!\tk^*L_{\l'}=0$.  If 
$\s\in W_{\l'}$ then $\tih_!\tk^*L_{\l'}^*\Bpq\{\bbq\la-2\ra,\bbq[-1]\}$.}
\nl
We have $\tk^*L_{\l'}^*=\bbq\bxt\ce$. Now $\ti h=\ti h'y$ where
$y:\TT\T(\UU_\s-\{1\})@>>>\TT\T(\UU_\s-\{1\})$ is $(t,\x)\m(tt_\x\i,\x)$ and 
$\tih':\TT\T(\UU_\s-\{1\})@>>>\TT$ is $(t,\x)\m t$. Clearly, $y_!(\bbq\bxt\ce)=\bbq\bxt\ce$. It 
remains to note that $\tih_!(\bbq\bxt\ce)$ is $0$ if $\s\n W_{\l'}$ and is 
$\Bpq\{\bbq\la-2\ra,\bbq[-1]\}$ if $\s\in W_\l$. (This follows from (a).)

\mpb

We show:

(c) {\it Assume that $\l\in\fs_\iy$ satisfies $\s\in W_\l$ and that $\o\in\{\dot\s,\dot\s\i\}$. 
Then we have canonically $L_\l^\o=L_\l^{\o\i}$.}
\nl
Define $\z:\TT@>>>\TT$ by $t\m\o^2t$. It is enough to show that $\z^*L_\l=L_\l$ canonically. For 
$t\in\TT$ we have $(\z^*L_\l)_t=(L_\l)_{\o^2 t}=(L_\l)_{\cha_\s(-1)}\ot(L_\l)_t$. Hence it is 
enough to show that we have canonically $(L_\l)_{\cha_\s(-1)}=\bbq$. It is also enough to show 
that $\cha_\s^*L_\l=\bbq$. This follows from $\a_\s\in R_\l$.

\subhead 3.6\endsubhead
Now assume that $w=w'=\s\in S$, $\o\in\k_0\i(\s)$, $\l,\l'\in\fs_\iy$ are such that $\s(\l')=\l$. 
In this subsection we show:

(a)  {\it If $\s\n W_\l$, then $L_\l^\o\cir L_{\l'}^{\o\i}=L_{\l'}^1\la-2\ra\ot\fL$.}

(b) {\it If $\s\in W_\l$, then}
$$L_\l^\o\cir L_{\l'}^{\o\i}\Bpq\{L_{\l'}^1\la-2\ra\ot\fL,
L_{\l'}^\o\la-2\ra\ot\fL,L_{\l'}^\o[-1]\ot\fL\}.$$
\nl
(Note that the conditions $\s\in W_\l$ and $\s\in W_{\l'}$ are equivalent.) With the notation of 
3.4, we have
$$\Xi=\{(x\UU,y\UU,z\UU)\in\tcb^3;x\i y\in\UU\o t\UU,y\i z\in\UU\o\i t'\UU\text{ for some $t,t'$ 
in }\TT\}.$$
If $(x\UU,y\UU,z\UU)\in\Xi$ then $x\i z\in\UU\o\UU\o\i w'{}\i(t)t'\UU$; in particular we 
have $x\i z\in\BB\cup\BB\o\BB$. Thus, $\Xi$ can be partitioned as $\tcb^I\cup\tcb^{II}$ where
$$\tcb^I=\{(x\UU,y\UU,z\UU)\in\Xi;x\i z\in\BB\}$$ is a closed subset and
$$\tcb^{II}=\{(x\UU,y\UU,z\UU)\in\Xi;x\i z\in\BB\o\BB\}$$ is an open subset. The map 
$p'_{02}:\Xi@>>>\tcb^2$ (see 3.4) restricts to maps
$$p^I_{02}:\tcb^I@>>>\tco_1, p^{II}_{02}:\tcb^{II}@>>>\tco_\s;$$ using 3.4(a) we deduce
$$L_\l^\o\cir L_{\l'}^{\o\i}\Bpq\{p^I_{02!}(c^*(L_\l\bxt L_{\l'})),\qua
p^{II}_{02!}(c^*(L_\l\bxt L_{\l'}))\}.$$
We show:
$$p^I_{02!}(c^*(L_\l\bxt L_{\l'}))=L_{\l'}^1\ot\fL\la-2\ra.\tag c$$
We have
$$\align&\tcb^I=\{(x\UU,y\UU,z\UU)\in\tcb^3;x\i y\in\UU\o t\UU,y\i z\in\UU\o\i t'\UU\\&
\text{for some $t,t'$ in }\TT,x\i z\in\BB\},\endalign$$
or equivalently 
$$\tcb^I=\{(x\UU,y\UU,z\UU)\in\tcb^3;x\i y\in\UU\o t\UU,x\i z\in\UU\s(t)t'\UU
\text{for some $t,t'$ in }\TT\}.$$ 
Let $Y=\{(x\UU,z\UU,t,t')\in\tcb\T\tcb\T\TT\T\TT;x\i z\in\UU\s(t)t'\UU\}$.
We define $d:\tcb^I@>>>Y$ by $(x\UU,y\UU,z\UU)\m(x\UU,z\UU,t,t')$ where $t,t'$ in $\TT$ are as in
the last formula for $\tcb^I$. The fibre of $d$ at $(x\UU,z\UU,t,t')\in Y$ can be identified with 
$\{y\UU;y\in x\UU\o t\UU\}$, an affine line. Thus, $d$ is an affine line bundle. We have a 
cartesian diagram
$$\CD
Y@>c^I>>\TT\T\TT  \\
@Vj^I VV       @VhVV  \\
\tco_1@>\tj^I>>\TT\\
\endCD$$
where $c^I:Y@>>>\TT\T\TT$ is $(x\UU,z\UU,t,t')\m(t,t')$, $j^I:Y@>>>\tco_1$ is
$(x\UU,z\UU,t,t')\m(x\UU,z\UU)$, $\tj^I=j^1:\tco_1@>>>\TT$, $h:\TT\T\TT@>>>\TT$ is
$(t,t')\m\s(t)t'$.
As in 3.4 we have $h_!(L_\l\bxt L_{\l'})=L_{\l'}\ot\fL$ (since $\s(\l')=\l$). It follows that
$$(j^I)!(c^I)^*(L_\l\bxt L_{\l'})=(\tj^I)^* h_!(L_\l\bxt L_{\l'})=(\tj^I)^* L_{\l'}\ot\fL.$$
Hence 
$$\align&p^I_{02!}(c^*(L_\l\bxt L_{\l'}))=(j^I)!d_!d^*(c^I)^*(L_\l\bxt L_{\l'})
=(j^I)!(c^I)^*(L_\l\bxt L_{\l'})\la-2\ra\\&=
(\tj^I)^* L_{\l'}\ot\fL\la-2\ra=L_{\l'}^1\ot\fL\la-2\ra.\endalign$$
This proves (c). Next we show that
$$p^{II}_{02!}(c^*(L_\l\bxt L_{\l'}))\text{ is $0$ if $\s\n W_{\l'}$ and is 
$\Bpq\{L_{\l'}^\o\la-2\ra,L_{\l'}^\o[-1]\}$ if }\s\in W_{\l'}.\tag d$$
We have
$$\align&\tcb^{II}=\{(x\UU,y\UU,z\UU)\in\tcb^3;x\i y\in\UU\o t\UU,y\i z\in\UU\o\i t'\UU\\&
\text{ for some $t,t'$ in }\TT,x\i z\in\UU\o t_1\UU\text{ for some }t_1\in\TT\}.\endalign$$
Let $(x\UU,z\UU)\in\tco_\s$. We can write uniquely $z=x\x_0\o t_1u_1$ where $\x_0\in\UU_\s$, 
$t_1\in\TT$, $u_1\in\UU$. The fibre $\Ph$ of $p^{II}_{02}$ at $(x\UU,z\UU)$ can be identified with 
$$\align&\{y\UU\in G/UU;x\i y\in\UU\o t\UU,y\i z\in\UU\o\i t'\UU\}\\&
=\{y\UU\in G/UU;x\i y\in\UU\o t\UU,y\i x\x_0\o t_1u_1\in\UU\o\i t'\UU\}.\endalign$$
Setting $x\i y=\x\o tu'$ where $\x\in\UU_\s$, we can identify
$$\align&\Ph=\{(t,t',\x)\in\TT\T\TT\T\UU_\s;u'{}\i t\i\o\i\x\i\x_0\o t_1\in\UU\o\i t'\UU\}\\&
=\{(t,t',\x)\in\TT\T\TT\T\UU_\s;\o\i\x\i\x_0\o\in\UU\o\i\s(t)t't_1\i\UU\}\\&
=\{(t,t',\x)\in\TT\T\TT\T(\UU_\s-\{\x_0\});t_{\x\i\x_0}=\s(t)t't_1\i\}\endalign$$
where for $\x_1\in\UU_s-\{1\}$ we define $t_{\x_1}\in\TT$ by $\o\i\x_1\i\o\in\UU\o\i t_{\x_1}\UU$.
Let 
$$Y'=\{(x\UU,z\UU,t,t',\x_1)\in\tcb\T\tcb\T\TT\T\TT\T(\UU_\s-\{1\});
x\i z\in\UU_\s\o\s(t)t't_{\x_1}\i\UU\},$$
$$Y'_1=\{(x\UU,z\UU,t'_1,\x_1)\in\tcb\T\tcb\T\TT\T(\UU_\s-\{1\});
x\i z\in\UU_\s\o t'_1t_{\x_1}\i\UU\}.$$
We see that $\tcb^{II}$ may be identified with $Y'$. (The identification is via 
$$(x\UU,z\UU,t,t',\x_1)\m(x\UU,x\x_0\x_1\i\o t\UU,z\UU)$$
 where 
$\x_0\in\UU_\s$ is given by $x\i z\in\x_0\o\TT\UU$.) Under this identification, $p^{II}_{02}$ 
becomes the composition $fj^{II}$ where $j^{II}:Y'@>>>Y'_1$ is 
$$(x\UU,z\UU,t,t',\x_1)\m(x\UU,z\UU,s(t)t',\x_1)$$
 and $f:Y'_1@>>>\tco_\s$ is
$$(x\UU,z\UU,t'_1,\x_1)\m(x\UU,z\UU);$$
 moreover, the local system $c^*(L_\l\bxt L_{\l'})$ on 
$\tcb^{II}$ becomes the local system \lb
$(c^{II})^*(L_\l\bxt L_{\l'})$ on $Y'$ where 
$c^{II}:Y'@>>>\TT\T\TT$ is $(x\UU,z\UU,t,t',\x_1)\m(t,t')$.
We have a diagram with cartesian squares
$$\CD
{}   @.              Y'@>c^{II}>>\TT\T\TT\\
 @.                @Vj^{II}VV    @VhVV  \\
\TT\T(\UU_\s-\{1\})@<\tj'<<Y'_1 @>\tj^{II}>>\TT\\
  @V\tih VV                @VfVV              @.\\
\TT@<j'<<              \tco_s @.  {}
\endCD$$
where $\tj^{II}:Y'_1@>>>\TT$ is $(x\UU,z\UU,t'_1,\x_1)\m t'_1$, $j':\tco_\s@>>>\TT$ is $j^\o$,
$\tj':Y'_1@>>>\TT\T(\UU_\s-\{1\})$ is $(x\UU,z\UU,t'_1,\x_1)\m(t'_1,\x_1)$, $h:\TT\T\TT@>>>\TT$ 
is $(t,t')\m\s(t)t'$ and $\tih'$ is as in 3.5.

Let $L'=(\tj^{II})^*L_{\l'}$ (a local system on $Y'_1$). Let $L''=j'{}^*L_{\l'}=L_{\l'}^\o$ (a 
local system on $\tco_\s$).
Define $\tf:Y'_1@>>>\TT$ by $(x\UU,z\UU,t'_1,\x_1)\m t_{\x_1}\i$. Let $\tL=\tf^*L_{\l'}$ (a local 
system on $Y'_1$). The stalk of $L'$ at $(x\UU,z\UU,t'_1,\x_1)\in Y'_1$ is $(L_{\l'})_{t'_1}$.
The stalk of $f^*L''$ at $(x\UU,z\UU,t'_1,\x_1)\in Y'_1$ is 
$(L_{\l'})_{t'_1t_{\x_1}\i}=(L_{\l'})_{t'_1}\ot(L_{\l'})_{t_{\x_1}\i}$. Thus we have 
$L'=f^*L''\ot\tL^*$. 

As in 3.4 we have $h_!(L_\l\bxt L_{\l'})=L_{\l'}\ot\fL$ (since $\s(\l')=\l$). Using the cartesian 
diagrams above, we see that
$$\align&p^{II}_{02!}(c^*(L_\l\bxt L_{\l'}))=f_!j^{II}_!(c^{II})^*(L_\l\bxt L_{\l'})\\&
=f_!j^{II}_!(c^{II})^*(L_\l\bxt L_{\l'})=
f_!(\tj^{II})^*h_!(L_\l\bxt L_{\l'})=f_!(\tj^{II})^*L_{\l'}\ot\fL)\\&
=f_!(L')\ot\fL=f_!(f^*L''\ot\tL^*)\ot\fL=L''\ot f_!(\tL^*)\ot\fL
=L''\ot f_!\tj'{}^*\tk^*(L_{\l'}^*)\\&
=L''\ot f_!\tj'{}^*\tk^*(L_{\l'}^*)=L''\ot j'{}^*\tih_!\tk^*(L_{\l'}^*)=
L''\ot j'{}^*\tih_!\tk^*(L_{\l'}^*).\endalign$$
Here $\tk$ is as in 3.5. Using 3.5(b) we see that this is $0$ if $\s\n W_{\l'}$ and is 
$\Bpq\{L''\la-2\ra,L''[-1]\}$ if $\s\in W_{\l'}$.
This proves (d). Now (a),(b) follow from (c),(d).

\subhead 3.7\endsubhead
Now assume that $w\in W$, $\s\in S$, $\o\in\{\dot\s,\dot\s\i\},\o'\in\k\i_0(w)$, 
$\l,\l'\in\fs_\iy$ are such that $w(\l')=\l$, $|\s w|<|w|$. We show:

(a) {\it If $\s\n W_\l$ then 
$L_\l^\o\cir L_{\l'}^{\o'}\ot\fL=L_{\l'}^{\o\o'}\la-2\ra\ot\fL\ot\fL$.}

(b) {\it If $\s\in W_\l$, then}
$$L_\l^\o\cir L_{\l'}^{\o'}\ot\fL\Bpq\{L_{\l'}^{\o\o'}\la-2\ra\ot\fL\ot\fL,
L_{\l'}^{\o'}\la-2\ra\ot\fL\ot\fL,L_{\l'}^{\o'}[-1]\ot\fL\ot\fL\}.$$
Using 3.4(b), we have $L_{\l'}^{\o'}\ot\fL=L_{(\s w)(\l')}^{\o\i}\cir L_{\l'}^{\o\o'}$. Hence 
$L_\l^\o\cir L_{\l'}^{\o'}\ot\fL=L_\l^\o\cir L_{(\s w)(\l')}^{\o\i}\cir L_{\l'}^{\o\o'}$.
We now apply 3.6(a),(b) to describe $L_\l^\o\cir L_{(\s w)(\l')}^{\o\i}$. If $\s\n W_\l$, we obtain
$$L_\l^\o\cir L_{\l'}^{\o'}\ot\fL=L_{(\s w)(\l')}^1\cir L_{\l'}^{\o\o'}\la-2\ra\ot\fL.$$
By 3.4(b) this equals $L_{\l'}^{\o\o'}\la-2\ra\ot\fl^{\ot2}$, proving (a). If $\s\in W_\l$, we 
obtain
$$\align&L_\l^\o\cir L_{\l'}^{\o'}\ot\fL   
\Bpq\{L_{(\s w)\l'}^1\cir L_{\l'}^{\o\o'}\la-2\ra\ot\fL,\\&
L_{(\s w)\l'}^{\o\i}\cir L_{\l'}^{\o\o'}\la-2\ra\ot\fL,L_{(\s w)\l'}^{\o\i}\cir L_{\l'}^{\o\o'}[-1]
\ot\fL\}.\endalign$$
(We have used that $L_{(\s w)\l'}^\o=L_{(\s w)\l'}^{\o\i}$, see 3.5(c).) We now substitute

$L_{(\s w)\l'}^1\cir L_{\l'}^{\o\o'}=L_{\l'}^{\o\o'}\ot\fL$, 
$L_{(\s w)\l'}^{\o\i}\cir L_{\l'}^{\o\o'}=L_{\l'}^{\o'}\ot\fL$, 
\nl
see 3.4(b); we obtain
$$L_\l^\o\cir L_{\l'}^{\o'}\ot\fL\Bpq \{L_{\l'}^{\o\o'}\la-2\ra\ot\fL\ot\fL,
L_{\l'}^{\o'}\la-2\ra\ot\fL\ot\fL,L_{\l'}^{\o'}[-1]\ot\fL\ot\fL\}.$$
This proves (b).

\subhead 3.8\endsubhead
Let $\cd^\spa\tcb^2$ be the subcategory of $\cd(\tcb^2)$ consisting of objects
which are restrictions of objects in the $G\T\TT^2$-equivariant derived category.
Let $\cm^\spa\tcb^2$ be the subcategory of $\cd^\spa\tcb^2$ consisting of objects which are 
perverse sheaves. Let $\cm^{\preceq}\tcb^2$  (resp. $\cm^{\prec}\tcb^2$) be the subcategory
of $\cm^\spa \tcb^2$ whose objects are perverse sheaves $L$ such that any 
composition factor of $L$ is of the form $\LL_\l^{\dw}$ for some 
$w\cdo\l\preceq\boc$ (resp. $w\cdo\l\prec\boc$). Let $\cd^{\preceq}\tcb^2$ (resp. 
$\cd^{\prec}\tcb^2$) be the subcategory of $\cd^\spa\tcb^2$ whose objects are complexes
$L$ such that $L^j$ is in $\cm^{\preceq}\tcb^2$ (resp. $\cm^{\prec}\tcb^2$) for any $j$.
We write $\cd_m()$ or $\cm_m()$ for the mixed version of any of the categories above.
Let $\cc^\spa\tcb^2$ be the subcategory of $\cm^\spa\tcb^2$ consisting of semisimple objects. 
Let $\cc_0^\spa \tcb^2$ be the subcategory of 
$\cm^\spa_m\tcb^2$ consisting of objects of pure of weight zero. 
Let $\cc^\boc \tcb^2$ be the subcategory of $\cm^\spa \tcb^2$ consisting of objects 
which are direct sums of objects of the form $\LL_\l^{\dw}$ with 
$w\cdo\l\in\boc$. Let $\cc_0^\boc\tcb^2$ be the subcategory of $\cc^\spa_0\tcb^2$ 
consisting of those $L\in\cc^\spa_0\tcb^2$ such that, as an object of $\cc^\spa\tcb^2$, $L$ 
belongs to $\cc^\boc\tcb^2$. For $L\in\cc_0^\spa\tcb^2$ let $\un{L}$ be the 
largest subobject of $L$ such that as an object of $\cc^\spa\tcb^2$, we have 
$\un{L}\in\cc^\boc\tcb^2$.

\subhead 3.9 \endsubhead
Let $r\ge1$. Let $\ww=(w_1,\do,w_r)\in W^r$, $\pmb\o=(\o_1,\o_2,\do,\o_r)$ be such that 
$\o_i\in\k\i_0(w_i)$ for $i=1,\do,r$ and $\pmb\l=(\l_1,\l_2,\do,\l_r)\in\fs_n^r$. We set
$$|\ww|=|w_1|+|w_2|+\do+|w_r|.$$
For $J\sub[1,r]$, let
$$\align&\tco^J_\ww=\{(x_0\UU,x_1\UU,\do,x_r\UU)\in\tcb^{r+1};\\&x_{i-1}\i x_i\UU\in
\bG_{w_i}\frl i\in J,x_{i-1}\i x_i\in G_{w_i}\frl i\in[1,r]-J\}.\endalign$$
Define $c:\tco_\ww^\emp@>>>\TT^r$ by
$$c(x_0\UU,x_1\UU,\do,x_r\UU)=((x_0\i x_1)_{\o_1},(x_1\i x_2)_{\o_2},\do,
(x_{r-1}\i x_r)_{\o_r}).$$
Let $M^{\pmb\o}_{\pmb\l}\in\cd_m(\tcb^{r+1})$ be the local system 
$c^*(L_{\l_1}\bxt\do\bxt L_{\l_r})$ on $\tco_\ww^\emp$ extended by $0$ on 
$\tcb^{r+1}-\tco_\ww^\emp$.  For $J\sub[1,r]$ we set 
$$M^{\pmb\o,J}_{\pmb\l}=p_{01}^*{}^1L\ot p_{12}^*{}^2L\ot\do
\ot p_{r-1,r}^*{}^rL\in\cd_m(\tcb^{r+1}),$$
$$L^{\pmb\o,J}_{\pmb\l}=p_{0r!}M^{\pmb\o,J}_{\pmb\l}\la|\ww|\ra
={}^1L\cir{}^2L\cir\do\cir{}^rL\la|\ww|\ra\in\cd_m(\tcb^2),$$
where ${}^iL$ is $L_{\l_i}^{\o_i\sha}$ for $i\in J$ and $L_{\l_i}^{\o_i}$ for 
$i\n J$. Note that $M^{\pmb\o,\emp}_{\pmb\l}=M^{\pmb\o}_{\pmb\l}$. Moreover, from
\cite{\MONO, 2.15} we have:

(a) {\it $M^{\pmb\o,J}_{\pmb\l}$ is the intersection cohomology complex of 
$\tco_\ww^J$ with coefficients in $M^{\pmb\o}_{\pmb\l}$.}

\mpb

Consider the free $\TT^{r-1}$-action on $\tcb^{r+1}$ given by
$$\align&(\t_1,\t_2,\do,\t_{r-1}):(x_0\UU,x_1\UU,\do,x_{r-1}\UU,x_r\UU)\m\\&
(x_0\UU,x_1\t_1\UU,\do,x_{r-1}\t_{r-1}\UU,x_r\UU).\endalign$$
Note that $\tco^J_\ww$ is stable under this $\TT^{r-1}$-action. We also have a free 
$\TT^{r-1}$-action on $\TT^r$ given by
$$\align&(\t_1,\t_2,\do,\t_{r-1}):(t_1,t_2,\do,t_r)\m\\&
(t_1\t_1,w_2\i(\t_1\i)t_2\t_2,w_3\i(\t_2\i)t_3\t_3,\do,
w_{r-1}\i(\t_{r-2}\i)t_{r-1}\t_{r-1},w_r\i(\t_{r-1}\i)t_r).\endalign$$
Let ${}'\tcb^{r+1}=\TT^{r-1}\bsl\tcb^{r+1}$. Let ${}'\tco^J_\ww=\TT^{r-1}\bsl\tco^J_\ww$ (a 
locally closed subvariety of ${}'\tcb^{r+1}$). Let ${}'\TT^r=\TT^{r-1}\bsl\TT^r$. Note that 
${}'\tco^\emp_\ww=\TT^{r-1}\bsl\tco^\emp_\ww$ is an open dense smooth 
irreducible subvariety of ${}'\tco^J_\ww$. Now $c:\tco_\ww^\emp@>>>\TT^r$ is compatible with the 
$\TT^{r-1}$-actions on $\tco_\ww^\emp,\TT^r$ hence it induces a map 
${}'c:{}'\tco_\ww^\emp@>>>{}'\TT^r$. The homomorphism $c':\TT^r@>>>\TT$ given by
$$(t_1,t_2,\do,t_r)\m t_1w_2(t_2)w_2w_3(t_3)\do w_2w_3\do w_r(t_r)$$
is constant on each orbit of the $\TT^{r-1}$-action on $\TT^r$ hence it induces a morphism 
${}'\TT^r@>>>\TT$ whose composition with ${}'c$ is denoted by $\bc:{}'\tco_\ww^\emp@>>>\TT$. Let 
${}'M^{\pmb\o,\emp}_{\pmb\l}$ be the local system $\bc^*L_{\l_1}$ on ${}'\tco_\ww^\emp$ extended 
by $0$ on ${}'\tcb^{r+1}-{}'\tco_\ww^\emp$. Let ${}'M^{\pmb\o,J}_{\pmb\l}\in\cd_m({}'\tcb^{r+1})$ 
be the intersection cohomology complex of ${}'\tco_\ww^J$ with coefficients in 
${}'M^{\pmb\o,\emp}_{\pmb\l}$ extended by $0$ on ${}'\tcb^{r+1}-{}'\tco_\ww^J$. Let 
$\bp_{0r}:{}'\tco_\ww^J@>>>\tcb^2$ be the map induced by $p_{0r}:\tco_\ww^J@>>>\tcb^2$. We define 
${}'L^{\pmb\o,J}_{\pmb\l}\in\cd_m^\spa\tcb^2$ as follows: 

if $\l_k=w_{k+1}(\l_{k+1})\text{ for }k=1,2,\do,r-1$, we set 
${}'L^{\pmb\o,J}_{\pmb\l}=\bp_{0r!}{}'M^{\pmb\o,J}_{\pmb\l}\la|\ww|\ra$;

otherwise, we set ${}'L^{\pmb\o,J}_{\pmb\l}=0$.

\subhead 3.10\endsubhead  
For $L,L'\in\cc^\boc_0\tcb^2$ we set
$$L\un\cir L'=\un{(L\cir L')^{\{a-\nu\}}}\in\cc^\boc_0\tcb^2.$$
(For the notation ${}^{\{i\}}$ see 0.2.)
By \cite{\MONO, 2.24}, $L,L'\m L\un{\cir}L'$ defines a monoidal structure on 
$\cc^\boc_0\tcb^2$. Hence if ${}^1L,{}^2,\do,{}^rL$ are in $\cc^\boc_0\tcb^2$ 
then ${}^1L\un{\cir}{}^2L\un{\cir}\do\un{\cir}{}^rL\in\cc^\boc_0\tcb^2$ is well defined.

\subhead 3.11\endsubhead  
Let $w\cdo\l\in I_n$ and let $\o\in\k\i(w),s\in\ZZ$. We show that we have canonically:

(a) $(\ee^s)^*L_\l^\o=L_{\ee^{-s}(\l)}^{\ee^{-s}(\o)}$,
    $(\ee^s)^*\LL_\l^\o=\LL_{\ee^{-s}(\l)}^{\ee^{-s}(\o)}$.
\nl
It is enough to prove the first of these equalities. 
Let $\x=(x\UU,y\UU)\in\tcb^2$. We have $x\i y\in\UU\ee^{-s}(\o)t\UU$ with $t\in\TT$ hence
$\ee^s(x)\i\ee^s(y)\in\UU\o\ee^s(t)\UU$. The stalk of $(\ee^s)^*L_\l^\o$ at $\x$ is equal to the 
stalk of $L_\l$ at $\ee^s(t)$ hence to the stalk of $(\ee^s)^*L_\l$ at $t$. The stalk of 
$L_{\ee^{-s}(\l)}^{\ee^{-s}(\o)}$ at $\x$ is equal to the stalk of $L_{\ee^{-s}(\l)}$ at $t$.
It remains to show that $(\ee^s)^*L_\l=L_{\ee^{-s}(\l)}$. This follows from the definitions.

\head 4. Sheaves on $Z_s$\endhead
\subhead 4.1\endsubhead
{\it In this section we fix $s\in\ZZ$.}
\nl
Now $\TT$ acts on $\tcb^2$ by $t:(x\UU,y\UU)\m(xt\UU,y\ee^s(t)\UU)$. Let $\TT\bsl_s\tcb^2$ be the 
set of orbits. Let
$$Z_s=\{(B,B',\g U_B);B\in\cb,B'\in\cb,\g U_B\in\tG_s/U_B;\g B\g\i=B'\}.$$
We define $\e_s:\tcb^2@>>>Z_s$ by $\e_s:(x\UU,y\UU)\m(x\BB x\i,y\BB y\i,y\t^s\UU x\i)$. Clearly, 
$\e_s$ induces a map $\TT\bsl_s\tcb^2@>>>Z_s$. We show:

(a) {\it $\e_s$ induces an isomorphism $\TT\bsl_s\tcb^2@>>>Z_s$.}
\nl
We show only that our map is bijective.
Let $(B,B',\g)\in\cb\T\cb\T\tG_s$ be such that $\g B\g\i=B'$. We can find $x\in G$ such that 
$B=x\BB x\i$. We set $y=\g x\t^{-s}\in G$. Then $\e_s$ carries the $\TT$-orbit of $(x\UU,y\UU)$ to 
$(B,\g B\g\i,\g x\UU x\i)=(B,B',\g U_B)$; thus our map is surjective. Now assume that
$x,x',y,y'$ in $G$ are such that 

$(x\BB x\i,y\BB y\i,y\t^s\UU x\i)=(x'\BB x'{}\i,y'\BB y'{}\i,y'\t^s\UU x'{}\i)$.
\nl
To complete the proof of (a) it is enough to show that $x'=xtu$, $y'=y\ee^s(t)u'$ for some $u,u'$ 
in $\UU$ and some $t\in\TT$. Since $x\i x'\in\BB$ we have $x'=xtu$ for some $u\in\UU$ and some 
$t\in\TT$. We have $y'\t^s\UU u\i t\i x\i=y\t^s\UU x\i$ hence $y'=y\ee^s(t)u'$ for some 
$u'\in\UU$. This completes the proof of (a).

\mpb

For $w\in W$ let $Z_s^w=\{(B,B',\g U_B)\in Z_s;(B,B')\in\co_w\}$. 
The closure of $Z_s^w$ in $Z_s$ is $\bZ_s^w=\{(B,B',gU_B);(B,B')\in\bco_w,g\in G,gBg\i=B'\}$.
We have $\e_s\i(Z_s^w)=\tco_w$, $\e_s\i(\bZ_s^w)=\btco_w$.

Let $\o\in\k_0\i(w)$ and let $\l\in\fs_\iy$ be such that $w\cdo\l\in I^s$. We have a diagram 
$\TT@<j^\o<<\tcb^2_w@>\e_s^w>>Z_s^w$ where $\e_s^w$ is the restriction of $\e_s$ and $j^\o$ is as 
in 3.1. The $\TT$-action on $\tcb^2$ described above is compatible under $j^\o$ with the 
$\TT$-action on $\TT$ given by $t:t'\m w\i(t\i)t'\ee^s(t)$. From \cite{\CDGVI, 28.2} we see that 
$L_\l$ is equivariant for the $\TT$-action on $\TT$ given by $t:t'\m w\i(\ee^{-s}(t_1))t't_1\i$.
(We use that $w\cdo\l\in I^s$.) Using the change of variable $t_1=\ee^s(t)\i$, we deduce that 
$L_\l$ is also equivariant for the $\TT$-action on $\TT$ 
given by $t:t'\m w\i(t\i)t'\ee^s(t)$. It follows that $(j^\o)^*L_\l$ is $\TT$-equivariant, so that
there is a well defined local system $\cl_{\l,s}^\o$ of rank $1$ on $Z_s^w$ such that 
$(\e_s^w)^*\cl_{\l,s}^\o=(j^\o)^*L_\l=L_\l^\o$. Let 
$\cl_{\l,s}^{\o\sha}$ be its extension to an intersection cohomology complex of $\bZ_s^w$, viewed 
as a complex on $Z_s$, equal to $0$ on $Z_s-\bZ_s^w$. 
We shall view $\cl_{\l,s}^\o$ as a constructible sheaf on $Z_s$ which is $0$ on $Z_s-Z_s^w$. Let
$$\Bbb L_{\l,s}^\o=\cl_{\l,s}^{\o\sha}\la|w|+\nu+\r\ra,$$
a simple perverse sheaf on $Z_s$. 

{\it In the remainder of this subsection we assume that $s\ne0$ and that we are in case A.}
\nl
Let $w\in W$ and let $X_s^w=\{B\in\cb;(B,\ee^s(B))\in\co_w\}$. When $s>0$, $X_s^w$ coincides with
the variety $X_w$ defined in \cite{\DL} in terms of the Frobenius map $\ee^s:G@>>>G$; when $s<0$, 
$X_s^w$ can be identified with the variety $X_{\ee^{-s}(w\i)}$ defined in \cite{\DL} in terms of 
the Frobenius map $\ee^{-s}:G@>>>G$. Note that the finite group $G^{\ee^s}=\{g\in G;\ee^s(g)=g\}$ 
acts by conjugation on $X_s^w$.

Let $\tX_s^w=\{x\UU\in G/\UU;x\i\ee^s(x)\in G_w\}$. We define $\ph:\tX_s^w@>>>X_s^w$ by
$x\UU\m x\BB x\i$. This is a principal $\TT$-bundle with $\TT$ acting on $\tX_s^w$ by
$t:x\UU\m xt\UU$. We define $j'_{\dw}:\tX_s^w@>>>\TT$ by $j'_{\dw}(x\UU)=(x\i\ee^s(x))_{\dw}$. Now
let $\l\in\fs_\iy$ be such that $w\cdo\l\in I^s$. Then there is a well defined local system 
$\cf_{\l,s}^{\dw}$ on $X_s^w$ such that $\ph^*\cf_{\l,s}^{\dw}=(j'_{\dw})^*L_\l$. (This is in fact
the restriction of $\cl_{\l,s}^{\dw}$ to $X_s^w$ under the imbedding $X_s^w@>>>Z_s^w$, 
$x\BB x\i\m(x\BB x\i,\ee^s(x)\BB\ee^s(x\i),\t^sx\UU x\i)$.) 
The local system $\cf_{\l,s}^{\dw}$ on $X_s^w$ is of the type considered in \cite{\DL}.
Note also that $\cf_{\l,s}^{\dw}$ has a natural $G^{\ee^s}$-equivariant structure. (It is the 
restriction of the $G$-equivariant structure of $\cl_{\l,s}^{\dw}$.) It follows that for $j\in\ZZ$,
$H^j_c(X_s^w,\cf_{\l,s}^{\dw})$ is naturally a $G^{\ee^s}$-module. (This representation of 
$G^{\ee^s}$ is one of the main themes of \cite{\DL}.) Let 
$\bX_s^w=\{B\in\cb;(B,\ee^s(B))\in\bco_w\}$. Then $X_s^w$ is open dense smooth in $\bX_s^w$ and
$G^{\ee^s}$ acts by conjugation on $\bX_s^w$. Hence for $j\in\ZZ$, the intersection cohomology 
space $IH^j(\bX_s^w,\cf_{\l,s}^{\dw})$ is naturally a $G^{\ee^s}$-module. 

If $\rr,\rr'$ are $G^{\ee^s}$-modules and $\rr$ is irreducible we denote by $(\rr:\rr')$ the 
multiplicity of $\rr$ in $\rr'$. Let $\Irr(G^{\ee^s})$ be the set of isomorphism classes of 
irreducible representations of $G^{\ee^s}$. From \cite{\DL, 7.7} it is known that for any 
$\rr\in\Irr(G^{\ee^s})$

(i) there exists $w\cdo\l\in I^s$ such that $(\rr:\op_jH^j_c(X_s^w,\cf_{\l,s}^{\dw}))\ne0$.
\nl
From \cite{\ORA, 2.8} we see using (i) that for any $\rr\in\Irr(G^{\ee^s})$

(ii) there exists $w\cdo\l\in I^s$ such that $(\rr:\op_jIH^j(X_s^w,\cf_{\l,s}^{\dw}))\ne0$.
\nl
By \cite{\DL, 6.3}, any $\rr\in\Irr(G^{\ee^s})$ determines a $W$-orbit $\fo$ on $\fs_\iy$: the set 
of all $\l\in\fs_\iy$ such that $(\rr:\op_jH^j_c(X_s^w,\cf_{\l,s}^{\dw}))\ne0$ for some $w\in W$ 
with $w\cdo\l\in I^s$ or equivalently (see \cite{\ORA, 2.8}) such that 
$(r:\op_jIH^j(\bX_s^w,\cf_{\l,s}^{\dw}))\ne0$ for some $w\in W$ with $w\cdo\l\in I^s$; we have 
necessarily $\ee^s(\fo)=\fo$. For any $\fo\in W\bsl\fs_\iy$ such that $\ee^s(\fo)=\fo$, let 
$\Irr_{\fo}(G^{\ee^s})$ be the set of all $\rr\in\Irr(G^{\ee^s})$ such that the $W$-orbit on 
$\fs_\iy$ determined by $\rr$ is $\fo$. With notation in 1.14 we have the following result:

(b) {\it There exists a pairing $\Irr_{\fo}(G^{\ee^s})\T\Irr_s(\HH^1_\fo)@>>>\bbq$, 
$(\rr,E)\m b_{\rr,E}$ such that for any $\rr\in\Irr_{\fo}(G^{\ee^s})$, any 
$z\cdo\l\in I^s\cap I_\fo$ and any $j\in\ZZ$ we have}
$$(\rr:IH^j(\bX_s^z,\cf_{\l,s}^{\dz}))=(-1)^j(j-|z|:\sum_{E\in\Irr_s(\HH^1_\fo}b_{\rr,E}
\tr(\ee_sc_{z\cdo\l},E^v).$$
In the case where $G$ has connected centre, (b) is just a reformulation on \cite{\ORA, 3.8(ii)}. A
proof similar to that in {\it loc.cit.} applies without the hypotesis on the centre.

\subhead 4.2\endsubhead
{\it In the remainder of this section let $\boc,a,\fo,n,\Ps$ be as in 3.1(a).}
\nl
The $G\T\TT^2$-action on $\tcb^2$ defined in 3.1 commutes with the $\TT$-action on $\tcb^2$ in 
4.1; hence it induces a $G\T\TT^2$-action on $\TT\bsl_s\tcb^2$. We define a $G\T\TT^2$-action on 
$Z_s$ by 
$$(g,t_1,t_2):(B,B',\g U_B)\m(gBg\i,gB'g\i,g\g x_0\ee^s(t_2^{-n})t_1^nx_0\i g\i U_{gBg\i})$$
where $x_0$ is any element of $G$ such that $x_0\BB x_0\i=B$. (The choice of $x_0$ does not 
matter; to see this, it is enough to show that for $b\in B$ we have 
$$\g x_0\ee^s(t_2^{-n})t_1^nx_0\i U_B=\g x_0b\ee^s(t_2^{-n})t_1^nb\i x_0\i U_B$$
which is immediate.) In this $G\T\TT^2$ action, the subgroup 
$\{(1,t_1,t_2)\in G\T\TT^2;t_1=\ee^s(t_2)\}$ acts trivially. Note that the bijection 
$\TT\bsl_s\tcb^2@>>>Z_s$ in 4.1(a) is compatible with the $G\T\TT^2$-actions. 

Let $w\in W,\o\in\k_0\i(w)$. Since the $G\T\TT^2$-action on $\tco_w$ is 
transitive, it follows that the $G\T\TT^2$-action on $Z_s^w$ is transitive. We show :

(a) {\it Let $\cl$ be an irreducible $G\T\TT^2$-equivariant local system on $Z_s^w$. Then $\cl$ is
isomorphic to $\cl_{\l,s}^\o$ for a unique $\l\in\fs_n$ such that $w\cdo\l\in I^s$.}
\nl
The local system $(\e_s^w)^*\cl$ on $\tco_w$ is irreducible and $G\T\TT^2$-equivariant hence, by 
3.1(c), is isomorphic to $L_\l^\o$ for a well defined $\l\in\fs_n$. Now the restriction of 
$(\e_s^w)^*\cl$ to any fibre of $\e_s^w$ is $\bbq$. On the other hand, the restriction of 
$L_\l^\o$ to the fibre of $\e_s^w$ passing through $(\UU,\o\UU)$ is (under an obvious 
identification with $\TT$) the inverse image of $L_\l$ under the map $\TT@>>>\TT$, 
$t\m w\i(t\i)\ee^s(t)$, hence it is $L_{w(\l\i)\ee^{-s}(\l)}$ which is $\bbq$ if and only if 
$w(\l)=\ee^{-s}\l$. We see that we must have $w(\l)=\ee^{-s}(\l)$. We have 
$(\e_s^w)^*\cl\cong(\e_s^w)^*\cl_{\l,s}^\o$ (both are $L_\l^\o$) hence $\cl\cong\cl_{\l,s}^\o$. 
This proves (a).

\mpb

We define $\fh:Z_s@>>>Z_{-s}$ by $(B,B',gU_B)\m(B',B,g\i U_{B'})$.  Note that 
$\fh\e_s=\e_{-s}\ti\fh:\tcb^2@>>>Z_{-s}$ with $\ti\fh$ as in 3.1. For $L\in\cd_m(Z_{-s})$ we set 
$L^\da=\fh^*L$.

\subhead 4.3\endsubhead
Let 
$$I^s_n=I_n\cap I^s.$$
Note that if $w\cdo\l\in I^s_n$ and $\o\in\k_0\i(w)$, then $\cl_{\l,s}^\o|_{Z_s^w}$, 
$\Bbb L_{\l,s}^\o$ can be regarded naturally as objects in the mixed derived category of pure 
weight zero. Moreover, $\cl_{\l,s}^\o|_{Z_s^w}$ (resp. $\Bbb L_{\l,s}^\o$) is (noncanonically) 
isomorphic to $\cl_{\l,s}^{\dw}|_{Z_s^w}$ (resp. $\Bbb L_{\l,s}^{\dw}$) 
in the mixed derived category. 

We define $\ti\e_s:\cd(Z_s)@>>>\cd(\tcb^2)$, $\ti\e_s:\cd_m(Z_s)@>>>\cd_m(\tcb^2)$ by
$$\ti\e_s(L)=\e_s^*(L)\la\r\ra.$$
From the definition we have 
$$\e_s^*\cl_{\l,s}^{\o\sha}=L_\l^{\o\sha},\qua\ti\e_s\Bbb L_{\l,s}^\o=\LL_\l^\o.$$
Let $\cd^\spa Z_s$ be the subcategory of $\cd(Z_s)$ consisting of objects which
are restrictions of objects in the $G\T\TT^2$-equivariant derived category. Let 
$\cm^\spa Z_s$ be the subcategory of $\cd^\spa Z_s$ consisting of objects which are perverse 
sheaves. Let $\cm^{\preceq}Z_s$ (resp. $\cm^{\prec}Z_s$) be the subcategory
of $\cd^\spa Z_s$ whose objects are perverse sheaves $L$ such that any 
composition factor of $L$ is of the form $\Bbb L_{\l,s}^{\dw}$ for some 
$w\cdo\l\in I^s_n$ such that $w\cdo\l\preceq\boc$ (resp. $w\cdo\l\prec\boc$). 
Let $\cd^{\preceq}Z_s$ (resp. $\cd^{\prec}Z_s$) be the subcategory of $\cd^\spa Z_s$ whose objects
are complexes $L$ such that $L^j$ is in $\cm^{\preceq}Z_s$ (resp. $\cm^{\prec}Z_s$) for any $j$.
We write $\cd_m()$ or $\cm_m()$ for the mixed version of any of the categories above.

Let $\cc^\spa Z_s$ be the subcategory of $\cm^\spa Z_s$ consisting of semisimple objects. Let 
$\cc_0^\spa Z_s$ be the subcategory of $\cm^\spa_mZ_s$ consisting of objects of 
pure of weight zero. Let $\cc^\boc Z_s$ be the subcategory of $\cm^\spa Z_s$ consisting of objects 
which are direct sums of objects of the form $\Bbb L_{\l,s}^{\dw}$ with 
$w\cdo\l\in\boc^s$. Let $\cc_0^\boc Z_s$ be the subcategory of $\cc^\spa_0Z_s$ 
consisting of those $L\in\cc^\spa_0Z_s$ such that, as an object of $\cc^\spa Z_s$,
$L$ belongs to $\cc^\boc Z_s$. For $L\in\cc_0^\spa Z_s$ let $\un{L}$ be the 
largest subobject of $L$ such that as an object of $\cc^\spa Z_s$, we have 
$\un{L}\in\cc^\boc Z_s$.

From 4.2(a) we see that, if $M\in\cm^\spa Z_s$, then any composition factor of $M$ is of the form 
$\Bbb L_{\l,s}^{\dw}$ for some $w\cdo\l\in I^s_n$. From the definitions we see that 
$M\m\ti\e_s M$ is a functor $\cd^\spa Z_s@>>>\cd^\spa\tcb^2$ and also
$\cd_m^\spa Z_s@>>>\cd_m^\spa\tcb^2$; moreover, it is a functor 
$\cm^\spa Z_s@>>>\cm^\spa\tcb^2$ and also $\cm_m^\spa Z_s@>>>\cm_m^\spa\tcb^2$. From the 
definitions we see that for $M\in\cm^\spa Z_s$ 

(a) {\it we have $M\in\cm^{\preceq}Z_s$ if and only if $\ti\e_s M\in\cm^\preceq\tcb^2$; we have 
$M\in\cm^{\prec}Z_s$ if and only if $\ti\e_s M\in\cm^\prec\tcb^2$.}
\nl
Note that if $X\in\cd(Z_s)$ and $j\in\ZZ$, then 
$$(\e_s^*X)^{j+\r}=\e_s^*(X^j)[\r].\tag b$$
Moreover, if $Y\in\cm_m(Z_s)$ and $j'\in\ZZ$ then 
$$gr_{j'}(\ti\e_s Y)=\ti\e_s(gr_{j'}Y).\tag c$$

\mpb

For $w\cdo\l\in I_n$ we show:

(d) {\it We have $w\cdot\l\in I^s_n$ if and only if $w\i\cdot w(\l\i)\in I^{-s}_n$.}
\nl
We must show that we have $w(\l)=\ee^{-s}(\l)$ if and only if $\l\i=\ee^s(w(\l\i))$. In other 
words, we must show that $\l(w\i(t))=\l(\t^st\t^{-s})$ for all $t\in\TT_n$ if and only if
$\l(t')=\l(w\i(\t^{-s}t'\t^s))$ for all $t'\in\TT_n$. Setting $t'=\t^st\t^{-s}$, we have 
$w\i(t)= w\i(\t^{-s}t'\t^s)$ and it remains to use that $t\m\t^st\t^{-s}$ is a bijection 
$\TT_n@>>>\TT_n$.

\mpb

For $w\cdot\l\in I^s_n$ we show:

(e) {\it Let $\o\in\k_0\i(w)$. We have canonically
$(\Bbb L_{\l,s}^\o)^\da=\Bbb L_{w(\l\i),-s}^{\o\i}$.}
\nl
(The equality in (e) makes sense in view of (d).) By \cite{\MONO, 2.2(a)} and with notation of 3.1
we have canonically
$\ti\fh^*\LL_\l^\o=\LL_{w(\l\i)}^{\o\i}$. Hence
$\e_{-s}^*\LL_{w(\l\i)}^{\o\i}=\e_{-s}^*\ti\fh^*\LL_\l^\o=\fh^*\e_s^*\LL_\l^\o$
so that $\ti\e_{-s}\LL_{w(\l\i)}^{\o\i}=\fh^*\ti\e_s\LL_\l^\o$ and (e) follows.

\subhead 4.4\endsubhead
Let $r,f$ be integers such that $0\le f\le r-3$. Let
$$\align&\cy\\&
=\{((x_0\UU,x_1\UU,\do,x_r\UU),\g)\in\tcb^{r+1}\T\tG_s;\g\in x_{f+3}\UU\t^s x_f\i,
\g\in x_{f+2}\BB\t^s x_{f+1}\i\}.\endalign$$
Define $\vt:\cy@>>>\tcb^{r+1}$ by
$((x_0\UU,x_1\UU,\do,x_r\UU),\g)\m(x_0\UU,x_1\UU,\do,x_r\UU).$
For $y',y''\in W$ let
$$\tcb^{r+1}_{[y',y'']}=\{(x_0\UU,x_1\UU,\do,x_r\UU)\in\tcb^{r+1};
x_f\i x_{f+1}\in G_{y'},x_{f+2}\i x_{f+3}\in G_{y''{}\i}\}.$$

We show:

(a) {\it Let $\x\in\tcb^{r+1}_{[y',y'']}$. If $\ee^s(y')\ne y''$ then $\vt\i(\x)=\emp$. If
$\ee^s(y')=y''$ then $\vt\i(\x)\cong\kk^{\nu-|y'|}$.}
\nl
We set $\x=(x_0\UU,x_1\UU,\do,x_r\UU)$. If $\vt\i(\x)\ne\emp$ then 
$x_f\i x_{f+1}\in G_{y'}$, $x_{f+2}\i x_{f+3}\in G_{y''{}\i}$ and 
$(x_{f+3}\UU\t^s x_f\i)\cap(x_{f+2}\BB\t^s x_{f+1}\i)\ne\emp$. Hence for some 
$u\in\UU$, $b\in\BB$ we have 
$$u\t^sx_f\i x_{f+1}=x_{f+3}\i x_{f+2}b\t^s\in\t^s G_{y'}\cap G_{y''}\t^s$$
so that $\ee^s(y')=y''$. If we assume that $\ee^s(y')=y''$, then $\vt\i(\x)$ can be identified with
$$\{\g\in\tG_s;\g\in x_{f+3}\UU\t^s x_f\i,\g\in x_{f+2}\BB\t^s x_{f+1}\i\}$$
hence with 
$$\{(u,b)\in\UU\T\BB;u\t^sx_f\i x_{f+1}=x_{f+3}\i x_{f+2}b\t^s\}.$$
We substitute $x_{f+3}\i x_{f+2}=u_0\ee^s(\dy')t_0u'_0$, $x_f\i x_{f+1}=u_1\dy't_1u'_1$, 
where $t_0\in\TT$, \lb $u_0,u'_0,u_1,u'_1\in\UU$. 
Then $\vt\i(\x)$ is identified with 
$\{(u,b)\in\UU\T\BB;u\t^su_1\dy't_1u'_1=u_0\ee^s(\dy')t_0u'_0b\t^s\}$. The map 
$(u,b)\m u_0\i u\ee^s(u_1)$ identifies this variety with
$\UU\cap\ee^s(\dy')\BB\ee^s(\dy')\i\cong\kk^{\nu-|y'|}$. This proves (a).

\mpb

Now $\TT^2$ acts freely on $\cy$ by
$$\align&(t_1,t_2):((x_0\UU,x_1\UU,\do,x_r\UU),\g)\m\\&((x_0\UU,x_1\UU,\do,x_f\UU,
x_{f+1}t_1\UU,x_{f+2}t_2\UU,x_{f+3}\UU,\do,x_r\UU),\g).\endalign$$
Let
$$\align&{}^!\cy=\TT\bsl\{((x_0\UU,x_1\UU,\do,x_r\UU),\g)\in\tcb^{r+1}\T\tG_s;\\&
\g\in x_{f+3}\UU\t^s x_f\i,\g\in x_{f+2}\UU\t^s x_{f+1}\i\}\endalign$$
where $\TT$ acts freely by
$$\align&t:((x_0\UU,x_1\UU,\do,x_r\UU),\g)\m\\&((x_0\UU,x_1\UU,\do,x_f\UU,x_{f+1}\ee^{-s}(t)\UU,
x_{f+2}t\UU,x_{f+3}\UU,\do,x_r\UU),\g).\endalign$$
Note that the obvious map $\b:{}^!\cy@>>>\TT^2\bsl\cy$ is an isomorphism. We define 
${}^!\et:{}^!\cy@>>>Z_s$ by
$$((x_0\UU,x_1\UU,\do,x_r\UU),\g)\m\e_s(x_{f+1}\UU,x_{f+2}\UU).$$
We define $\pmb\t:\cy@>>>{}^!\cy$ as the composition of the obvious map
$\cy@>>>\TT^2\bsl\cy$ with $\b\i$. Let $\et={}^!\et\pmb\t:\cy@>>>Z_s$. We have
$$\et((x_0\UU,x_1\UU,\do,x_r\UU),\g)=\e_s(x_{f+1}t\i\UU,x_{f+2}t'{}\i\UU)$$
where $t,t'$ in $\TT$ are such that $\g\in x_{f+2}t'{}\i\UU\t^s tx_{f+1}\i$. 

\subhead 4.5\endsubhead
Let $z\cdo\l\in I^s_n$. Let $P=\et^*\cl_{\l,s}^{\dz\sha}$. Let $p_{ij}:\tcb^{r+1}@>>>\tcb^2$ be 
the projection to the $ij$ coordinates. We have the following result:
$$\vt_!P\Bpq\{p_{f,f+1}^*L_{\ee^{-s}(\l)}^{\ee^{-s}(\dy)}\ot p_{f+1,f+2}^*L_\l^{\dz\sha}\ot
p_{f+2,f+3}^*L_{y(\l)}^{\dy\i}\la2|y|-2\nu\ra;y\in W\}.\tag a$$
Define $e:\tcb^{r+1}@>>>\tcb^4$ by
$$(x_0\UU,x_1\UU,\do,x_r\UU)\m(x_f\UU,x_{f+1}\UU,x_{f+2}\UU,x_{f+3}\UU).$$
Then (a) is obtained by applying $e^*$ to the statement similar to (a) in
which $\{0,1,\do,r\}$ is replaced by $\{f,f+1,f+2,f+3\}$. Thus it is enough to prove (a) in the 
special case where $r=3,f=0$. In the remainder of the proof we assume that $r=3,f=0$.

For any $y',y''$ in $W$ let $\vt_{y',y''}:\vt\i(\tcb^4_{[y',y'']})@>>>\tcb^4$ be the restriction 
of $\vt$. Let $P^{y',y''}$ be the restriction of $P$ to $\vt\i(\tcb^4)_{[y',y'']}$. Clearly, we 
have
$$\vt_!P\Bpq\{(\vt_{y',y''})_!P^{y',y''};(y',y'')\in W^2\}.$$
Since $\vt\i(\tcb^{r+1}_{[y',y'']})=\emp$ when $\ee^s(y')\ne y''$, see 4.4(a), we deduce that
$$\vt_!P\Bpq\{(\vt_{\ee^{-s}(y),y\i})_!P^{\ee^{-s}(y),y\i};y\in W\}.$$
Hence to prove (a) it is enough to show for any $y\in W$ that
$$\vt_{y!}P_y=p_{01}^*L_{\ee^{-s}(\l)}^{\ee^{-s}(\dy)}\ot p_{12}^*L_\l^{\dz\sha}\ot 
p_{23}^*L_{y(\l)}^{\dy\i}\la2|y|-2\nu\ra,$$
where we write $\vt_y,P_y$ instead of $\vt_{\ee^{-s}(y),y\i},P^{\ee^{-s}(y),y\i}$. 
Using $z(\l)=\ee^{-s}(\l)$ we can replace $p_{01}^*L_{\ee^{-s}(\l)}^{\ee^{-s}(\dy)}$ by
$p_{01}^*L_{z(\l)}^{\ee^{-s}(\dy)}$. Thus it is enough to show for any $y\in W$ that
$$\vt_{y!}P_y=p_{01}^*L_{z(\l)}^{\ee^{-s}(\dy)}\ot p_{12}^*L_\l^{\dz\sha}\ot 
p_{23}^*L_{y(\l)}^{\dy\i}\la2|y|-2\nu\ra.\tag b$$
We have a cartesian diagram
$$\CD
\tV_y@>\tb>>\ti\cv_y\\
@VVV     @VVV\\
V_y@>b>>\cv_y
\endCD$$
where 
$$V_y=\{(x_0\UU,x_1\UU,x_2\UU,x_3\UU)\in\tcb^4;x_0\i x_1\in G_{\ee^{-s}(y)},
x_1\i x_2\in G_z,x_2\i x_3\in G_{y\i}\},$$
$$\align&\cv_y=\TT\bsl\{(x_0\UU,x_1\UU,x_2\UU,x_3\UU)\in\tcb^4;x_0\i x_1\in G_{\ee^{-s}(y)},
x_1\i x_2\in G_z,\\&x_2\i x_3\in G_{y\i},\ee^s((x_0\i x_1)_{\ee^{-s}(\dy)})=(x_3\i x_2)_{\dy}\}
\endalign$$
with $\TT$ acting (freely) by 
$$t:(x_0\UU,x_1\UU,x_2\UU,x_3\UU)\m(x_0\UU,x_1\ee^{-s}(t)\UU,x_2t\UU,x_3\UU),$$ 
$\tV_y=\vt\i(V_y)$ and
$$\align&\ti\cv_y=\TT\bsl\{((x_0\UU,x_1\UU,x_2\UU,x_3\UU),\g)\in\tcb^4\T\tG_s;x_0\i x_1\in 
G_{\ee^{-s}(y)},x_1\i x_2\in G_z, \\& x_2\i x_3\in G_{y\i},\g\in x_3\UU\t^s x_0\i,
\g\in x_2\UU\t^s x_1\i\}\endalign$$
with $\TT$ acting (freely) by 
$$t:((x_0\UU,x_1\UU,x_2\UU,x_3\UU).\g)\m((x_0\UU,x_1\ee^{-s}(t)\UU,x_2t\UU,x_3\UU),\g);$$ 
we have 
$$b(x_0\UU,x_1\UU,x_2\UU,x_3\UU)=\TT-\text{orbit of }(x_0\UU,x_1t\UU,x_2t'\UU,x_3\UU)$$
where $t,t'$ in $\TT$ are such that $\ee^s((x_0\i x_1t)_{\ee^{-s}(\dy)})=(x_3\i x_2t')_{\dy}$,
$$\tb((x_0\UU,x_1\UU,x_2\UU,x_3\UU),\g)=\TT-\text{orbit of }
((x_0\UU,x_1t\UU,x_2t'\UU,x_3\UU),\g)$$
where $t,t'$ in $\TT$ are such that $\g\in x_2t'\UU\t^st\i x_1\i$;  
the vertical maps are the obvious ones. We also have a cartesian diagram
$$\CD
\tV'_y@>\tb'>>\ti\cv'_y\\
@VVV     @VVV\\
V'_y@>b'>>\cv'_y
\endCD$$
where $\tV'_y,\ti\cv'_y,V'_y,\cv'_y$ are defined in the same way as
$\tV_y,\ti\cv_y,V_y,\cv_y$ but the condition $x_1\i x_2\in G_z$ is replaced by
the condition $x_1\i x_2\in\bG_z$; the maps $\tb',b'$ are given by the same 
formulas as $\tb,b$; the vertical maps are the obvious ones.

Let $j:V'_y@>>>\tcb^4$ be the inclusion. It is enough to show that
$$j^*\vt_{y!}P_y=j^*(p_{01}^*L_{z(\l)}^{\ee^{-s}(\dy)}\ot p_{12}^*L_\l^{\dz\sha}\ot 
p_{23}^*L_{y(\l)}^{\dy\i})\la2|y|-2\nu\ra.$$
By definition, $P|_{\tV'_y}$ is the inverse image of $\cl_{\l,s}^{\dz\sha}$ under the composition 
of $\tb'$ with $\ti\cv'_y@>>>\cv'_y@>{}^!\et_y>>Z_s$ where the first map is the obvious one and
$${}^!\et_y(x_0\UU,x_1\UU,x_2\UU,x_3\UU)=\e_s(x_1\UU,x_2\UU).$$
Hence $P|_{\tV'_y}$ is the inverse image of $\cl_{\l,s}^{\dz\sha}$ under the composition of 
$\et_y:={}^!\et_yb'$ with the obvious map $\vt'_y:\tV'_y@>>>V'_y$. Since $\vt_y$ is an affine 
space bundle with fibres of dimension $\nu-|y|$, it follows that 
$j^*\vt_{y!}P_y=\et_y^*\cl_{\l,s}^{\dz\sha}\la2|y|-2\nu\ra$. Thus it is enough to show that 
$$\et_y^*\cl_{\l,s}^{\dz\sha}
=j^*(p_{01}^*L_{z(\l)}^{\ee^{-s}(\dy)}\ot p_{12}^*L_\l^{\dz\sha}\ot p_{23}^*L_{y(\l)}^{\dy\i}).$$
Since $\et_y$ is smooth as a map to $\bZ_s^z$, we see that
$\et_y^*\cl_{\l,s}^{\dz\sha}$ is the intersection cohomology complex of $V'_y$ 
with coefficients in the local system $(\et^0_y)^*\cl_{\l,s}^{\dz}$ on $V_y$; here 
$\et^0_y:V_y@>>>Z_s^z$ is the restriction of $\et_y:V'_y@>>>\bZ_s^z$. By 3.9(a),
$$j^*(p_{01}^*L_{z(\l)}^{\ee^{-s}(\dy)}\ot p_{12}^*L_\l^{\dz\sha}\ot p_{23}^*L_{y(\l)}^{\dy\i})$$
is the intersection cohomology complex of $V'_y$ with coefficients in the 
local system 
$$\tL=j^*(p_{01}^*L_{z(\l)}^{\ee^{-s}(\dy)}\ot p_{12}^*L_\l^{\dz}\ot p_{23}^*L_{y(\l)}^{\dy\i})$$
on $V_y$. It is then enough to show that $\tL=(\et^0_y)^*\cl_{\l,s}^{\dz}$.

Let $\x=(x_0\UU,x_1\UU,x_2\UU,x_3\UU)\in V_y$. From the definition of $\et^0_y$ we see that the 
stalk $((\et^0_y)^*\cl_{\l,s}^{\dz})_\x$ is equal to 
$$(\cl_{\l,s}^{\dz})_{\e_s(x_1t_1\i,x_2t_2\i)}=(L_\l)_{t_0}$$
where $t_0\in\TT$, $t_1\in\TT$, $t_2\in\TT$ are such that
$t_0=(t_1x_1\i x_2t_2\i)_{\dz}$,  
$$\ee^s((x_0\i x_1t_1\i)_{\ee^{-s}(\dy)})=(x_3\i x_2t_2\i)_{\dy},$$
We can choose $t_1,t_2$ so that
$$(x_0\i x_1t_1\i)_{\ee^{-s}(\dy)}=1, (x_3\i x_2t_2\i)_{\dy}=1;$$
thus we can assume that $t_1=(x_0\i x_1)_{\ee^{-s}(\dy)}$, $t_2=(x_3\i x_2)_{\dy}=1$.

The stalk $\tL_\x$ is $(L_{z(\l)})_{t'_1}\ot(L_\l)_{t'_2}\ot(L_{y(\l)})_{t'_3}$ where 
$$t'_1=(x_0\i x_1)_{\ee^{-s}(\dy)}\in\TT,t'_2=(x_1\i x_2)_{\dz}\in\TT,
t'_3=(x_2\i x_3)_{\dy\i}\in\TT.$$
It is enough to show that $(\et_y^*\cl_{\l,s}^{\dz})_\x=\tL_\x$, or that 
$$(t_1x_1\i x_2t_2\i)_{\dz}=z\i(t'_1)t'_2y\i(t'_3)$$
where $t_1,t_2,t'_1,t'_2,t'_3$ are as above. We have $t_1=t'_1$ and $x_3\i x_2\in\UU\dy t_2\UU,$ 
hence 
$$x_2\i x_3\in\UU t_2\i\dy\i\UU=\UU\dy\i y(t_2\i)\UU,$$ 
so that $t'_3=y(t_2\i)$ and $t_2\i=y\i(t'_3)$. We have
$$t_1x_1\i x_2t_2\i\in t_1\UU\dz t'_2\UU t_2\i=\UU\dz z\i(t_1)t'_2t_2\i\UU,$$
so that
$$(t_1x_1\i x_2t_2\i)_{\dz}=z\i(t_1)t'_2t_2\i=z\i(t'_1)t'_2y\i(t'_3),$$
as required. This completes the proof of (b) hence that of (a).

\subhead 4.6\endsubhead
Let 
$$(w_1,w_2,\do,w_f,w_{f+2},w_{f+4},\do,w_r)\in W^{r-2},$$
$$(\l_1,\l_2,\do,\l_f,\l_{f+2},\l_{f+4},\do,\l_r)\in\fs_n^{r-2}.$$ 
We set $z=w_{f+2},\l=\l_{f+2}$. We assume that $z(\l)=\ee^{-s}(\l)$. Let $P$ be as in 4.5. Let 
$$P'=\ot_{i\in[1,r]-\{f+1,f+2,f+3\}}p_{i-1,i}^*L_{\l_i}^{\dw_i\sha}\in\cd_m(\tcb^{r+1}),$$
$\tP=P\ot\vt^*P'\in\cd_m(\cy)$. For any $y\in W$ we set
$$\ww_y=(w_1,w_2,\do,w_f,\ee^{-s}(y),w_{f+2},y\i,w_{f+4},\do,w_r)\in W^r,$$
$$\pmb\o_y=(\dw_1,\dw_2,\do,\dw_f,\ee^{-s}(\dy),\dw_{f+2},\dy\i,\dw_{f+4},\do,\dw_r),$$
$$\pmb\l_y=(\l_1,\l_2,\do,\l_f,\ee^{-s}(\l_{f+2}),\l_{f+2},y(\l_{f+2}),\l_{f+4},\do,\l_r)
\in\fs_n^r.$$
We set $\Xi=\vt_!\tP$. We have:
$$\Xi\Bpq\{M^{\pmb\o_y,[1,r]-\{f+1,f+3\}}_{\pmb\l_y}\la2|y|-2\nu\ra;y\in W\}\tag a$$
in $\cd_m(\tcb^{r+1})$. This follows immediately from 4.5(a) since $\Xi=P'\ot\vt_!(P)$.

\subhead 4.7\endsubhead
We preserve the setup of 4.6. Let $\cs=\sqc_{\ww'}\tco_{\ww'}^\emp$ where the union is over all
$\ww'=(w'_1,\do,w'_r)\in W^r$ such that $w'_i=w_i$ for $i\n\{f+1,f+3\}$. This is a locally closed 
subvariety of $\tcb^{r+1}$. For $y\in W$ let $R_y$ be the restriction of 
$M^{\pmb\o_y,\emp}_{\pmb\l_y}$ to $\tco^\emp_{\ww_y}$ extended by $0$ on 
$\cs-\tco^\emp_{\ww_y}$ (a constructible sheaf on $\cs$). From the definitions we have
$$M^{\pmb\o_y,[1,r]-\{f+1,f+3\}}_{\pmb\l_y}|_\cs=R_y.$$
From 4.6(a) we deduce $\Xi|_{\cs}\Bpq\{R_y\la2|y|-2\nu\ra;y\in W\}$.
We now restrict further to $\tco^\emp_{\ww_y}$ (for $y\in W$); we obtain
$$\Xi|_{\tco^\emp_{\ww_y}}
\Bpq\{R_{y'}\la2|y'|-2\nu\ra|_{\tco^\emp_{\ww_y}};y'\in W\}.$$
In the right hand side we have $R_{y'}\la2|y'|-2\nu\ra|_{\tco^\emp_{\ww_y}}=0$ 
if $y'\ne y$. It follows that $\Xi|_{\tco^\emp_{\ww_\l}}=R_y\la2|y|-2\nu\ra|_{\tco^\emp_{\ww_y}}$.
Since $R_y|_{\tco^\emp_{\ww_y}}$ is a local system we deduce for $y\in W$ the following result.

(a) {\it Let $h\in\ZZ$. If $h=2\nu-2|y|$ then
$\ch^h\Xi|_{\tco^\emp_{\ww_y}}=R_y|_{\tco^\emp_{\ww_y}}(|y|-\nu)$. If $h\ne2\nu-2|y|$, then 
$\ch^h\Xi|_{\tco^\emp_{\ww_y}}=0$.}  

\subhead 4.8\endsubhead
We preserve the setup of 4.6. We set
$$k=3\nu+(r+1)\r+\sum_{i\in[1,r]-\{f+1,f+3\}}|w_i|.\tag a$$ 
For $y\in W$ we set
$$K_y=M^{\pmb\o_y,[1,r]-\{f+1,f+3\}}_{\pmb\l_y}\la|\ww_y|+\nu+(r+1)\r\ra,$$
$$\tK_y=M^{\pmb\o_y,[1,r]}_{\pmb\l_y}\la|\ww_y|+\nu+(r+1)\r\ra.$$
From 4.6(a) we deduce:
$$\Xi\la k\ra\Bpq\{K_y;y\in W\}.\tag b$$
We show:

(c) {\it For any $j>0$ we have $(\Xi\la k\ra)^j=0$. Equivalently, $\Xi^j=0$ for any $j>k$.}
\nl
Using (b) we see that it is enough to show that for any $y\in W$ we have $(K_y)^j=0$ for any 
$j>0$. Now $\tK_y$ is a (simple) perverse sheaf hence for any $j$ we have 
$\dim\supp\ch^j\tK_y\le-j$. Moreover $K_y$ is obtained by restricting $\tK_y$ to an open subset 
of its support and then extending the result (by zero) on the complement of this subset in 
$\tcb^{r+1}$. Hence $\supp\ch^jK_y\sub\supp\ch^i\tK_y$ so that $\dim\supp\ch^iK_y\le-j$. Since
this holds for any $j$ we see that $(K_y)^j=0$ for any $j>0$. 

\subhead 4.9\endsubhead
We preserve the notation of 4.6. We show:

(a) {\it Let $j\in\ZZ$ and let $X$ be a composition factor of $\Xi^j$. Then\lb
$X\cong M^{\pmb\o',[1,r]}_{\pmb\l'}\la|\ww'|+\nu+(r+1)\r\ra$ for some 
$$\ww'=(w'_1,w'_2,\do,w'_r)\in W^r,\pmb\l'=(\l'_1,\l'_2,\do,\l'_r)\in\fs_n^r$$
such that $w'_i=w_i$, $\l'_i=\l_i$ for $i\in[1,r]-\{f+1,f+3\}$ and such that 
$$\l'_{f+1}=w'_{f+2}(\l'_{f+2}),\l'_{f+2}=w'_{f+3}(\l'_{f+3}).$$
Here $\pmb\o'=(\dw'_1,\dw'_2,\do,\dw'_r)$.}
\nl
From 4.6(a) we see that, for some $y\in W$, $X$ is a composition factor of 
$$(M^{\pmb\o_y,[1,r]-\{f+1,f+3\}}_{\pmb\l_y}\la2|y|-2\nu\ra)^j$$
where $\pmb\o_y,\pmb\l_y$ are as in 4.6.
Using this and \cite{\MONO, 2.18(b)} we see that 
$$X\cong M^{\pmb\o',[1,r]}_{\pmb\l'}\la|\ww'|+\nu+(r+1)\r\ra$$
for some 
$$\ww'=(w'_1,w'_2,\do,w'_r)\in W^r,\pmb\l'=(\l'_1,\l'_2,\do,\l'_r)\in\fs_n^r$$
such that $w'_i=w_i$, $\l'_i=\l_i$ for $i\in[1,r]-\{f+1,f+3\}$; here\lb
$\pmb\o'=(\dw'_1,\dw'_2,\do,\dw'_r)$. It remains to show that we have automatically
$$\l'_{f+1}=w'_{f+2}(\l'_{f+2}), \l'_{f+2}=w'_{f+3}(\l'_{f+3}).$$ To see this we note that 
$(M^{\pmb\o_y,[1,r]-\{f+1,f+3\}}_{\pmb\l_y}\la2|y|-2\nu\ra)^j$
is equivariant for the $\TT^2$-action
$$\align&(t_1,t_2):(x_0\UU,x_1\UU,\do,x_r\UU)\m  \\&
(x_0\UU,x_1\UU,\do,x_f\UU,x_{f+1}t_1\UU,x_{f+2}t_2\UU,x_{f+3}\UU,\do,x_r\UU)\endalign$$
hence so are its composition factors and this implies that the equalities above for 
$\l'_{f+1},\l'_{f+2}$ do hold.

\subhead 4.10\endsubhead
From 4.8(c) we see that we have a distinguished triangle $(\Xi',\Xi,\Xi^k[-k])$ where 
$\Xi'\in\cd_m(\tcb^{r+1})$ satisfies $(\Xi')^j=0$ for all $j\ge k$. We show:

(a) {\it Let $j\in\ZZ$ and let $K$ be one of $\Xi,\Xi^j,\Xi'$. For any $\ww'\in W^r$ and any 
$h\in\ZZ$, $\ch^hK|_{\tco^\emp_{\ww'}}$ is a local system.}
\nl
We prove (a) for $K=\Xi$ or $K=\Xi^j$. Using 4.6(a), we see that it is enough to show that 
$\ch^h(M^{\pmb\o_y,[1,r]-\{f+1,f+3\}}_{\pmb\l_y})|_{\tco^\emp_{\ww'}}$ is a 
local system for any $h$ and that
$\ch^h((M^{\pmb\o_y,[1,r]-\{f+1,f+3\}}_{\pmb\l_y})^j)|_{\tco^\emp_{\ww'}}$ 
is a local system for any $h$ and any $j$.  This follows by an argument entirely similar to that 
in the proof of \cite{\MONO, 3.10}.

Now (a) for $K=\Xi'$ follows from (a) for $\Xi$ and $\Xi^k[-k]$ using the long exact sequence for 
cohomology sheaves of $(\Xi',\Xi,\Xi^k[-k])$ restricted to $\tco^\emp_{\ww'}$.

\mpb

We show:

(b) {\it Let $(y',y'')\in W^2$, $j=2\nu-|y'|-|y''|$. Let 
$$\ww_{y',y''}=(w_1,w_2,\do,w_f,y',w_{f+2},y''{}\i,w_{f+3},\do,w_r)\in W^r.$$
The induced homomorphism $\ch^j\Xi|_{\tco^\emp_{\ww_{y',y''}}}@>>>
\ch^{j-k}(\Xi^k)|_{\tco^\emp_{\ww_{y',y''}}}$ is an isomorphism.}
\nl
We have an exact sequence of constructible sheaves
$$\ch^j\Xi'|_{\tco^\emp_{\ww_{y',y''}}}@>>>\ch^j\Xi|_{\tco^\emp_{\ww_{y',y''}}}@>>>
\ch^{j-k}(\Xi^k)|_{\tco^\emp_{\ww_{y',y''}}}@>>>\ch^{j+1}\Xi'|_{\tco^\emp_{\ww_{y',y''}}}.$$
Hence it is enough to show that 
$\ch^{j'}\Xi'|_{\tco^\emp_{\ww_{y',y''}}}=0$ if $j'\ge j$. Assume that \lb
$\ch^{j'}\Xi'|_{\tco^\emp_{\ww_{y',y''}}}\ne0$ for some $j'\ge j$. Since 
$\ch^{j'}\Xi'|_{\tco^\emp_{\ww_{y',y''}}}$ is a local system (see (a)), we deduce that 
$\tco^\emp_{\ww_{y',y''}}$ is contained in $\supp(\ch^{j'}\Xi')$. We have $(\Xi'[k-1])^{\tj}=0$ 
for any $\tj>0$ hence $\dim\supp(\ch^{j''}\Xi'[k-1])\le-j''$ for any $j''$. Taking $j''=j'-k+1$, 
we deduce that 
$$\dim\tco^\emp_{\ww_{y',y''}}\le\dim\supp(\ch^{j'}\Xi')\le-j'+k-1\le-j+k-1$$
hence
$$|\ww_{y',y''}|+\nu+(r+1)\r\le-j+k-1.$$
We have $|\ww_{y',y''}|+\nu+(r+1)\r=-j+k$ hence $-j+k\le-j+k-1$, contradiction. This proves (b).

\subhead 4.11\endsubhead
For $(y',y'')\in W^2$ we set
$$\pmb\o_{y',y''}=(\dw_1,\dw_2,\do,\dw_f,\dy',\dw_{f+2},\dy''{}\i,\dw_{f+3},\do,\dw_r)\in W^r,$$ 
$$\pmb\l_{y',y''}=(\l_1,\l_2,\do,\l_f,\ee^{-s}(\l_{f+2}),\l_{f+2},y''(\l_{f+2}),
\l_{f+4},\do,\l_r)\in\fs_n^r,$$
$$K_{y',y''}=M^{\pmb\o_{y',y''},\emp}_{\pmb\l_{y',y''}}\la|\ww_{y',y''}|+\nu+(r+1)\r\ra\in
\cm_m(\tcb^{r+1}),$$
$$\tK_{y',y''}=M^{\pmb\o_{y',y''},[1,r]}_{\pmb\l_{y',y''}}\la|\ww_{y',y''}|+\nu+(r+1)\r\ra\in
\cm_m(\tcb^{r+1}).$$
Note that when $y'=\ee^{-s}(y),y''=y$, $\ww_{y',y''},\pmb\o_{y',y''}$, $\pmb\l_{y',y''}$ and
$\tK_{y',y''}$ become $\ww_y,\pmb\o_y,\pmb\l_y$ (see 4.6) and $\tK_y$ (see 4.8). We show that we 
have canonically
$$gr_0(\Xi^k(k/2))=\op_{y\in W}\tK_y.\tag a$$
Since $gr_0(\Xi^k(k/2))$ is a semisimple perverse sheaf of pure weight zero, it is a direct sum 
of simple perverse sheaves, necessarily of the form described in 4.9(a). Thus we have canonically
$$gr_0(\Xi^k(k/2))=\op_{(y',y'')\in W^2}V_{y',y''}\ot\tK_{y',y''}$$
where $V_{y',y''}$ are mixed $\bbq$-vector spaces of pure weight $0$. By \cite{\BBD, 5.1.14}, 
$\Xi$ is mixed of weight $\le0$ hence $\Xi^k(k/2)$ is mixed of weight $\le0$. Hence we have an 
exact sequence in $\cm_m(\tcb^{r+1})$:
$$0@>>>\cw\i(\Xi^k(k/2))@>>>\Xi^k(k/2)@>>>gr_0(\Xi^k(k/2))@>>>0\tag a$$
that is,
$$0@>>>\cw\i(\Xi^k(k/2))@>>>\Xi^k(k/2)@>>>
\op_{(y',y'')\in W^2}V_{y',y''}\ot\tK_{y',y''}@>>>0.$$
(Here $\cw\i(?)$ denotes the part of weight $\le-1$ of a mixed perverse sheaf.) Hence for any 
$(\ty',\ty'')\in W^2$ we have an exact sequence of (mixed) cohomology sheaves restricted to 
$\tco^\emp_{\ww_{\ty',\ty''}}$ (where $h=2\nu-|\ty'|-|\ty''|-k$): 
$$\align&\ch^h(\cw\i(\Xi^k(k/2)))@>\a>>\ch^h(\Xi^k(k/2))@>>>
\op_{(y',y'')\in W^2}V_{y',y''}\ot\ch^h(\tK_{y',y''})@>>>\\&
\ch^{h+1}(\cw\i(\Xi^k(k/2))).\tag b\endalign$$
Moreover, by 4.10(b), we have an equality of local systems on $\tco^{\emp}_{\ww_{\ty',\ty''}}$:
$$\ch^h(\Xi^k(k/2))=\ch^{h+k}(\Xi(k/2))=\ch^{2\nu-|y'|-|y''|}(\Xi(k/2))$$
and this is $R_{y}(k/2+|y|-\nu)$ if $\ty'=\ee^{-s}(y),\ty''=y$ (see 4.7(a)) and is $0$
if $\ty'\ne\ee^{-s}(\ty'')$ (see 4.4(a)) hence is pure of weight $-k-|\ty'|-|\ty''|+\nu=h$.
On the other hand, $\ch^h(\cw\i(\Xi^k(k/2)))$ is mixed of weight $\le h-1$; it
follows that $\a$ in (b) must be zero.

Assume that $\ch^h(\tK_{y',y''})$ is not identically zero on $\tco^\emp_{\ww_{\ty',\ty''}}$. Then,
by 4.10(a), $\tco^\emp_{\ww_{\ty',\ty''}}$ is contained in $\supp\ch^h(\tK_{y',y''})$ which has 
dimension $\le-h$ (resp. $<-h$ if $(y',y'')\ne(\ty',\ty'')$); hence 
$-h=\dim\tco^\emp_{\ww_{\ty',\ty''}}$ is $\le-h$ (resp. $<-h$); we see that we must have 
$(y',y'')=(\ty',\ty'')$ and we have $\ch^h(\tK_{y',y''})=\ch^h(K_{y',y''})$ on 
$\tco^\emp_{\ww_{\ty',\ty''}}$.

Assume that $\ch^{h+1}(\cw\i(\Xi^k(k/2)))$ is not identically $0$ on 
$\tco^\emp_{\ww_{\ty',\ty''}}$. Then, by 4.10(a), $\tco^\emp_{\ww_{\ty',\ty''}}$ is contained in 
$\supp\ch^{h+1}(\cw\i(\Xi^k(k/2)))$ which has dimension $\le-h-1$; hence 
$-h=\dim\tco^\emp_{\ww_{\ty',\ty''}}\le-h-1$, a contradiction. 
We see that (b) becomes an isomorphism of local systems on $\tco^\emp_{\ww_{\ty',\ty''}}$:
$$0=V_{\ty',\ty''}\ot K_{\ty',\ty''}\text{ if } \ee^s(\ty')\ne\ty'',$$
$$R_{\ty''}(-h/2)@>\si>>V_{\ty',\ty''}\ot\ch^h(K_{\ty',\ty''})\text{ if }\ee^s(\ty')=\ty''.$$
When $\ee^s(\ty')=\ty'$ we have $\ch^h(K_{\ty',\ty''})=R_{\ty''}(-h/2)$ as local systems
on $\tco^\emp_{\ww_{\ty',\ty''}}$. It follows that $V_{\ty',\ty''}$ is $\bbq$ if
$\ee^s(\ty')=\ty''$ and is $0$ if $\ee^s(\ty')\ne\ty''$. This proves (a).

\subhead 4.12\endsubhead
Let $h\in[1,r]$. Let ${}_h\cd^{\preceq}\tcb^{r+1}$ (resp. ${}_h\cd^{\prec}\tcb^{r+1}$) be the 
subcategory of $\cd\tcb^{r+1}$ consisting of objects $K$ such that for any 
$j\in\ZZ$, any composition factor of $K^j$ is of the form 
$M^{\pmb\o,[1,r]}_{\pmb\l}\la|\ww|+\nu+(r+1)\r\ra$ for some $\ww=(w_1,\do,w_r)\in W^r$, 
$\pmb\l=(\l_1,\l_2,\do,\l_r)\in\fs_n^r$ such that $w_h\cdo\l_h\preceq\boc$ (resp. 
$w_h\cdo\l_h\prec\boc$). (Here $\pmb\o=(\dw_1,\dw_2,\do,\dw_r)$.)

Let ${}_h\cm^{\preceq}\tcb^{r+1}$ be the subcategory of ${}_h\cd^{\preceq}\tcb^{r+1}$ consisting 
of perverse sheaves. Let ${}_h\cm^{\prec}\tcb^{r+1}$ be the subcategory of 
${}_h\cd^{\prec}\tcb^{r+1}$ consisting of perverse sheaves.

If $K\in\cm_m(\tcb^{r+1})$ is pure of weight $0$ and is also in ${}_h\cd^{\preceq}\tcb^{r+1}$, we 
denote by $\un K$ the sum of all simple subobjects of $K$ (without mixed structure) which are not 
in ${}_h\cd^{\prec}\tcb^{r+1}$.

\subhead 4.13\endsubhead
Let $Z_s@<\et<<\cy@>\vt>>\tcb^4$ be as in 4.4 with $r=3,f=0$. We define 
$\fb:\cd(Z_s)@>>>\cd(\tcb^2)$ and $\fb:\cd_m(Z_s)@>>>\cd_m(\tcb^2)$ by 
$$\fb(L)=p_{03!}\vt_!\eta^*L.$$
We show:

(a) {\it If $L\in\cd^{\preceq}(Z_s)$ then $\fb(L)\in\cd^{\preceq}\tcb^2$.}

(b) {\it If $L\in\cd^{\prec}(Z_s)$ then $\fb(L)\in\cd^{\prec}\tcb^2$.}

(c) {\it If $L\in\cm^{\preceq}(Z_s)$ and $h>5\r+2\nu+2a$ then 
$(\fb(L))^h\in\cm^{\prec}\tcb^2$.}
\nl
We can assume that $L=\Bbb L_{\l,s}^{\dz}$ where $z\cdo\l\in I^s_n$, 
$z\cdo\l\preceq\boc$. Applying 4.5(a) with $P=\et^*\cl_{\l,s}^{\dz\sha}$ we see that
$$\fb(\cl_{\l,s}^{\dz\sha})\Bpq\{L^{\ee^{-s}(\dy),\dz,\dy\i,\{2\}}_{\ee^{-s}(\l),\l,y(\l)}
\la-|z|-2\nu\ra;y\in W\},$$
hence
$$\fb(\Bbb L_{\l,s}^{\dz\sha})\Bpq\{L^{\ee^{-s}(\dy),\dz,\dy\i,\{2\}}_{\ee^{-s}(\l),\l,y(\l)}
\la-\nu+\r\ra;y\in W\}.$$
To prove (a) it is enough to show that for any $y\in W$ we have
$$L^{\ee^{-s}(\dy),\dz,\dy\i,\{2\}}_{\ee^{-s}(\l),\l,y(\l)}\in\cd^{\preceq}\tcb^2.$$
When $z\cdo\l\in\boc$ this follows from \cite{\MONO, 2.10(a)}. When $z\cdo\l\prec\boc$ this 
again follows from \cite{\MONO, 2.10(a)}, applied to the two-sided cell containing $z\cdo\l$
instead of $\boc$. The same argument proves (b).
To prove (c) we can assume that $z\cdo\l\in\boc$; it is enough to prove that for any $y\in W$ we 
have 
$$(L^{\ee^{-s}(\dy),\dz,\dy\i,\{2\}}_{\ee^{-s}(\l),\l,y(\l)}\la-\nu+\r\ra)^h\in\cm^{\prec}\tcb^2$$
if $h>5\r+2\nu+2a$ or that 
$$(L^{\ee^{-s}(\dy),\dz,\dy\i,\{2\}}_{\ee^{-s}(\l),\l,y(\l)})^j\in\cm^{\prec}\tcb^2$$ 
if $j>6\r+\nu+2a$. This follows from \cite{\MONO, 2.20(a)}. This completes the proof of 
(a),(b),(c).

We define $\un\fb:\cc^\boc_0(Z_s)@>>>\cc^\boc_0(\tcb^2)$ by   
$$\un\fb(L)=\un{gr_{5\r+2\nu+2a}((\fb(L))^{5\r+2\nu+2a})}((5\r+2\nu+2a)/2).$$ 
We show:

(d) {\it Let $z\cdo\l\in\boc^s$. If $\ee^s(\boc)=\boc$, then
$$\un\fb(\Bbb L_{\l,s}^{\dz})=\op_{y\in W;y\cdo\l\in\boc}
\LL_{\ee^{-s}(\l)}^{\ee^{-s}(\dy)}\un{\cir}\LL_\l^{\dz}\un{\cir}\LL_{y(\l)}^{\dy\i}.$$
If $\ee^s(\boc)\ne\boc$, then $\un\fb(\Bbb L_{\l,s}^{\dz})=0$.}
\nl
We shall apply the method of \cite{\CONV, 1.12} with $\Ph:\cd_m(Y_1)@>>>\cd_m(Y_2)$ replaced by 
$p_{03!}:\cd_m(\tcb^4)@>>>\cd_m(\tcb^2)$ and with $\cd^{\preceq}(Y_1)$, $\cd^{\preceq}(Y_2)$ 
replaced by ${}_2\cd^{\preceq}(\tcb^2)$, ${}_2\cd^{\preceq}(\tcb^4)$, see 4.12. We shall take 
$\XX$ in {\it loc.cit.} equal to $\vt_!\et^*\Bbb L_{\l,s}^{\dz}$. The conditions of {\it loc.cit.}
are satisfied: those concerning $\XX$ are satisfied with $c'=2\nu+3\r$. (For $h>|z|+3\nu+4\r$ we 
have $\Xi^h=0$ that is $(\XX[-|z|-\nu-\r])^h=0$, with $\Xi$ as in 4.8(c). Hence if $j>2\nu+3\r$ we 
have $\XX^j=0$.) The conditions concerning $p_{03!}$ are satisfied with $c=2\r+2a$. (This follows 
from \cite{\MONO, 2.20(a)}.) Since $\fb(\Bbb L_{\l,s}^{\dz})=p_{03!}\XX$ and $c+c'=5\r+2\nu+2a$, 
we see that
$$\un\fb(\Bbb L_{\l,s}^{\dz})=\un{gr_{2\r+2a}(p_{03!}((\un{gr_{2\nu+3\r}
((\vt_!\et^*\Bbb L_{\l,s}^{\dz})^{2\nu+3\r})}((2\nu+3\r)/2)))^{2\r+2a})}(\r+a).$$
Using 4.11(a), we see that (with $\Xi$ as in 4.11(a) and $k=|z|+3\nu+4\r$) we have
$$\align&\un{gr_{2\nu+3\r}((\vt_!\et^*\Bbb L_{\l,s}^{\dz})^{2\nu+3\r})}((2\nu+3\r)/2)\\&=
\un{gr_{2\nu+3\r}((\Xi\la|z|+\nu+\r\ra)^{2\nu+3\r})}((2\nu+3\r)/2)\\&=\un{gr_0(\Xi^k(k/2)}
=\op_{y\in W}M^{\ee^{-s}(\dy),\dz,\dy\i,[1,3]}_{\ee^{-s}(\l),\l,y(\l)}\la2|y|+|z|+\nu+4\r\ra.
\endalign$$
Hence
$$\align&\un\fb(\Bbb L_{\l,s}^{\dz})=\un{gr_{2\r+2a}(\op_{y\in W}(p_{03!}
M^{\ee^{-s}(\dy),\dz,\dy\i,[1,3]}_{\ee^{-s}(\l),\l,y(\l)}\la2|y|+|z|+\nu+4\r\ra)^{2\r+2a})}(\r+a)
\\&=\un{gr_{2\r+2a}(\op_{y\in W}(L^{\ee^{-s}(\dy),\dz,\dy\i,[1,3]}
_{\ee^{-s}(\l),\l,y(\l)})^{6\r+\nu+2a}((\nu+4\r)/2))}(\r+a).\endalign$$
Using \cite{\MONO, 2.26(a)}, we see that in the last direct sum, the contribution of $y\in W$ is 
$0$ unless $y\cdo\l\in\boc$ and $\ee^{-s}(y)\cdo\ee^{-s}(\l)\in\boc$. We see that the last direct 
sum is zero unless $\ee^s(\boc)=\boc$. If $\ee^s(\boc)=\boc$, for the terms corresponding to $y$
such that $y\cdo\l\in\boc$, we may apply \cite{\MONO, 2.24(a)}. Now (d) follows.

\subhead 4.14\endsubhead
We set $\ZZ_\boc=\{s'\in\ZZ;\ee^{s'}(\boc)=\boc\}$. This is a subgroup of $\ZZ$.
In the remainder of this section we assume that $s\in\ZZ_\boc$.

Let $Z_s@<{}^!\et<<{}^!\cy$ be as in 4.4 with $r=3,f=0$. Let ${}^!\tcb^4$ be the space of 
orbits of the free $\TT^2$-action on $\tcb^4$ given by 
$$(t_1,t_2):(x_0\UU,x_1\UU,x_2\UU,x_3\UU)\m(x_0\UU,x_1t_1\UU,x_2t_2\UU,x_3\UU);$$
let ${}^!\vt:{}^!\cy@>>>{}^!\tcb^4$ be the map induced by $\vt$. We define 
$\fb':\cd(Z_s)@>>>\cd(\tcb^2)$ and $\fb':\cd_m(Z_s)@>>>\cd_m(\tcb^2)$ by 
$$\fb'(L)=p_{03!}{}^!\vt_!{}^!\et^*L.$$
(The map ${}^!\tcb^4@>>>\tcb^2$ induced by $p_{03}:\tcb^4@>>>\tcb^2$ is
denoted again by $p_{03}$.) Let $\t:\cy@>>>{}^!\cy$ be as in 4.4 (it is a 
principal $\TT^2$-bundle). We have the following results.

(a) {\it If $L\in\cd^{\preceq}(Z_s)$, then $\fb'(L)\in\cd^{\preceq}\tcb^2$.}

(b) {\it If $L\in\cd^{\prec}(Z_s)$, then $\fb'(L)\in\cd^{\prec}\tcb^2$.}

(c) {\it If $L\in\cm^{\preceq}(Z_s)$ and $h>\r+2\nu+2a$, then $(\fb'(L))^h\in\cm^{\prec}\tcb^2$.}
\nl
We can assume that $L=\Bbb L_{\l,s}^{\dz}$ where $z\cdo\l\in I^s_n$, 
$z\cdo\l\preceq\boc$. A variant of the proof of 4.5(a) gives:
$$\fb'(\cl_{\l,s}^{\dz\sha})\Bpq\{{}'L^{\ee^{-s}(\dy),\dz,\dy\i,\{2\}}_{\ee^{-s}(\l),\l,y(\l)}
\la-|z|-2\nu\ra;y\in W\},$$
hence
$$\fb'(\Bbb L_{\l,s}^{\dz\sha})\Bpq\{{}'L^{\ee^{-s}(\dy),\dz,\dy\i,\{2\}}_{\ee^{-s}(\l),\l,y(\l)}
\la-\nu+\r\ra;y\in W\}.$$
To prove (a) it is enough to show that for any $y\in W$ we have
$${}'L^{\ee^{-s}(\dy),\dz,\dy\i,\{2\}}_{\ee^{-s}(\l),\l,y(\l)}\in\cd^{\preceq}\tcb^2.$$
When $z\cdo\l\in\boc$ this follows from \cite{\MONO, 2.10(c)}. When $z\cdo\l\prec\boc$ this 
again follows from \cite{\MONO, 2.10(c)}, applied to the two-sided cell containing $z\cdo\l$
instead of $\boc$. The same argument proves (b). To prove (c) we can assume that $z\cdo\l\in\boc$;
it is enough to prove that for any $y\in W$ we have 
$$({}'L^{\ee^{-s}(\dy),\dz,\dy\i,\{2\}}_{\ee^{-s}(\l),\l,y(\l)}\la-\nu+\r\ra)^h\in
\cm^{\prec}\tcb^2$$
if $h>\r+2\nu+2a$ or that 
$({}'L^{\ee^{-s}(\dy),\dz,\dy\i,\{2\}}_{\ee^{-s}(\l),\l,y(\l)})^j\in\cm^{\prec}\tcb^2$
if $j>2\r+\nu+2a$. This follows from \cite{\MONO, 2.20(c)}. This completes the proof of 
(a),(b),(c).

We define $\un{\fb'}:\cc^\boc_0(Z_s)@>>>\cc^\boc_0(\tcb^2)$ by
$$\un{\fb'}(L)=\un{gr_{\r+2\nu+2a}((\fb'(L))^{\r+2\nu+2a})}((\r+2\nu+2a)/2).$$ 
In the remainder of this subsection we fix $z\cdo\l\in \boc^s$ and we set $L=\Bbb L_{\l,s}^{\dz}$.
We show:

(d) {\it We have canonically $\un{\fb'}(L)=\un\fb(L)$.}
\nl
The method of proof is similar to that of \cite{\MONO, 2.22(a)}. It is based on the fact that
$$\fb(L)=\fb'(L)\ot\fL^{\ot2}$$
which follows from the definitions. We define $\car_{i,j}$ for $i\in[0,2\r+1]$ and $\cp_{i,j}$ for
$i\in[0,2\r]$ as in \cite{\MONO, 2.17}, but replacing $L^J,{}'L^J,r,\d$ by $\fb(L),\fb'(L),3,2\r$.
In particular, we have
$$\cp_{i,j}=\cx_{4\r-i}(i-2\r)\ot(\fb'(L))^{-4\r+i+j}\text{ for }i\in[0,2\r]$$
where $\cx_{4\r-i}$ is a free abelian group of rank $\bin{2\r}{i}$ and $\cx_{4\r}=\ZZ$. We have 
for any $j$ an exact sequence analogous to \cite{\MONO, 2.17(a)}:
$$\do@>>>\cp_{i,j-1}@>>>\car_{i+1,j}@>>>\car_{i,j}@>>>\cp_{i,j}@>>>
\car_{i+1,j+1}@>>>\car_{i,j+1}@>>>\do,\tag e$$
and we have
$$\car_{0,j}=(\fb(L))^j,\qua \cp_{0,j}=(\fb'(L))^{j-4\r}(-2\r).$$
We show:

(f) {\it If $i\in[0,2\r+1]$ then $\car_{i,j}\in\cm^{\preceq}\tcb^2$.}

(g) {\it If $i\in[0,2\r+1]$, $j>6\r-i+\nu+2a$ then 
$\car_{i,j}\in\cm^{\prec}\tcb^2$.}
\nl
We prove (f),(g) by descending induction on $i$ as in \cite{\MONO, 2.21}. If $i=2\r+1$ then,
since $\car_{2\r+1,j}=0$, there is nothing to prove. Now assume that $i\in[0,2\r]$. Assume that 
$\l'\cdot w$ is such that $\LL_{\l'}^{\dw}$ is a composition factor of $\car_{i,j}$ (without the 
mixed structure). We must show that $w\cdo\l'\preceq\boc$ and that, if $j>6\r-i+\nu+2a$, then 
$w\cdo\l'\prec\boc$. Using (e), we see that $\LL_{\l'}^{\dw}$ is a composition factor of 
$\car_{i+1,j}$ or of $\cp_{i,j}$. In the first case, using the induction hypothesis we see that 
$w\cdo\l'\preceq\boc$ and that, if $j>6\r-i+\nu+2a$ (so that $j>6\r-i-1+\nu+2a$), then 
$w\cdo\l'\prec\boc$. In the second case, $\LL_{\l'}^{\dw}$ is a composition factor of 
$(\fb'(L))^{-4\r+i+j}$. Using (a),(c), we see that $w\cdo\l'\preceq\boc$ and that, if 
$j>6\r-i+\nu+2a$ (so that $-4\r+i+j>\nu+2\r+2a$), then $w\cdo\l'\prec\boc$. This proves (f),(g).

We show:

(h) {\it Assume that $i\in[0,2\r+1]$. Then $\car_{i,j}$ is mixed of weight $\le j-i$.}
\nl 
We argue as in \cite{\MONO, 2.22} by descending induction on $i$. If $i=2\r+1$ there is nothing to
prove. Assume now that $i\le2\r$. By Deligne's theorem, $\fb'(L)$ is mixed of weight $\le0$; hence
$(\fb'(L))^{-4\r+i+j}$ is mixed of weight $\le-4\r+i+j$ and 
$\cx_{4\r-i}(i-2\r)\ot(\fb'(L))^{-4\r+i+j}$ is mixed of weight $\le-4\r+i+j-2(i-2\r)=j-i$. In 
other words, $\cp_{i,j}$ is mixed of weight $\le j-i$. Thus in the exact sequence 
$\car_{i+1,j}@>>>\car_{i,j}@>>>\cp_{i,j}$ coming from (e) in which $\car_{i+1,j}$ is mixed of 
weight $\le j-i-1<j-i$ (by the induction hypothesis) and $\cp_{i,j}$ is mixed of weight $\le j-i$,
we must have that $\car_{i,j}$ is mixed of weight $\le j-i$. This proves (h).

\mpb

We now prove (d). From (e) we deduce an exact sequence 
$$gr_j(\car_{1,j})@>>>gr_j(\car_{0,j})@>>>gr_j(\cp_{0,j})@>>>gr_j(\car_{1,j+1}).$$
By (h) we have $gr_j(\car_{1,j})=0$. We have $gr_j(\car_{0,j})=gr_j(\fb(L)^j)$,
$gr_j(\cp_{0,j})=gr_j((\fb'(L))^{-4\r+j}(-2\r))$. Moreover, by (g) we have 
$\car_{1,j+1}\in\cd^{\prec}\tcb^2$ since $j+1>6\r-1+\nu+2a$. It follows that 
$gr_j(\car_{1,j+1})\in\cd^{\prec}\tcb^2$. Thus the exact sequence above induces an isomorphism as 
in (d).

\mpb

Let $p'_{ij}:\tcb^3@>>>\tcb^2$ be the projection to the $ij$-coordinate, where  $ij$ is
$12$, $23$ or $13$. Let 
$$R=\TT\bsl\{(x_0\UU,x_1\UU,x_2\UU,x'_3\UU,\g)\in\tcb^4\T G_s;\g\in x_2\UU\t^sx_1\i\}$$ 
where $\TT$ acts freely by 
$$t:(x_0\UU,x_1\UU,x_2\UU,x'_3\UU,\g)\m(x_0\UU,x_1\ee^{-s}(t)\UU,x_2t\UU,x'_3\UU,\g).$$
We have cartesian diagrams
$$\CD
R@>d_1>>{}'\cy\T\tcb^2  \\
@Vc_1VV    @Vs_1VV\\
\tcb^3@>p'>>\tcb^2\T\tcb^2
\endCD$$

$$\CD
R@>d_2>>\tcb^2\T{}'\cy  \\
@Vc_2VV    @Vs_2VV\\
\tcb^3@>p'>>\tcb^2\T\tcb^2
\endCD$$
where

$d_1(x_0\UU,x_1\UU,x_2\UU,x_3\UU,\g)=((x_0\UU,x_1\UU,x_2\UU,\g x_0\t^{-s}\UU,\g),
(\g x_0\t^{-s}\UU,x_3\UU))$,

$$\align&d_2(x_0\UU,x_1\UU,x_2\UU,x_3\UU,\g)\\&=((x_0\UU,\g\i x_3\t^s\UU),
(\g\i x_3\t^s\UU,x_1\UU,x_2\UU,x_3\UU,\g)),\endalign$$

$c_1(x_0\UU,x_1\UU,x_2\UU,x_3\UU,\g)=(x_0\UU,\g x_0\t^{-s}\UU,x_3\UU)$,

$c_2(x_0\UU,x_1\UU,x_2\UU,x_3\UU,\g)=(x_0\UU,\g\i x_3\t^s\UU,x_3\UU)$,

$p'=(p'_{12},p'_{23})$, $s_1=p_{03}{}'\vt\T1$, $s_2=1\T p_{03}{}'\vt$.
\nl
It follows that $p'{}^*s_{1!}= c_{1!}d_1^*$, $p'{}^*s_{2!}= c_{2!}d_2^*$.
Now let $L\in\cd(Z_s)$, $L'\in\cd(\tcb^2)$, $\tL'\in\cd(\tcb^2)$,
We have $\et^*L\bxt L'\in\cd({}'\cy\T\tcb^2$, $\tL'\bxt\et^*L\in\cd(\tcb^2\T{}'\cy)$. We have 
$$p'_{12}{}^*\fb'(L)\ot p'_{23}{}^*L'=
p'{}^*s_{1!}(\et^*L\bxt L')=c_{1!}d_1^*(\et^*L\bxt L')=c_{1!}(e_1^*L\bxt e'_1{}^*L'),$$
$$p'_{12}{}^*\tL'\ot p'_{23}{}^*\fb'(L)=
p'{}^*s_{2!}(\tL'\bxt \et^*L)=c_{2!}d_2^*(\tL'\bxt\et^*L)=c_{2!}(e'_2{}^*\tL'\bxt e_1{}^*L),$$
where 

$e_1:R@>>>Z_s$ is $(x_0\UU,x_1\UU,x_2\UU,x_3\UU,\g)\m\e_s(x_1\UU,x_2\UU)$,

$e'_1:R@>>>\tcb^2$ is $(x_0\UU,x_1\UU,x_2\UU,x_3\UU,\g)\m(\g x_0\t^{-s}\UU,x_3\UU)$,

$e'_2:R@>>>\tcb^2$ is $(x_0\UU,x_1\UU,x_2\UU,x_3\UU,\g)\m(x_0\UU,\g\i x_3\t^s\UU)$.

Applying $p'_{13!}$ we see that
$$\fb'(L)\cir L'=\tc_!(e_1^*L\bxt e'_1{}^*L), \tL'\cir\fb'(L)=\tc_!(e'_2{}^*L\bxt e_1{}^*L),$$
where $\tc:R@>>>\tcb^2$ is $(x_0\UU,x_1\UU,x_2\UU,x_3\UU,\g)\m(x_0\UU,x_3\UU)$.

We define $\ee:\tcb^2@>>>\tcb^2$ by $\ee(x\UU,y\UU)=(\ee(x)\UU,\ee(y)\UU)$. We show:
 
(i) {\it If in addition $L'\in\cm(\tcb^2)$ is $G$-equivariant, then we have canonically} 
$$\fb'(L)\cir L'=(\ee^{s*}L')\cir\fb'(L).$$
We take $\tL'=\ee^{s*}L'$. It is enough to show that 
$\tc_!(e_1^*L\bxt e'_1{}^*L')=\tc_!(e'_2{}^*\tL'\bxt e_1^*L)$. Hence it is enough to show that we 
have canonically $e'_1{}^*L'=e'_2{}^*\tL'$ that is, $e'_1{}^*L'=e''_2{}^*L'$ where 
$e''_2=\ee^s e'_2:R@>>>\tcb^2$. We identify $\tG_s$ with $G$ by $\g\m g$ where $\g=g\t^s$. Then 
$e'_1:R@>>>\tcb^2$ is $(x_0\UU,x_1\UU,x_2\UU,x_3\UU,\g)\m(g\ee^s(x_0)\UU,x_3\UU)$,
$e''_2:R@>>>\tcb^2$ is $(x_0\UU,x_1\UU,x_2\UU,x_3\UU,\g)\m(\ee^s(x_0)\UU,g\i x_3\UU)$. The 
equality $e'_1{}^*L'=e''_2{}^*L'$ follows from the $G$-equivariance of $L'$. This proves (i).

\mpb

We show:

(j) {\it If $L\in\cc^\boc_0Z_s$, $L'\in\cc^\boc\tcb^2$, then we have canonically 
$\un\fb(L)\un\cir L'=(\ee^{s*}L')\un\cir\un\fb(L)$.}
\nl
By (d), it is enough to prove that $\un\fb'(L)\un\cir L'=(\ee^{s*}L')\un\cir\un\fb'(L)$. Using (i)
together with (a),(b),(c) and results in \cite{\MONO, 2.23}, we see that both sides are equal to 
$$\align&\un{gr_{\r+\nu+3a}(\tc_!(e_1^*L\ot e'_1{}^*L'))^{\r+\nu+3a}}((\r+\nu+3a)/2)\\&
=\un{gr_{\r+\nu+3a}\tc_!(e_1^*L\ot e''_2{}^*L'))^{\r+\nu+3a}}((\r+\nu+3a)/2).\endalign$$

\subhead 4.15\endsubhead
Let 
$$\fZ_s=\{(z_0\UU,z_1\UU,z_2\UU,z_3\UU),\g)\in\tcb^4\T\tG_s;\g\in z_2\BB\t^sz_1\i\}.$$
Define $\ti\vt:\fZ_s@>>>\tcb^4$ by
$((z_0\UU,z_1\UU,z_2\UU,z_3\UU),\g)\m(z_0\UU,z_1\UU,z_2\UU,z_3\UU)$. Let 
$${}'\cy=\{((x_0\UU,x_1\UU,x_2\UU,x_3\UU,x_4\UU),\g)\in\tcb^5\T\tG_s;\g\in x_3\UU\t^s x_0\i, 
\g\in x_2\BB\t^s x_1\i\},$$
$${}''\cy=\{((x_0\UU,x_1\UU,x_2\UU,x_3\UU,x_4\UU),\g)\in\tcb^5\T\tG_s;\g\in x_4\UU\t^s x_1\i, 
\g\in x_3\BB\t^s x_2\i\}.$$
Define ${}'\vt:{}'\cy@>>>\tcb^5$, ${}''\vt:{}''\cy@>>>\tcb^5$ by
$$((x_0\UU,x_1\UU,x_2\UU,x_3\UU,x_4\UU),\g)\m(x_0\UU,x_1\UU,x_2\UU,x_3\UU,x_4\UU).$$
We have isomorphisms ${}'\fc:{}'\cy@>\si>>\fZ_s$, ${}''\fc:{}''\cy@>\si>>\fZ_s$ given by
$${}'\fc:((x_0\UU,x_1\UU,x_2\UU,x_3\UU,x_4\UU),\g)\m((x_0\UU,x_1\UU,x_2\UU,x_4\UU),\g),$$
$${}''\fc:((x_0\UU,x_1\UU,x_2\UU,x_3\UU,x_4\UU),\g)\m((x_0\UU,x_2\UU,x_3\UU,x_4\UU),\g).$$
Define ${}'d:\tcb^5@>>>\tcb^4$, ${}''d:\tcb^5@>>>\tcb^4$ by
$${}'d:(x_0\UU,x_1\UU,x_2\UU,x_3\UU,x_4\UU)\m(x_0\UU,x_1\UU,x_2\UU,x_4\UU),$$
$${}''d:(x_0\UU,x_1\UU,x_2\UU,x_3\UU,x_4\UU)\m(x_0\UU,x_2\UU,x_3\UU,x_4\UU).$$
We fix $w,u$ in $W$ and $\l,\l'$ in $\fs_n$. We assume that $w\cdot\l\in I^s_n$. The smooth 
subvarieties
$${}'\cu=\{((x_0\UU,x_1\UU,x_2\UU,x_3\UU,x_4\UU),\g)\in{}'\cy;x_1\i x_2\in G_w,
x_3\i x_4\in G_{\ee^s(u)}\},$$
$$\cu=\{((x_0\UU,x_1\UU,x_2\UU,x_3\UU),\g)\in\fZ_s;x_1\i x_2\in G_w,x_0\i g\i x_3\in G_u\},$$   
$${}''\cu=\{((x_0\UU,x_1\UU,x_2\UU,x_3\UU,x_4\UU),\g)\in{}''\cy;x_2\i x_3\in G_w,x_0\i x_1\in G_u
\},$$   
of ${}'\cy,\fZ_s,{}''\cy$ correspond to each other under the isomorphisms
${}'\cy@>{}'\fc>>\fZ_s@<{}''\fc<<{}''\cy$. Moreover, the maps ${}'\s:{}'\cu@>>>Z_s$, 
$\s:\cu@>>>Z_s$, ${}''\s:{}''\cu@>>>Z_s$ given by
$$((x_0\UU,x_1\UU,x_2\UU,x_3\UU,x_4\UU),\g)\m\e_s(x_1\UU,x_2\UU),$$    
$$((x_0\UU,x_1\UU,x_2\UU,x_3\UU),\g)\m\e_s(x_1\UU,x_2\UU),$$
$$((x_0\UU,x_1\UU,x_2\UU,x_3\UU,x_4\UU),\g)\m\e_s(x_2\UU,x_3\UU),$$
correspond to each other under the isomorphisms ${}'\cy@>{}'\fc>>\fZ_s@<{}''\fc<<{}''\cy$.

Also, the maps ${}'\ti\s:{}'\cu@>>>\tco_{\ee^s(u)}$, $\ti\s:\cu@>>>\tco_{\ee^s(u)}$, given by
$$((x_0\UU,x_1\UU,x_2\UU,x_3\UU,x_4\UU),\g)\m(x_3\UU,x_4\UU),$$
$$((x_0\UU,x_1\UU,x_2\UU,x_3\UU),\g)\m(\g x_0\t^{-s}\UU,x_3\UU)$$
correspond to each other under the isomorphism ${}'\cy@>{}'\fc>>\fZ_s$ and the 
maps $\ti\s_1:\cu@>>>\tco_u$, ${}''\ti\s:{}''\cu@>>>\tco_u$ given by
$$((x_0\UU,x_1\UU,x_2\UU,x_3\UU),\g)\m(x_0\UU,\g\i x_3\t^s\UU),$$
$$((x_0\UU,x_1\UU,x_2\UU,x_3\UU,x_4\UU),\g)\m(x_0\UU,x_1\UU),$$
correspond to each other under the isomorphism $\fZ_s@<{}''\fc<<{}''\cy$. It 
follows that the local systems ${}'\s^*\cl_{\l,s}^{\dw}$, $\s^*\cl_{\l,s}^{\dw}$, 
${}''\s^*\cl_{\l,s}^{\dw}$ correspond to each other under the isomorphisms 
${}'\cy@>{}'\fc>>\fZ_s@<{}''\fc<<{}''\cy$; the local systems 
${}'\ti\s^*L_{\ee^s(\l')}^{\ee^s\du)}$, $\ti\s^*L_{\ee^s(\l')}^{\ee^s(\du)}$ 
correspond to each other under the isomorphism ${}'\cy@>{}'\fc>>\fZ_s$; the local systems 
$\ti\s_1^*L_{\l'}^{\du}$, ${}''\ti\s^*L_{\l'}^{\du}$ correspond to each other 
under the isomorphism $\fZ_s@<{}''\fc<<{}''\cy$. Moreover, by the $G$-equivariance of 
$L_{\l'}^{\du}$, we have as in the proof of 4.14(i):
$\ti\s^*L_{\ee^s(\l')}^{\ee^s(\du)}=\ti\s_1^*(L_{\l'}^{\du})$. 

Let ${}'K,K,{}''K$ be the intersection cohomology complex of the closure of ${}'\cu,\cu,{}''\cu$ 
respectively with coefficients in the local system 
$${}'\s^*\cl_{\l,s}^{\dw}\ot{}'\ti\s^*L_{\ee^s(\l')}^{\ee^s(\du)},
\s^*\cl_{\l,s}^{\dw}\ot\ti\s^*L_{\ee^s(\l')}^{\ee^s(\du)}=
\s^*\cl_{\l,s}^{\dw}\ot\ti\s_1^*(L_{\l'}^{\du}),
{}''\s^*\cl_{\l,s}^{\dw}\ot{}''\ti\s^*L_{\l'}^{\du},$$
on ${}'\cu,\cu,{}''\cu$ (respectively), extended by $0$ on the complement of this closure in
${}'\cy,\fZ_s,{}''\cy$. We see that ${}'K,K,{}''K$ correspond to each other under the isomorphisms
${}'\cy@>{}'\fc>>\fZ_s@<{}''\fc<<{}''\cy$. Hence we have ${}'\fc_!({}'K)=K={}''\fc_!({}''K)$. 
Using this and the commutative diagram
$$\CD
{}'\cy@>{}'\fc>>\fZ_s@<{}''\fc<<{}''\cy\\
@V{}'\vt VV   @V\ti\vt VV     @V{}''\vt VV   \\
\tcb^5@>{}'d>>\tcb^4@<{}''d<<\tcb^5
\endCD$$
we see that 
$${}'d_!{}'\vt_!({}'K)={}''d_!{}''\vt_!({}''K).\tag a$$
(Both sides are equal to $\ti\vt_!K$.)

\subhead 4.16\endsubhead
In this subsection we study the functor ${}'d_!:\cd_m(\tcb^5)@>>>\cd_m(\tcb^4)$. Let 
$\ww=(w_1,w_2,w_3,w_4)$, $\pmb\l=(\l_1,\l_2,\l_3,\l_4)\in\fs_n^4$,
$\pmb\o=(\o_1,\o_2,\o_3,\o_4)$ (with $\o_i\in\k\i_0(w_i)$). Assume that
$w_4\cdo\l_4\preceq\boc$. Let $K=M_{\pmb\l}^{\pmb\o,[1,4]}\la|\ww|+5\r+\nu\ra\in\cd_m(\tcb^5)$.
As in \cite{\MONO, 3.16}, properties (a)(b),(c),(d) hold:

(a) {\it If $h>a+\r$ then $({}'d_!K)^h\in{}'\cm^{\prec}(\tcb^4)$. Moreover,}
$$\align&\un{gr_{a+\r}(({}'dK)^{a+\r})}((a+\r)/2)=\op_{y'\in W;y'{}\i\cdo\l_4\in\boc}
\Hom_{\cc^\boc\tcb^2}(\LL_{\l_4}^{\dy'{}\i},\LL_{\l_3}^{\o_3}\un{\cir}
\LL_{\l_4}^{\o_4})\\&\ot M_{\l_1,\l_2,\l_4}^{\o_1,\o_2,\dy'{}\i,[1,3]}\la
|w_1|+|w_2|+|y'|+4\r+\nu\ra.\endalign$$

(b) {\it If $K\in{}_4\cd^{\preceq}(\tcb^5)$ then ${}'d_!(K)\in{}_4\cd^{\preceq}(\tcb^4)$.}

(c) {\it If $K\in{}_4\cd^{\prec}(\tcb^5)$ then ${}'d_!(K)\in{}_4\cd^{\prec}(\tcb^4)$.}

(d) {\it If $K\in{}_4\cm^{\preceq}(\tcb^5)$ and $h>a+\r$ then 
$({}'d_!(K))^h\in{}_4\cm^{\prec}(\tcb^4)$.}

\subhead 4.17\endsubhead
In this subsection we study the functor ${}''d_!:\cd_m(\tcb^5)@>>>\cd_m(\tcb^4)$.
Let $\ww=(w_1,w_2,w_3,w_4)$, $\pmb\l=(\l_1,\l_2,\l_3,\l_4)\in\fs_n^4$,
$\pmb\o=(\o_1,\o_2,\o_3,\o_4)$ (with $\o_i\in\k\i_0(w_i)$). Assume that
$w_1\cdo\l_1\preceq\boc$. Let $K=M_{\pmb\l}^{\pmb\o,[1,4]}\la|\ww|+5\r+\nu\ra\in\cd_m(\tcb^5)$.
As in \cite{\MONO, 3.17}, properties (a)(b),(c),(d) hold:

(a) {\it If $h>a+\r$ then $({}''d_!K)^h\in{}'\cm^{\prec}(\tcb^4)$. Moreover,}
$$\align&\un{gr_{a+\r}(({}''d_!K)^{a+\r})}((a+\r)/2)=\op_{y'\in W;y'\cdo\l_2\in\boc}
\Hom_{\cc^\boc\tcb^2}(\LL_{\l_2}^{\dy'},\LL_{\l_1}^{\o_1}\un{\cir}\LL_{\l_2}^{\o_2})\\&
\ot M_{\l_2,\l_3,\l_4}^{\dy',\o_3,\o_4,[1,3]}\la|w_3|+|w_4|+|y'|+4\r+\nu\ra.\endalign$$

(b) {\it If $K\in{}_1\cd^{\preceq}(\tcb^5)$ then ${}''d_!(K)\in{}_1\cd^{\preceq}(\tcb^4)$.}

(c) {\it If $K\in{}_1\cd^{\prec}(\tcb^5)$ then ${}''d_!(K)\in{}_1\cd^{\prec}(\tcb^4)$.}

(d) {\it If $K\in{}_1\cm^{\preceq}(\tcb^5)$ and $h>a+\r$ then 
$({}''d_!(K))^h\in{}_1\cm^{\prec}(\tcb^4)$.}

\subhead 4.18\endsubhead
Let $w\cdo\l\in I^s_n$, $u\cdo\l'\in\boc$. We shall apply the method of \cite{\CONV, 1.12} with 
$\Ph:\cd_m(Y_1)@>>>\cd_m(Y_2)$ replaced by ${}'d_!:\cd_m(\tcb^5)@>>>\cd_m(\tcb^4)$ and with 
$\cd^{\preceq}(Y_1)$, $\cd^{\preceq}(Y_2)$ replaced by ${}_4\cd^{\preceq}(\tcb^5)$, 
${}_4\cd^{\preceq}(\tcb^4)$, see 4.15. We shall take $\XX$ in {\it loc.cit.} equal to 
$\Xi={}'\vt_!({}'K)$ as in 4.15, $(w_2,w_4)=(w,\ee^s(u))$, $(\l_2,\l_4)=(\l,\ee^s(\l'))$. The 
conditions of {\it loc.cit.} are satisfied: those concerning $\XX$ are satisfied with 
$c'=k=|w|+|u|+3\nu+5\r$ (see 4.8(c)); those concerning $\Ph$ are satisfied with $c=a+\r$ (see 
4.16). We see that
$$\align&\un{gr_{a+\r+k}(({}'d_!{}'\vt_!({}'K))^{a+\r+k})}((a+\r+k)/2)\\&=
\un{gr_{a+\r}(({}'d_!\un{gr_k(({}'\vt_!({}'K))^k)}(k/2))^{a+\r})}((a+\r)/2).\endalign$$
Using 4.11(a), we have:
$$\align&gr_k({}'\vt_!({}'K))^k)(k/2)   
=\op_{y\in W}M_{\ee^{-s}(\l),\l,y(\l),\ee^s(\l')}^{\ee^{-s}(\dy),\dw,\dy\i,\ee^s(\du),[1,4]}
\la2|y|+|w|+|u|+5\r+\nu\ra\\&=\un{gr_k({}'\vt_!({}'K))^k}(k/2).\endalign$$
Hence, using 4.16(a), we have
$$\align&\un{gr_{a+\r}(({}'d_!\un{gr_k(({}'\vt_!({}'K))^k)}(k/2))^{a+\r})}((a+\r)/2)\\&=
\op_{y\in W}\op_{y'\in W;y'{}\i\cdo \ee^s(\l')\in\boc}\Hom_{\cc^\boc\tcb^2}
(\LL_{\ee^s(\l')}^{\dy'{}\i},\LL_{y(\l)}^{\dy\i}\un{\cir}\LL_{\ee^s(\l')}^{\ee^s(\du)})\\&\ot 
M_{\ee^{-s}(\l),\l,\ee^s(\l')}^{\ee^{-s}(\dy),\dw,\dy'{}\i,[1,3]}\la|y|+|w|+|y'|+4\r+\nu\ra.
\endalign$$  
Since $y'{}\i\cdo \ee^s(\l')\in\boc$, $\ee^s(u)\cdo\ee^s(\l')\in\boc$ (recall that 
$\ee^s\boc=\boc$), for $y\in W$ we have 
$$\Hom_{\cc^\boc\tcb^2}(\LL_{\ee^s(\l')}^{\dy'{}\i},\LL_{y(\l)}^{\dy\i}\un{\cir}
\LL_{\ee^s(\l')}^{\ee^s(\du)})=0$$ 
unless $\ee^s(\l')=y'(\l)$ (see \cite{\MONO, 4.6(b)}) and $y\i\cdo y(\l)\in\boc$ (see 
\cite{\MONO, 2.26(a)}) or equivalently, $y\cdo\l\in\boc$. Thus we have
$$\align&\un{gr_{a+\r+k}(({}'d_!{}'\vt_!({}'K))^{a+\r+k})}((a+\r+k)/2)\\&=
\op_{y\in W;y\cdo\l\in\boc}\op_{y'\in W;y'{}\i\cdo y'(\l)\in\boc}
\Hom_{\cc^\boc\tcb^2}(\LL_{y'(\l)}^{\dy'{}\i},
\LL_{y(\l)}^{\dy\i}\un{\cir}\LL_{\ee^s(\l')}^{\ee^s(\du)})\\&
\ot M_{\ee^{-s}(\l),\l,y'(\l)}^{\ee^{-s}(\dy),\dw,\dy'{}\i,[1,3]}\la|y|+|w|+|y'|+4\r+\nu\ra.tag a
\endalign$$

\subhead 4.19\endsubhead
In the setup of 4.18 we shall apply the method of \cite{\CONV, 1.12} with 
$\Ph:\cd_m(Y_1)@>>>\cd_m(Y_2)$ replaced by ${}''d_!:\cd_m(\tcb^5)@>>>\cd_m(\tcb^4)$ and with 
$\cd^{\preceq}(Y_1)$, $\cd^{\preceq}(Y_2)$ replaced by ${}_1\cd^{\preceq}(\tcb^5)$, 
${}_1\cd^{\preceq}(\tcb^4)$, see 4.15. We shall take $\XX$ in {\it loc.cit.} equal to 
$\Xi={}''\vt_!({}''K)$ as in 4.15, $(w_1,w_3)=(u,w)$, $(\l_1,\l_3)=(\l',\l)$. The conditions of 
{\it loc.cit.} are satisfied: those concerning $\XX$ are satisfied with $c'=k=|w|+|u|+3\nu+5\r$ 
(see 4.8(c)); those concerning $\Ph$ are satisfied with $c=a+\r$ (see 4.17). We see that
$$\align&\un{gr_{a+\r+k}(({}''d_!{}''\vt_!({}''K))^{a+\r+k})}((a+\r+k)/2)\\&
=\un{gr_{a+\r}(({}''d_!\un{gr_k(({}''\vt_!({}''K))^k)}(k/2))^{a+\r})}((a+\r)/2).\endalign$$
Using 4.11(a), we have:
$$\align&gr_k({}''\vt_!({}''K))^k)(k/2)=\op_{y'\in W}M_{\l',\ee^{-s}(\l),\l,y'(\l)}^{\du,
\ee^{-s}(\dy'),\dw,\dy'{}\i,[1,4]}\la2|y'|+|w|+|u|+5\r+\nu\ra\\&
=\un{gr_k({}''\vt_!({}''K))^k}(k/2).\endalign$$
Hence, using 4.17(a), we have
$$\align&\un{gr_{a+\r}(({}''d_!\un{gr_k(({}''\vt_!({}''K))^k)}(k/2))^{a+\r})}
((a+\r)/2)\\&=\op_{y'\in W}\op_{y_1\in W;y_1\cdo\ee^{-s}(\l)\in\boc}
\Hom_{\cc^\boc\tcb^2}(\LL_{\ee^{-s}(\l)}^{\dy_1},\LL_{\l'}^{\du}\un{\cir}
\LL_{\ee^{-s}(\l)}^{\ee^{-s}(\dy')})\\&
\ot M_{\ee^{-s}(\l),\l,y'(\l)}^{\dy_1,\dw,\dy'{}\i,[1,3]}\la|y_1|+|w|+|y'|+4\r+\nu\ra.\endalign$$
Since $u\cdo\l'\in\boc$, for $y'\in W$ we have 
$$\Hom_{\cc^\boc\tcb^2}(\LL_{(\ee^{-s}(\l)}^{\dy_1},\LL_{\l'}^{\du}\un{\cir}
\LL_{\ee^{-s}(\l)}^{\ee^{-s}(\dy')})=0$$
unless $\ee^s(\l')=y'(\l)$ (see \cite{\MONO, 4.6(b)}) and $y'(\l)=\ee^s(\l')$ (see 
\cite{\MONO, 2.26(a)}). Thus we have
$$\align&\un{gr_{a+\r+k}(({}''d_!{}''\vt_!({}''K))^{a+\r+k})}((a+\r+k)/2)\\&=
\op_{y'\in W;y'\cdo\l\in\boc}\op_{y_1\in W;y_1\cdo \ee^{-s}(\l)\in\boc}
\Hom_{\cc^\boc\tcb^2}(\LL_{\ee^{-s}(\l)}^{\dy_1},\LL_{\l'}^{\du}\un{\cir}
\LL_{\ee^{-s}(\l)}^{\ee^{-s}(\dy')})\\&
\ot M_{\ee^{-s}(\l),\l,y'(\l)}^{\dy_1,\dw,\dy'{}\i,[1,3]}\la|y_1|+|w|+|y'|+4\r+\nu\ra.\endalign$$
Setting $y_1=\ee^{-s}y$ and using that $\ee^{-s}y\cdo\ee^{-s}(\l)\in\boc$ if and only if
$y\cdo\l\in\boc$, we can rewrite this as follows:
$$\align&\un{gr_{a+\r+k}(({}''d_!{}''\vt_!({}''K))^{a+\r+k})}((a+\r+k)/2)\\&=
\op_{y'\in W;y'\cdo\l\in\boc}\op_{y\in W;y\cdo\l\in\boc}
\Hom_{\cc^\boc\tcb^2}(\LL_{\ee^{-s}(\l)}^{\ee^{-s}\dy},\LL_{\l'}^{\du}\un{\cir}
\LL_{\ee^{-s}(\l)}^{\ee^{-s}(\dy')})\\&
\ot M_{\ee^{-s}(\l),\l,y'(\l)}^{\ee^{-s}\dy,\dw,\dy'{}\i,[1,3]}\la|y|+|w|+|y'|+4\r+\nu\ra.\tag a
\endalign$$

\subhead 4.20\endsubhead
Let $y_1\cdo\l_1\in\boc$, $y_2\cdo\l_2\in\boc$, $y_3\cdo\l_3\in\boc$. From \cite{\MONO, 3.20} we
see that:

(a) {\it we have canonically}
$$\Hom_{\cc^\boc\tcb^2}(\LL_{y_2(\l_2)}^{\dy_2\i},
\LL_{y_1(\l_1)}^{\dy_1\i}\un{\cir}\LL_{\l_3}^{\dy_3})
=\Hom_{\cc^\boc\tcb^2}(\LL_{\l_1}^{\dy_1},\LL_{\l_3}^{\dy_3}\un{\cir}\LL_{\l_2}^{\dy_2}).$$

In the setup of 4.18, we apply 4.18(a), 4.19(a) to $w\cdo\l$, $u\cdo\l'$ and we use the equality
$$\align&\un{gr_{a+\r+k}(({}'d_!{}'\vt_!({}'K))^{a+\r+k})}((a+\r+k)/2)\\&=
\un{gr_{a+\r+k}(({}''d_!{}''\vt_!({}''K))^{a+\r+k})}((a+\r+k)/2)\endalign$$
which comes from ${}'d_!{}'\vt_!({}'K)={}''d_!{}''\vt_!({}''K)$, see 4.15(a); we obtain
$$\align&\op_{y\in W;y\cdo\l\in\boc}\op_{y'\in W;y'\cdo \l\in\boc}
\Hom_{\cc^\boc\tcb^2}(\LL_{y'(\l)}^{\dy'{}\i},
\LL_{y(\l)}^{\dy\i}\un{\cir}\LL_{\ee^s(\l')}^{\ee^s(\du)})\\&
\ot M_{\ee^{-s}(\l),\l,y'(\l)}^{\ee^{-s}(\dy),\dw,\dy'{}\i,[1,3]}\la|y|+|w|+|y'|+4\r+\nu\ra\\&
=\op_{y'\in W;y'\cdo\l\in\boc}\op_{y\in W;y\cdo \l\in\boc}
\Hom_{\cc^\boc\tcb^2}(\LL_{\ee^{-s}(\l)}^{\ee^{-s}\dy},\LL_{\l'}^{\du}\un{\cir}
\LL_{\ee^{-s}(\l)}^{\ee^{-s}(\dy')})\\&
\ot M_{\ee^{-s}(\l),\l,y'(\l)}^{\ee^{-s}\dy,\dw,\dy'{}\i,[1,3]}\la|y|+|w|+|y'|+4\r+\nu\ra.\tag b
\endalign$$

\subhead 4.21\endsubhead
We assume that $w\cdo\l,u\cdo\l'$ in 4.18 satisfy in addition $w\cdo\l\in\boc$. We apply $p_{03!}$ 
and $\la N\ra$ for some $N$ to the two sides of 4.20(b). (Recall that $p_{03}:\tcb^4@>>>\tcb^2$.) 
We obtain
$$\align&\op_{y\in W;y\cdo\l\in\boc}\op_{y'\in W;y'\cdo\l\in\boc}
\Hom_{\cc^\boc\tcb^2}(\LL_{y'(\l)}^{\dy'{}\i},
\LL_{y(\l)}^{\dy\i}\un{\cir}\LL_{\ee^s(\l')}^{\ee^s(\du)})\\&
\ot \LL_{\ee^{-s}(\l)}^{\ee^{-s}(\dy)}\cir\LL_\l^{\dw}\cir\LL_{y'(\l)}^{\dy'{}\i}\\&
=\op_{y\in W;y\cdo \l\in\boc}\op_{y'\in W;y'\cdo\l\in\boc}
\Hom_{\cc^\boc\tcb^2}(\LL_{\ee^{-s}(\l)}^{\ee^{-s}\dy},\LL_{\l'}^{\du}\un{\cir}
\LL_{\ee^{-s}(\l)}^{\ee^{-s}(\dy')})\\&
\ot\LL_{\ee^{-s}(\l)}^{\ee^{-s}\dy}\cir\LL_{\l}^{\dw}\cir\LL_{y'(\l)}^{\dy'{}\i}.\endalign$$
Applying $\un{()^{\{2(a-\nu)\}}}$ to both sides and using \cite{\MONO, 2.24(a)} we obtain
$$\align&\op_{y\in W;y\cdo\l\in\boc}\op_{y'\in W;y'\cdo\l\in\boc}
\Hom_{\cc^\boc\tcb^2}(\LL_{y'(\l)}^{\dy'{}\i},
\LL_{y(\l)}^{\dy\i}\un{\cir}\LL_{\ee^s(\l')}^{\ee^s(\du)})\\&
\ot\LL_{\ee^{-s}(\l)}^{\ee^{-s}(\dy)}\un\cir\LL_\l^{\dw}\un\cir\LL_{y'(\l)}^{\dy'{}\i}\\&
=\op_{y\in W;y\cdo \l\in\boc}\op_{y'\in W;y'\cdo\l\in\boc}
\Hom_{\cc^\boc\tcb^2}(\LL_{\ee^{-s}(\l)}^{\ee^{-s}\dy},\LL_{\l'}^{\du}\un{\cir}
\LL_{\ee^{-s}(\l)}^{\ee^{-s}(\dy')})\\&
\ot\LL_{\ee^{-s}(\l)}^{\ee^{-s}\dy}\un\cir\LL_\l^{\dw}\un\cir\LL_{y'(\l)}^{\dy'{}\i},\endalign$$
or equivalently
$$\align&\op_{y\in W;y\cdo\l\in\boc}\LL_{\ee^{-s}(\l)}^{\ee^{-s}(\dy)}\un\cir\LL_\l^{\dw}\un\cir
\LL_{y(\l)}^{\dy\i}\un{\cir}\LL_{\ee^s(\l')}^{\ee^s(\du)}\\&
=\op_{y'\in W;y'\cdo\l\in\boc}
\LL_{\l'}^{\du}\un{\cir}\LL_{\ee^{-s}(\l)}^{\ee^{-s}(\dy')}\un\cir\LL_\l^{\dw}\un\cir
\LL_{y'(\l)}^{\dy'{}\i}.\endalign$$
Using 4.13(d), this can be rewritten as follows:
$$\un\fb(\Bbb L_{\l,s}^{\dw})\un{\cir}\LL_{\ee^s(\l')}^{\ee^s(\du)}
=\LL_{\l'}^{\du}\un{\cir}\un\fb(\Bbb L_{\l,s}^{\dw}).\tag a.$$
Another identification of the two sides in (a) is given by 4.14(j) with $L=\Bbb L_{\l,s}^{\dw}$, 
$L'=\LL_{\l'}^{\du}$ (note that $\un\fb(L)=\un\fb'(L)$ by 4.14(d)). In fact, 
the arguments in 4.13-4.20 and in this subsection show that

(b) {\it these two identifications of the two sides of (a) coincide.}

\subhead 4.22\endsubhead
Let $s',s''\in\ZZ$. Let 
$$\align&V=\{(B_0,B_1,B_2,\g U_{B_0},\g'U_{B_1});\\&(B_0,B_1,B_2)\in\cb^3,\g\in\tG_{s'},
\g'\in\tG_{s''},\g B_0\g\i=B_1,\g'B_1\g'{}\i=B_2\}.\endalign$$
Define $p_{01}:V@>>>Z_{s'}$, $p_{12}:V@>>>Z_{s''}$, $p_{02}:V@>>>Z_{s'+s''}$ by 
$$p_{01}:(B_0,B_1,B_2,\g U_{B_0},\g'U_{B_1})\m(B_0,B_1,\g U_{B_0}),$$
$$p_{12}:(B_0,B_1,B_2, g U_{B_0},\g'U_{B_1})\m(B_1,B_2,\g'U_{B_1}),$$
$$p_{02}:(B_0,B_1,B_2,\g U_{B_0},\g'U_{B_1})\m(B_0,B_2,\g'\g U_{B_0}).$$
For $L\in\cd(Z_{s'}),L'\in\cd(Z_{s''})$ we set 
$$L\bul L'=p_{02!}(p_{01}^*L\ot p_{12}^*L')\in\cd(Z_{s'+s''}).$$ 
This operation defines a monoidal structure on $\sqc_{s'\in\ZZ}\cd(Z_{s'})$.
Hence if ${}^1L\in\cd(Z_{s_1},{}^2L\in\cd(Z_{s_2}),\do,
{}^rL\in\cd(Z_{s_r})$, then ${}^1L\bul{}^2L\bul\do\bul{}^rL\in\cd(Z_{s_1+\do+s_r})$ is well 
defined. Note that, if $L\in\cd_m(Z_{s'}),L'_m\in\cd(Z_{s''})$ then we have naturally
$L\bul L'\in\cd_m(Z_{s'+s''})$. We show: 

(a) {\it For $L\in\cd(Z_{s'}),L'\in\cd(Z_{s''})$ we have canonically 
$\e_{s'+s''}^*(L\bul L')=\e_{s'}^*(L)\cir\e^*_{s''}(L')$.}
\nl
Let 
$$Y=\{(x\UU,y\UU,\g U_{x\BB x\i});x\UU\in\tcb,y\UU\in\tcb;\g\in\tG_{s'}\}.$$
Define $j:Y@>>>\tcb^2$, $j_1:Y@>>>Z_{s'}$, $j_2:Y@>>>Z_{s''}$ by
$$\align &j(x\UU,y\UU,\g U_{x\BB x\i})=(x\UU,y\UU),\\&
j_1(x\UU,y\UU,\g U_{x\BB x\i})=(x\BB x\i,\g x\BB x\i\g\i,\g U_{x\BB x\i}),\\&
j_2(x\UU,y\UU,\g U_{x\BB x\i})=(\g x\BB x\i\g\i,y\BB y\i,y\UU\t^{s'+s''}x\i\g\i).\endalign$$
From the definitions we have 
$$\e_{s'+s''}^*(L\bul L')=j_!(j_1^*(L)\ot j_2^*(L'))=\e_{s'}^*(L)\cir\e^*_{s''}(L')$$
and (a) follows.

\subhead 4.23\endsubhead
Let $s'\in\ZZ_\boc$.  Let $L\in\cd^\spa Z_s,L'\in\cd^\spa Z_{s'}$. We show:

(a) {\it If $L\in\cd^{\preceq}Z_s$ or $L'\in\cd^{\preceq}Z_{s'}$ then 
$L\bul L'\in\cd^{\preceq}Z_{s+s'}$. If $L\in\cd^{\prec}Z_s$ or $L'\in\cd^{\prec}Z_{s'}$ then
$L\bul L'\in\cd^{\prec}Z_{s+s'}$.} 
\nl
For the first assertion of (a) we can assume that $L=\Bbb L_{\l,s}^{\dw}$,
$L'=\Bbb L_{\l',s'}^{\dw'}$ with $w\cdo\l\in I^s_n,w'\cdo\l'\in I^{s'}_n$ and either 
$w\cdo\l\preceq\boc$ or $w'\cdo\l'\preceq\boc$. Assume that $w_1\cdo\l_1\in I^{s+s'}_n$ and 
$\Bbb L_{\l_1,s+s'}^{\dw_1}$ is a composition factor of $(L\bul L')^j$. Then 
$\LL_{\l_1}^{\dw_1}=\ti\e_{s+s'}\Bbb L_{\l_1,s+s'}^{\dw_1}$ is a composition factor of
$$\align&\e^*_{s+s'}(L\bul L')^j\la\r\ra=(\e^*_{s+s'}(L\bul L'))^{j+\r}(\r/2)
=(\e^*_sL\cir\e^*_{s'}L')^{j+\r}(\r/2)\\&
=(\e^*_sL\la\r\ra\cir\e^*_{s'}L'\la\r\ra)^{j-\r}(-\r/2)
=(\LL_\l^{\dw}\cir\LL_{\l'}^{\dw'})^{j-\r}(\r/2).\endalign$$
From \cite{\MONO, 2.23(b)} we see that $w_1\cdo\l_1\preceq\boc$. This proves the first
assertion of (a). The second assertion of (a) can be reduced to the first assertion.

We show:

(b) {\it Assume that $L\in\cm^\spa Z_s,L'\in\cm^\spa Z_{s'}$ and that either 
$L\in\cd^{\preceq}Z_s$ or $L'\in\cd^{\preceq}Z_{s'}$. If $j>a+\r-\nu$ then 
$(L\bul L')^j\in\cm^{\prec}Z_{s+s'}$.}
\nl
We can assume that $L=\Bbb L_{\l,s}^{\dw}$, $L'=\Bbb L_{\l',s'}^{\dw'}$ with 
$w\cdo\l\in I^s_n,w'\cdo\l'\in I^{s'}_n$ and either $w\cdo\l\in\boc$ or $w'\cdo\l'\in\boc$. Assume
that $w_1\cdo\l_1\in I^{s+s'}_n$ and that $\Bbb L_{\l_1,s+s'}^{\dw_1}$ is a composition factor of 
$(L\bul L')^j$. Then as in the proof of (a), $\LL_{\l_1}^{\dw_1}$ is a composition factor of 
$$\ti e_{s+s'}(L\bul L')^j=(\LL_\l^{\dw}\cir\LL_{\l'}^{\dw'})^{j-\r}(-\r/2).$$
Since $j-\r>a-\nu$ we see from \cite{\MONO, 2.23(a)} that $w_1\cdo\l_1\prec\boc$. This proves (b).

\subhead 4.24\endsubhead
Let $s'\in\ZZ_\boc$. For $L\in\cc^\boc_0Z_s,L'\in\cc^\boc_0Z_{s'}$ we set
$$L\un{\bul}L'=\un{(L\bul L')^{\{a+\r-\nu\}}}\in\cc^\boc_0Z_{s+s'}.$$
Using 4.23(a),(b) we see as in \cite{\MONO, 2.24} that for 
$L\in\cc^\boc_0Z_s,L'\in\cc^\boc_0Z_{s'},L''\in\cc^\boc_0Z_{s''}$ we have
$$L\un{\bul}(L'\un{\bul}L'')=(L\un{\bul}L')\un{\bul}L''=\un{(L\bul L'\bul L'')^{\{2a+2\r-2\nu\}}}.$$
We see that $L,L'\m L\un{\bul}L'$ defines a monoidal structure on 
$\sqc_{s'\in\ZZ_\boc}\cc^\boc_0Z_{s'}$. Hence if 
${}^1L\in\cc^\boc_0Z_{s_1},{}^2L\in\cc^\boc_0Z_{s_2},\do,{}^rL\in\cc^\boc_0Z_{s_r}$,
then ${}^1L\un{\bul}{}^2L\un{\bul}\do\un{\bul}{}^rL\in\cc^\boc_0Z_{s_1+\do+s_r}$ is well defined; 
we have
$${}^1L\un{\bul}{}^2L\un{\bul}\do\un{\bul}{}^rL=
\un{({}^1L\bul{}^2L\bul\do\bul{}^rL)^{\{(r-1)(a+\r-\nu)\}}}.\tag a$$
For $L\in\cc^\boc_0Z_s,L'\in\cc^\boc_0Z_{s'}$ we have $\ti\e_s L,\ti\e_{s'}L'\in\cc^\boc_0\tcb^2$.
We show:
$$\ti\e_{s+s'}(L\un{\bul}L')=(\ti\e_sL)\un{\cir}(\ti\e_{s'}L').\tag b$$
It is enough to show that
$$\align&\e^*_{s+s'}(gr_0((L\bul L')^{a+\r-\nu})((a+\r-\nu)/2))[\r](\r/2)\\&
=gr_0((\e^*_sL[\r](\r/2)\cir\e^*_{s'}L'[\r](\r/2))^{a-\nu})((a-\nu)/2))).\endalign$$
The left hand side is equal to
$$gr_0(\e^*_{s+s'}((L\bul L')^{a+\r-\nu})((a+\r-\nu)/2))[\r](\r/2))$$
hence it is enough to show:
$$\align&\e^*_{s+s'}((L\bul L')^{a+\r-\nu})((a+\r-\nu)/2))[\r](\r/2)\\&
=(\e^*_sL[\r](\r/2)\cir\e^*_{s'}L'[\r](\r/2))^{a-\nu}((a-\nu)/2))\endalign$$
that is,
$$\e^*_{s+s'}((L\bul L')^{a+\r-\nu})[\r]=(\e^*_sL[\r]\cir\e^*_{s'}L'[\r])^{a-\nu},$$
or, after using 4.3(b):
$$(\e^*_{s+s'}(L\bul L'))^{a+2\r-\nu}=(\e^*_sL\cir\e^*_{s'}L')^{a+2\r-\nu}.$$
It remains to use that $\e^*_{s+s'}(L\bul L')=\e^*_sL\cir\e^*_{s'}L'$, see 4.22(a).

\subhead 4.25\endsubhead
In the setup of 4.14 let
$${}^\di\cy=\TT^2\bsl\{((x_0\UU,x_1\UU,x_2\UU,x_3\UU),\g)\in\tcb^4\T\tG_s;
\g\in x_3\UU\t^s x_0\i,\g\in x_2\UU\t^s x_1\i\}$$
where $\TT^2$ acts freely by
$$(t_1,t_2):((x_0\UU,x_1\UU,x_2\UU,x_3\UU),\g)\m((x_0t_1\UU,x_1t_2\UU,x_2t_2\UU,x_3t_1\UU),\g).$$
We define ${}^\di\eta:{}^\di\cy@>>>Z_s$ by
$$((x_0\UU,x_1\UU,x_2\UU,x_3\UU),\g)\m\e_s(x_1\UU,x_2\UU).$$
We define $d:{}^\di\cy@>>>Z_s$ by
$$((x_0\UU,x_1\UU,x_2\UU,x_3\UU),\g)\m\e_s(x_0\UU,x_3\UU).$$
We define $\fb'':\cd(Z_s)@>>>\cd(Z_s)$ and $\fb'':\cd_m(Z_s)@>>>\cd_m(Z_s)$ by
$$\fb''(L)=d_!({}^\di\et)^*L.$$
From the definitions it is clear that
$$\fb'(L)=\e_s^*\fb''(L).\tag a$$
Using (a) we see that 4.14(a),(b),(c) imply the following statements.

(b) {\it If $L\in\cd^{\preceq}(Z_s)$, then $\fb''(L)\in\cd^{\preceq}Z_s$. If 
$L\in\cd^{\prec}(Z_s)$ then $\fb''(L)\in\cd^{\prec}Z_s$.}

(c) {\it If $L\in\cm^{\preceq}(Z_s)$ and $h>2\nu+2a$ then 
$(\fb''(L))^h\in\cm^{\prec}\tcb^2$.}
\nl
We define $\un{\fb''}:\cc^\boc_0(Z_s)@>>>\cc^\boc_0(Z_s)$ by
$$\un{\fb''}(L)=\un{gr_{2\nu+2a}((\fb''(L))^{2\nu+2a})}(\nu+a).$$ 
Using results in 4.3 we see that, if $L\in\cc^\boc_0Z_s$, then

(d) $\un\fb'(L)=\ti\e_s(\un{\fb''}(L))$.

\head 5. The monoidal category $\cc^\boc\tcb^2$ \endhead
\subhead 5.1\endsubhead
{\it In this section, $\boc,a,\fo,n,\Ps$ are as in 3.1(a).}
\nl
Define $\d:\tcb@>>>\tcb^2$ by $x\UU\m(x\UU,x\UU)$. For $w\cdo\l\in\boc$ we set
$$\b_{w\cdo\l}=\ch^{-a+|w|}(\d^*(L_\l^{\dw\sha}))((-a+|w|)/2).$$
By \cite{\MONO, 4.1} we have

(a) $\dim\b_{w\cdo\l}=1$ if $w\cdo\l\in\DD_\boc$, $\dim\b_{w\cdo\l}=0$ if $w\cdo\l\n\DD_\boc$.
\nl
We set
$$\bold1'=\op_{d\cdo\l\in\DD_\boc}\b^*_{d\cdo\l}\ot\LL_\l^{\dot d}\in\cc^\boc_0\tcb^2.$$
Here $\b^*_{d\cdo\l}$ is the vector space dual to $\b_{d\cdo\l}$.

\subhead 5.2\endsubhead
For $L\in\cd_m(\tcb^2)$ we set $L^\da=\ti\fh^*L$ where $\ti\fh:\tcb^2@>>>\tcb^2$ is as in 3.1. By 
\cite{\MONO, 4.4(b)}, we have:

(a) {\it If $L\in\cc^\boc_0\tcb^2$ then $\fD(L^\da)\in\cc^\boc_0\tcb^2$.
If $L\in\cc^\boc\tcb^2$ then $\fD(L^\da)\in\cc^\boc\tcb^2$.}

\subhead 5.3\endsubhead
The bifunctor $\cc^\boc_0\tcb^2\T\cc^\boc_0\tcb^2@>>>\cc^\boc_0\tcb^2$, $L,L'\m L\un{\circ}L'$ in
3.10 gives rise to a bifunctor $\cc^\boc\tcb^2\T\cc^\boc\tcb^2@>>>\cc^\boc\tcb^2$ denoted again by 
$L,L'\m L\un{\cir}L'$ as follows. Let $L\in\cc^\boc\tcb^2$, $L'\in\cc^\boc\tcb^2$; by replacing
if necessary $\Ps$ by a power, we choose mixed 
structures of pure weight $0$ on $L,L'$, we define $L\un{\cir}L'$ as in 3.10 in terms of these 
mixed structures and we then disregard the mixed structure on $L\un{\cir}L'$. The resulting object
of $\cc^\boc\tcb^2$ is denoted again by $L\un{\cir}L'$; it is independent of the choice of $\Ps$
which defines the mixed structures.

Similarly for $s,s'$ in $\ZZ_\boc$,  the bifunctor 
$\cc^\boc_0Z_s\T\cc^\boc_0Z_{s'}@>>>\cc^\boc_0Z_{s+s'}$, $L,L'\m L\un{\bul}L'$ in 4.24 gives rise 
to a bifunctor $\cc^\boc Z_s\T\cc^\boc Z_{s'}@>>>\cc^\boc Z_{s+s'}$ denoted again by 
$L,L'\m L\un{\bul}L'$. Moreover, $\un\fb:\cc^\boc_0Z_s@>>>\cc^\boc_0\tcb^2$ in 4.13 can be also
viewed as a functor $\un\fb:\cc^\boc Z_s@>>>\cc^\boc\tcb^2$.

The operation $L\un{\bul}L'$ (resp. $L\un{\circ}L'$) makes 
$\sqc_{s\in\ZZ_\boc}\cc^\boc Z_s$ (resp. $\cc^\boc\tcb^2$) into a monoidal abelian category (see 
4.24, 3.10). By \cite{\MONO, 4.5(a)}, we have:

(a) {\it For $L,L'$ in $\cc^\boc\tcb^2$ we have canonically}
$$\Hom_{\cc^\boc\tcb^2}(\bold1',L\un\circ L')=\Hom_{\cc^\boc\tcb^2}(\fD(L'{}^\da),L).$$

\subhead 5.4\endsubhead
We set 
$$\bold1=\op_{d\cdo\l\in\DD_\boc}\b_{d\cdo\l}\ot\LL_\l^{\dot d\i}\in\cc^\boc_0\tcb^2.\tag a$$
Here $\b_{d\cdo\l}$ is as in 5.1. By \cite{\MONO, 4.7(g)},

(a) {\it $\bold1=\bold1'$ is a unit object of the monoidal category $\cc^\boc\tcb^2$.}
\nl
By \cite{\MONO, 4.8}, this monoidal category has a natural rigid structure.

\subhead 5.5\endsubhead
{\it In the remainder of this section we fix $s\in\ZZ_\boc$.}
\nl
In this case, $(\ee^s)^*$ defines an equivalence of categories $\cc^\boc\tcb^2@>>>\cc^\boc\tcb^2$;
this follows from 3.11(a).

By analogy with \cite{\URE, 6.2} and slightly extending a definition in \cite{\MUG, 3.1}, we 
define an {\it $\ee^s$-half-braiding} for an object $\cl\in\cc^\boc\tcb^2$, as a collection 
$e_\cl=\{e_\cl(L);L\in\cc^\boc\tcb^2\}$ where $e_\cl(L)$ is an isomorphism
$(\ee^s)^*(L)\un\cir\cl@>\si>>\cl\un\cir L$ such that $e_\cl(\bold1)=Id_\cl$ and such that 
(i),(ii) below hold:

(i) If $L@>t>>L'$ is a morphism in $\cc^\boc\tcb^2$ then the diagram
$$\CD
(\ee^s)^*(L)\un\cir\cl@>e_\cl(L)>>\cl\un\cir L\\
@V(\ee^s)^*(t)\un\bul1VV                  @V1\un\bul tVV   \\
(\ee^s)^*(L')\un\cir\cl@>e_\cl(L')>>\cl\un\cir L'\\
\endCD$$
is commutative.

(ii) If $L,L'\in\cc^\boc\tcb^2$ then $e_\cl(L\un\cir L'):(\ee^s)^*(L\un\cir L')\un\cir\cl@>>>
\cl\un\cir(L\un\cir L')$ is equal to the composition 
$$(\ee^s)^*(L)\un\cir(\ee^s)^*(L')\un\cir\cl@>1\un\cir e_\cl(L')>>
(\ee^s)^*(L)\un\cir\cl\un\cir L'@>e_\cl(L)\un\cir1>>\cl\un\cir L\un\cir L'.$$
(When $s=0$ this reduces to the definition of a half-braiding for $\cl$ given in
\cite{\MUG, 3.1}.)

Let $\cz^\boc_{\ee^s}$ be the category whose objects are the pairs $(\cl,e_\cl)$ where
$\cl$ is an object of $\cc^\boc\tcb^2$ and $e_\cl$ is an $\ee^s$-half-braiding for $\cl$. For
$(\cl,e_\cl),(\cl',e_{\cl'})$ in $\cz^\boc_{\ee^s}$ we define
$\Hom_{\cz^\boc_{\ee^s}}((\cl,e_\cl),(\cl',e_{\cl'}))$ to be the vector space consisting of all
$t\in\Hom_{\cc^\boc\tcb^2}(\cl,\cl')$ such that for any $L\in\cc^\boc\tcb^2$ the diagram
$$\CD
(\ee^s)^*(L)\un\cir\cl@>e_\cl(L)>>\cl\un\cir L\\
@V1\un\cir tVV                  @V t\un\cir1VV   \\
(\ee^s)^*(L)\un\cir\cl'@>e_\cl(L')>>\cl'\un\cir L\\
\endCD$$
is commutative. We say that $\cz^\boc_{\ee^s}$ is the {\it $\ee^s$-centre} of $\cc^\boc\tcb^2$.
By a variation of a result of \cite{\MUG}, \cite{\ENO} (which concerns the usual centre), the additive 
category $\cz^\boc_{\ee^s}$ is semisimple, with finitely many isomorphism classes of simple objects. 
By a variation of a general result on semisimple rigid monoidal categories in
\cite{\ENO, Proposition 5.4}, for any $L\in\cc^\boc\tcb^2$ one can define directly an
$\ee^s$-half-braiding on the object
$$\ci_s(L)=\op_{y\cdo\l\in\boc}(\ee^s)^*(\LL_\l^{\dy})\un\cir L\un\cir\LL_{y(\l)}^{\dy\i}
=\op_{y\cdo\l\in\boc}\LL_{\ee^{-s}(y)}^{\ee^{-s}(\l)}\un\cir L\un\cir\LL_{y(\l)}^{\dy\i}$$
of $\cc^\boc\tcb^2$ such that, denoting by $\ov{\ci_s(L)}$ the corresponding object of 
$\cz^\boc_{\ee^s}$, we have canonically
$$\Hom_{\cc^\boc\tcb^2}(L,L')=\Hom_{\cz^\boc_{\ee^s}}(\ov{\ci_s(L)},L')\tag a$$
for any $L'\in\cz^\boc_{\ee^s}$. 
(We use that for $y\cdo\l\in\boc$, the dual of the simple object $\LL_\l^{\dy}$ is 
$\LL_{y(\l)}^{\dy\i}$, see \cite{\MONO, 4.4(c)}; we also use 3.11(a).)
The $\ee^s$-half-braiding  on $\ci_s(L)$ can be described as follows: for any $X\in\cc^\boc\tcb^2$ 
we have canonically
$$\align&(\ee^s)^*(X)\un\cir\ci_s(L)=\op_{y\cdo\l\in\boc}
(\ee^s)^*(X)\un\cir(\ee^s)^*(\LL_{\l}^{\dy})\un\cir L\un\cir\LL_{y(\l)}^{\dy\i}\\&
=\op_{y\cdo\l\in\boc,z\cdo\l'\in\boc}\Hom_{\cc^\boc\tcb^2}
((\ee^s)^*(\LL_{\l'}^{\dz}),(\ee^s)^*(X\un\cir\LL_\l^{\dy}))\ot
(\ee^s)^*(\LL_{\l'}^{\dz})\un\cir L\un\cir\LL_{y(\l)}^{\dy\i}\\&
=\op_{y\cdo\l\in\boc,z\cdo\l'\in\boc}\Hom_{\cc^\boc\tcb^2}
(\LL_{\l'}^{\dz},X\un\cir\LL_\l^{\dy})\ot
(\ee^s)^*(\LL_{\l'}^{\dz})\un\cir L\un\cir\LL_{y(\l)}^{\dy\i}\\&
=\op_{y\cdo\l\in\boc,z\cdo\l'\in\boc}
\Hom_{\cc^\boc\tcb^2}(\LL_{y(\l)}^{\dy\i},\LL_{z(\l')}^{\dz\i}\ot X)
\ot(\ee^s)^*(\LL_{\l'}^{\dz})\un\cir L\un\cir\LL_{y(\l)}^{\dy\i}\\&
=\op_{z\cdo\l'\in\boc}(\ee^s)^*(\LL_{\l'}^{\dz})\un\cir L\un\cir \LL_{z(\l')}^{\dz\i}\ot X
=\ci_s(L)\un\cir X.
\endalign$$
(The fourth equality uses 4.20(a); we have also used 3.11(a).) We show:

(b) {\it If $z\cdo\l\in\boc$ and $\ci_s(\LL_\l^{\dz})\ne0$ then $z\cdo\l\in\boc^s$.}
\nl
For some $y\cdo\l'\in\boc$ we have $\LL_{\ee^{-s}(\l')}^{\ee^{-s}(\dy)}\un\cir\LL_\l^{\dz}\ne0$ 
(hence $\ee^{-s}(\l')=z(l)$) and $\LL_\l^{\dz}\un\cir\LL_{y(\l')}^{\dy\i}\ne0$ (hence $\l=\l'$). 
It follows that $z(\l)=\ee^{-s}(\l)$ and (b) is proved.

\subhead 5.6\endsubhead
By 4.13(d), for $z\cdo\l\in\boc^s$ we have canonically 
$$\un\fb(\Bbb L_{\l,s}^{\dz})=\ci_s(\LL_\l^{\dz})\tag a$$ 
as objects of $\cc^\boc\tcb^2$. Here $\un\fb:\cc^\boc Z_s@>>>\cc^\boc\tcb^2$ is as in 5.3.
Now $\ci_s(\LL_\l^{\dz})$ has a natural $\ee^s$-half-braiding (by 5.5) and 
$\un\fb(\Bbb L_{\l,s}^{\dz})$ has a natural $\ee^s$-half-braiding (by 4.14(j)). By 4.21(b),

(b) {\it these two $\ee^s$-half-braidings are compatible with the identification (a).}
\nl
In view of (a),(b) we can reformulate 5.5(a) as follows.

\proclaim{Theorem 5.7}For any $z\cdo\l\in\boc^s$, $L'\in\cz^\boc_{\ee^s}$, we have canonically
$$\Hom_{\cc^\boc\tcb^2}(\LL_\l^{\dz},L')=\Hom_{\cz^\boc_{\ee^s}}
(\ov{\un{\fb}(\Bbb L_{\l,s}^{\dz})},L')\tag a$$
where $\ov{\un{\fb}(\Bbb L_{\l,s}^{\dz})}$ is $\un{\fb}(\Bbb L_{\l,s}^{\dz})$ viewed as an object 
of $\cz^\boc_{\ee^s}$ with the $\ee^s$-half-braiding given by 4.14(j).
\endproclaim

\subhead 5.8\endsubhead
We set
$$\bold1'_0=\op_{d\cdo\l\in\DD_\boc}\b^*_{d\cdo\l}\ot\Bbb L_{\l,0}^{\dot d}\in\cc^\boc Z_0.$$
From the definitions we have $\ti\e_0\bold1'_0=\bold1'$. Since $\bold1'=\bold1$, we have
also $\ti\e_0\bold1'_0=\bold1$. We show:

(a) {\it For $L\in\cc^\boc Z_{-s},L'\in\cc^\boc Z_s$ we have}
$$\Hom_{\cm(Z_0)}(\bold1'_0,L\un\bul L')=\Hom_{\cm(Z_{-s})}(\fD(L'{}^\da),L).$$
We can assume that $L=\Bbb L_{\l,-s}^{\dw}$, $L'=\Bbb L_{\l',s}^{\dw'}$ where 
$w\cdo\l\in\boc^{-s}$, $w'\cdo\l'\in\boc^s$. Using the fully faithfulness of 
$\ti\e_0:\cm(Z_0)@>>>\cm\tcb^2$, $\ti\e_{-s}:\cm(Z_{-s})@>>>\cm\tcb^2$, and the equality 
$\ti\e_0\bold1'_0=\bold1$, we see that it is enough to prove that
$$\Hom_{\cm(\tcb^2)}(\bold1,\ti\e_0(L\un\bul L'))=
\Hom_{\cm\tcb^2}(\ti\e_{-s}(\fD(L'{}^\da)),\ti\e_{-s}(L)).$$
From 4.3 we have $\ti\e_{-s}(L)=\LL_\l^{\dw}$, $\ti\e_s(L')=\LL_{\l'}^{\dw'}$,
$\ti\e_{-s}(\Bbb L_{w'(\l'),-s}^{\dw'{}\i})=\LL_{w'(\l')}^{\dw'{}\i}$.

From 4.3(e) we have 
$$\ti\e_{-s}(\fD(L'{}^\da))=\ti\e_{-s}(\fD(\Bbb L_{w'(\l'{}\i),-s}^{\dw'{}\i}))=
\ti\e_{-s}(\Bbb L_{w'(\l'),-s}^{\dw'{}\i})= \LL_{w'(\l')}^{\dw'{}\i}.$$
(We have use that $\fD(\Bbb L_{w'(\l'{}\i),-s}^{\dw'{}\i})=\Bbb L_{w'(\l'),-s}^{\dw'{}\i}$ which
follows from \cite{\MONO, 4.4(a)}.) Using 4.24(b), we have 
$$\ti\e_0(L\un{\bul}L')=(\ti\e_{-s}L)\un{\cir}(\ti\e_sL')=
\LL_\l^{\dw}\un{\cir}\LL_{\l'}^{\dw'}.$$
Hence it is enough to prove
$$\Hom_{\cm\tcb^2}(\bold1,\LL_\l^{\dw}\un{\cir}\LL_{\l'}^{\dw'})=
\Hom_{\cm\tcb^2}(\LL_{w'(\l')}^{\dw'{}\i},\LL_\l^{\dw}).$$
This follows from \cite{\MONO, 4.5(a)}.

\head 6. Truncated induction, truncated restriction, truncated convolution \endhead
\subhead 6.1\endsubhead
{\it In this section  we fix $s\in\ZZ$. }
\nl
Let $\dZ_s=\{(B,B',\g)\in\cb\T\cb\T\tG_s;\g B\g\i=B'\}$. We have a diagram 
$$Z_s@<f<<\dZ_s@>\p>>\tG_s\tag a$$ 
where $f(B,B',\g)=(B,B',\g U_B)$, $\p(B,B',\g)=\g$. 
Note that $G$ acts on $Z_s$ by $g:(B,B',\g U_B)\m(gBg\i,gB'g\i,g\g g\i U_{gBg\i})$, on $\dZ_s$ by 
$g:(B,B',\g)\m(gBg\i,gB'g\i,g\g g\i)$, on $\tG_s$ by $g:\g\m g\g g\i$; moreover, $f$ and $\p$ are
compatible with these $G$-actions. We define $\c:\cd(Z_s)@>>>\cd(\tG_s)$ by
$$\c(L)=\p_!f^*L.$$
For any $w\cdo\l\in I$ we define $\fR_{\l,s}^{\dw}\in\cd(\tG_s)$, $R_{\l,s}^{\dw}\in\cd(\tG_s)$ by 
$$\fR_{\l,s}^{\dw}=\c(\cl_{\l,s}^{\dw}), R_{\l,s}^{\dw}=\c(\cl_{\l,s}^{\dw\sha}),
\text{ if }w\cdo\l\in I^s,$$    
$$\fR_\l^{\dw}=0, R_\l^{\dw}=0\text{ if }w\cdo\l\n I^s.$$

\mpb

Assume now that $s\ne0$ and that we are in case A. In this case, the conjugation $G$-action on $\tG_s$ is transitive, see
2.1, and the stabilizer of $\t^s$ for this $G$-action is the finite group 
$G^{\ee^s}=\{g\in G;\ee^s(g)=g\}$. 

With the notation of 4.1, for $w\in W$ we have isomorphisms
$$X^w_s@>\si>>\p\i(\t^s)\cap f\i(Z^w_s), \bX^w_s@>\si>>\p\i(\t^s)\cap f\i(\bZ^w_s)$$
given by $B\m(B,\ee^s(B),\t^s)$. Using this, and the transitivity of the $G$-action on $\tG_s$, we
see that for $w\cdo\l\in I^s$ and for $j\in\ZZ$, 
$(\fR_{\l,s}^{\dw})^j[-\D]$ (resp. $(R_{\l,s}^{\dw})^j[-\D]$) is the $G$-equivariant local system 
on $\tG_s$ whose stalk at $\t^s$ is $H^{j-\D}_c(X_s^z,\cf_{\l,s}^{\dw})[\D]$ (resp. 
$IH^{j-\D}(\bX_s^z,\cf_{\l,s}^{\dw})[\D]$) with the $G^{\ee^s}$-action considered in 4.1.

\mpb

We return to the general case.
We say that a simple perverse sheaf $A$ on $\tG_s$ is a {\it character sheaf} if the following 
equivalent conditions are satisfied:

(i) there exists $w\cdo\l\in I$ such that $(A:\op_j(\fR_{\l,s}^{\dw})^j)\ne0$;

(ii) there exists $w\cdo\l\in I$ such that $(A:(R_{\l,s}^{\dw})^j)\ne0$.
\nl
In case A with $s\ne0$, if $A$ satisfies either (i) or (ii), then
 it must be $G$-equivariant, hence $A[-D]$
must be a $G$-equivariant local system whose stalk at $\t^s$ viewed as a $G^{\ee^s}$-module
is irreducible, so that in this case the equivalence of (i),(ii) follows from the equivalence of
(i),(ii) in 4.1. In case A with $s=0$ the equivalence of (i),(ii) follows from 
\cite{\CSIII, 12.7}; a similar proof applies in case B (see also \cite{\CDGVI, 28.13}). 

A character sheaf $A$ determines a $W$-orbit $\fo$ on $\fs_\iy$: the set of $\l\in\fs_\iy$ such 
that $(A:\op_j(\fR_{\l,s}^{\dw})^j)\ne0$ for some $w\in W$ (or equivalently 
$(A:\op_j(R_{\l,s}^{\dw})^j)\ne0$ for some $w\in W$); we have necessarily
$\ee^s(\fo)=\fo$. 
In case A with $s\ne0$ this follows from 4.1. In case A with $s=0$ this follows from 
\cite{\CSIII, 11.2(a), 12.7}; a similar proof applies in case B.

We now fix $\fo\in W\bsl\fs_\iy$ such that $\ee^s(\fo)=\fo$.
We say that $A$ is an {\it $\fo$-character sheaf} if the $W$-orbit on 
$\fs_\iy$ determined by $A$ is $\fo$. Let $CS_{\fo,s}$ be a set of representatives for the 
isomorphism classes of $\fo$-character sheaves on $\tG_s$.
In case A with $s\ne0$ we have a natural bijection $CS_{\fo,s}\lra\Irr_{\fo}(G^{\ee^s})$ 
(notation of 4.1);
to $A\in CS_{\fo,s}$ corresponds the stalk of the $G$-equivariant local system $A[-\D]$ at $\t^s$,
viewed as an irreducible $G^{\ee^s}$-module.

\mpb

Let $\fo\in W\bsl\fs_\iy$ be such that $\ee^s(\fo)=\fo$. With notation in 2.4 we have the 
following result.

(b) {\it There exists a pairing $CS_{\fo,s}\T\Irr_s(\HH^1_{\fo})@>>>\bbq$, $(A,E)\m b_{A,E}$ such 
that for any $A\in CS_{\fo,s}$, any $z\cdo\l\in I$ with $\l\in\fo$ and any $j\in\ZZ$ we have}
$$(A:(R_{\l,s}^{\dz})^j)=(-1)^{j+\D}(j-\D-|z|;\sum_{E\in\Irr_s(\HH^1_{\fo})}b_{A,E}
\tr(\ee_sc_{z\cdo\l},E^v)).$$
Assume first that $z\cdo\l\in I^s$. 
In case A with $s\ne0$, (b) follows from 4.1(b). In case A with $s=0$, (b) is a 
reformulation of \cite{\CSIII, 14.11}, see \cite{\MONO, 5.1}. In case B, (b) can be deduced
from \cite{\CDGVII, 34.19} and the quasi-rationality result \cite{\CDGVIII, 39.8}. (In
{\it loc.cit.} there is the assumption that the adjoint group of $G$ is simple, which was 
made to simplify the arguments.)

Next we assume tha $z\cdo\l\in I-I^s$. Then the left hand side of (a) is zero; hence it
is enough to show that $\tr(\ee_sc_{z\cdo\l},E^v)=0$ for any $E\in\Irr_s(\HH^1_\fo)$. We have a 
direct sum decomposition $E^v=\op_{\l'\in\fs_\iy}1_{\l'}E^v$. It is enough to show that
for $\l'\in\fs_\iy$ we have $\ee_sc_{z\cdo\l}(1_{\l'}E^v)\sub 1_{\l''}E^v$ where $\l''\in\fs_\iy$,
$\l''\ne\l'$. We can assume that $\l'=\l$. We have 
$$\ee_sc_{z\cdo\l}(1_\l E^v)\sub\ee_s(1_{z(\l)}E^v)=1_{\ee^s(z(\l)}E^v.$$
It is enough to show that $\ee^s(z(\l))\ne\l$ that is, $z(\l)\ne\ee^{-s}(\l)$; this follows from 
$z\cdo\l\n I^s$.

\mpb

Given $A\in CS_{\fo,s}$, there is a unique two-sided cell $\boc_A$ of $I$ such that $b_{A,E}=0$ 
whenever $E\in\Irr_s(\HH^1_\fo)$ satisfies $\boc_E\ne\boc_A$. 
In case A with $s\ne0$ this follows from results in \cite{\ORA}, under the assumption
that the centre of $G$ is connected; but the argument in \cite{\ORA} extends to the general 
case. In case A with $s=0$ this follows from \cite{\CSIII, 16.7}. 
In case B this follows from \cite{\CDGIX, \S41}. We have necessarily $\boc_A\sub I_\fo$. As 
in \cite{\CDGIX, 41.8}, \cite{\CDGX, 44.18}, we see that:

(c) {\it We have $(A:\op_j(R_{\l,s}^{\dz})^j)\ne0$ for some $z\cdo\l\in\boc_A$; conversely, 
if 
$z\cdo\l\in I$ is such that $(A:\op_j(R_{\l,s}^{\dz})^j)\ne0$, then $\boc_A\preceq z\cdo\l$.}
\nl
Let $a_A$ be the value of the $a$-function on $\boc_A$. If $z\cdo\l\in I^s$, 
$E\in\Irr_s(\HH^1_\fo)$ satisfy $\tr(\ee_sc_{z\cdo\l},E^v)\ne0$ then 
$\boc_E\preceq z\cdo\l$; if in addition we have $z\cdo\l\in\boc_E$ then from the definitions we 
have
$$\tr(\ee_sc_{z\cdo\l},E^v)=\sum_{h\ge0}c_{z\cdo\l,E,h,s}v^{a_E-h}$$
where $c_{z\cdo\l,E,h,s}\in\bbq$ is zero for large $h$, 
$c_{z\cdo\l,E,0,s}=\tr(\ee_st_{z\cdo\l},E^\iy)$ and $a_E$ is as in 1.13. Hence from (b) we see 
that for $A\in CS_{\fo,s}$ and $z\cdo\l\in I_\fo$, $j\in\ZZ$, the following holds:

(d) {\it We have $(A:(R_{\l,s}^{\dz})^j)=0$ unless $\boc_A\preceq z\cdo\l$; if $z\cdo\l\in\boc_A$,
then
$$(A:(R_{\l,s}^{\dz})^j)=(-1)^{j+\D}(j-\D-|z|;
\sum_{E\in\Irr_s(\HH^1_\fo);\boc_E=\boc_A;h\ge0}b_{A,E}c_{z\cdo\l,E,h,s}v^{a_A-h})$$
which is $0$ unless $j-\D-|z|\le a_A$.}

\mpb

{\it In the remainder of this section let $\boc,a,n,\Ps$ be as in 3.1(a). We assume that $w\cdo\l\in\boc\implies\l\in\fo$.}
\nl
Note that $\c$ can be also viewed as a functor $\c:\cd_m(Z_s)@>>>\cd_m(\tG_s)$.

Let $\cm^{\preceq}\tG_s$  (resp. $\cm^{\prec}\tG_s$) be the category of perverse sheaves on 
$\tG_s$ whose composition factors are all of the form $A\in CS_{\fo,s}$ with $\boc_A\preceq\boc$ 
(resp. $\boc_A\prec\boc$). Let $\cd^{\preceq}\tG_s$  (resp. $\cd^{\prec}\tG_s$) be the subcategory
of $\cd(\tG_s)$ whose objects are complexes $K$ such that $K^j$ is in $\cm^{\preceq}\tG_s$ (resp. 
$\cm^{\prec}\tG_s$) for any $j$. Let $\cd^{\preceq}_m\tG_s$ (resp. $\cd^{\prec}_m\tG_s$) 
be the subcategory of $\cd_m(\tG_s)$ whose objects are also in $\cd^{\preceq}\tG_s$ 
(resp. $\cd^{\prec}\tG_s$). 

Let $z\cdo\l\in I_\fo$. From (d) we deduce:

(e) {\it If $z\cdo\l\preceq\boc$, then $(R_{\l,s}^{\dz})^j\in\cm^{\preceq}\tG_s$ for all 
$j\in\ZZ$.}

(f) {\it If $z\cdo\l\in\boc$ and $j>a+\D+|z|$ then $(R_{\l,s}^{\dz})^j\in\cm^{\prec}\tG_s$.}

(g) {\it If $z\cdo\l\prec\boc$ then $(R_{\l,s}^{\dz})^j\in\cm^{\prec}\tG_s$ for all $j\in\ZZ$.}

\subhead 6.2\endsubhead
Let $CS_{\boc,s}=\{A\in CS_{\fo,s};\boc_A=\boc\}$. For any $z\cdo\l\in I$ we set 
$$n_z=a(z)+\D+|z|.$$
Let $A\in CS_{\boc,s}$ and let $z\cdo\l\in\boc$. We have
$$(A:(R_{\l,s}^{\dz})^{n_z})
=(-1)^{a+|z|}\sum_{E\in\Irr_s(\HH^1_\fo)}b_{A,E}\tr(\ee_st_{z\cdo\l},E^\iy).\tag a$$
Indeed, from 6.1(b) we have
$$(A:(R_{\l,s}^{\dz})^{n_z})=(-1)^{a+|z|}\sum_{E\in\Irr_s(\HH^1_\fo)}b_{A,E}
(a;\tr(\ee_sc_{z\cdo\l},E^v))$$
and it remains to use that $(a;\tr(\ee_sc_{z\cdo\l},E^v))=\tr(\ee_st_{z\cdo\l},E^\iy)$. We show:

(b) {\it For any $A\in CS_{\boc,s}$ there exists $E\in\Irr_s(\HH^1_\fo)$ such that $b_{A,E}\ne0$
hence $\boc_E=\boc$.}
\nl
Assume that this is not so. Then, using 6.1(b), for any $z\cdo\l\in I_\fo$ we have \lb
$(A:\op_j(R_{\l,s}^{\dz})^j)=0$. This contradicts the assumption that $A\in CS_{\fo,s}$. We show:

(c) {\it For any $A\in CS_{\boc,s}$ there exists $z\cdo\l\in\boc$ such that 
$(A:(R_{\l,s}^{\dz})^{n_z})\ne0$.}
\nl
Assume that this is not so. Then, using (a), we see that
$$\sum_{E\in\Irr_s(\HH^1_\fo);\boc_E=\boc}b_{A,E}\tr(\ee_st_{z\cdo\l},E^\iy)=0$$
for any $z\cdo\l\in\boc$. 
If $z\cdo\l\in I_\fo-\boc$ then the last sum is automatically zero
since $t_{z\cdo\l}$ acts as $0$ on $E^\iy$ for each $E$ in the sum. Thus we have
$$\sum_{E\in\Irr_s(\HH^1_\fo);\boc_E=\boc}b_{A,E}\tr(\ee_st_{z\cdo\l},E^\iy)=0$$
for any $z\cdo\l\in I_\fo$.
In the last sum the condition $\boc_E=\boc$ is automatically satisfied if $b_{A,E}\ne0$. Thus we 
have
$$\sum_{E\in\Irr_s(\HH^1_\fo)}b_{A,E}\tr(\ee_st_{z\cdo\l},E^\iy)=0$$
for any $z\cdo\l\in I_\fo$. 
By a general argument (see for example \cite{\CDGVII, 34.14(e)}), the linear functions 
$t_{z\cdo\l}\m\tr(\ee_st_{z\cdo\l},E^\iy)$, $\JJ_\fo@>>>\bbq$ (for various 
$E$ as in the last sum) are linearly independent. It follows that $b_{A,E}=0$ for each $E$ as in
the last sum. This contradicts (b).

We show:

(d) {\it Let $z\cdo\l\in\boc$ be such that $(R_{\l,s}^{\dz})^{n_z}\ne0$. Then 
$z\cdo\l\underset\text{left}\to\si\ee^s(z\i)\cdo\ee^s(z(\l))$ and
$z\cdo\l\underset\text{left}\to\si\ee^s(z\i)\cdo\l$.}
\nl
Using (a) we see that there exists $E\in\Irr_s(\HH^1_\fo)$ such that 
$\tr(\ee_st_{z\cdo\l},E^\iy)\ne0$. We have $E^\iy=\op_{d\cdo\l_1\in\DD\cap\fo}t_{d\cdo\l_1}E^\iy$.
We define $d\cdo\l_1\in\DD\cap\fo$ by the condition that 
$z\cdo\l\underset\text{left}\to\si d\cdo\l_1$. We define $d'\cdo\l'_1\in\DD\cap\fo$ by the 
condition that $z\i\cdo z(\l)\underset\text{left}\to\si d'\cdo\l'_1$. Now 
$t_{z\cdo\l}:E^\iy@>>>E^\iy$ maps the summand
$t_{d\cdo\l_1}E^\iy$ into the summand $t_{d'\cdo\l'_1}E^\iy$ and all other summands to zero.
Moreover, $\ee_s$ maps $t_{d'\cdo\l'_1}E^\iy$ into $t_{\ee^s(d')\cdo\ee^s(\l'_1)}E^\iy$. Hence 
$\ee_st_{z\cdo\l}:E^\iy@>>>E^\iy$ maps the summand
$t_{d\cdo\l_1}E^\iy$ into the summand $t_{\ee^s(d')\cdo\ee^s(\l'_1)}E^\iy$ and all other
summands to zero. Since $\tr(\ee_st_{z\cdo\l},E^\iy)\ne0$ it follows that\lb
$t_{d\cdo\l_1}E^\iy=t_{\ee^s(d')\cdo\ee^s(\l'_1)}E^\iy\ne0$. Since 
$\ee^s(d')\cdo\ee^s(\l'_1)\in\DD\cap\fo$, it follows that\lb $d\cdo\l_1=\ee^s(d')\cdo\ee^s(\l'_1)$.
Since $\ee^s(z\i)\cdo \ee^s(z(\l))\underset\text{left}\to\si \ee^s(d')\cdo\ee^s(\l'_1)$, we see 
that\lb $z\cdo\l\underset\text{left}\to\si\ee^s(z\i)\cdo\ee^s(z(\l)$. To complete the proof, it 
remains to note that $\ee^s(z(\l))=\l$ that is $z\cdo\l\in I^s$. This follows from the fact that 
$(R_{\l,s}^{\dz})^{n_z}\ne0$. 

We show:

(e) {\it If $CS_{\boc,s}\ne\emp$ then $\ee^s(\boc)=\boc$.}
\nl
Using (c) and the hypothesis we see that there exists $z\cdo\l\in\boc$ such that 
$(R_{\l,s}^{\dz})^{n_z}\ne0$. Using (d), we see that $\ee^s(z\i)\cdo\ee^s(z(\l))\in\boc$. Since 
$z\i\cdo z(\l)\in\boc$ (see Q10 in 1.9) we have also $\ee^s(z\i)\cdo\ee^s(z(\l))\in\ee^s(\boc)$. 
Thus, $\boc\cap\ee^s(\boc)\ne\emp$. It follows that $\ee^s(\boc)=\boc$.

\subhead 6.3\endsubhead
{\it Until the end of 6.7 we assume that $s\in\ZZ_\boc$.}
\nl
We show:

(a) {\it If $L\in\cd^{\preceq}Z_s$ then $\c(L)\in\cd^{\preceq}\tG_s$. If 
$L\in\cd^{\prec}Z_s$ then $\c(L)\in\cd^{\prec}\tG_s$.}

(b)  {\it If $L\in\cm^{\preceq}Z_s$ and $j>a+\nu$ then 
$(\c(L))^j\in\cm^{\prec}\tG_s$.}
\nl
It is enough to prove (a),(b) assuming in addition that $L=\Bbb L_{\l,z}^{\dz}$ 
where $z\cdo\l\in I^s$, $z\cdo\l\preceq\boc$. Then (a) follows from 6.1(e),(g). In the setup of 
(b) we have
$$(\c(\Bbb L_{\l,s}^{\dz}))^j=(R_\l^{\dz})^{j+|z|+\nu+\r}((|z|+\nu+\r)/2)$$
and this is in $\cm^{\prec}G$ since $j+|z|+\nu+\r>a+\D+|z|$, see 6.1(f).

\subhead 6.4\endsubhead
Let $\cc^\spa\tG_s$ be the subcategory of $\cm(\tG_s)$ consisting of semisimple objects. Let 
$\cc_0^\spa\tG_s$ be the subcategory of $\cm_m(\tG_s)$ consisting of objects of 
pure of weight zero. Let $\cc^\boc\tG_s$ be the subcategory of $\cm(\tG_s)$ consisting of objects 
which are direct sums of objects in $CS_{\boc,s}$. Let $\cc_0^\boc\tG_s$ be the subcategory of 
$\cc^\spa_0\tG_s$ consisting of those $K$ such that, as an object of $\cc^\spa\tG_s$, $K$ belongs 
to $\cc^\boc\tG_s$. For $K\in\cc_0^\spa\tG_s$ let $\un{K}$ be the largest subobject of $K$ such 
that as an object of $\cc^\spa\tG_s$, we have $\un{K}\in\cc^\boc\tG_s$.

\subhead 6.5\endsubhead
For $L\in\cc^\boc_0Z_s$ we set 
$$\un\c(L)=\un{(\c(L))^{a+\nu}}((a+\nu)/2)
=\un{(\c(L))^{\{a+\nu\}}}\in\cc^\boc_0\tG_s.$$  
(The last equality uses that $\p$ in 6.1 is proper hence it preserves purity.) The functor 
$\un\c:\cc^\boc_0Z_s@>>>\cc^\boc_0\tG_s$ is called {\it truncated induction}. For 
$z\cdo\l\in\boc^s$ we have
$$\un\c(\Bbb L_{\l,s}^{\dz})=\un{(R_{\l,s}^{\dz})^{n_z}}(n_z/2).\tag a$$
Indeed,
$$\align&\un\c(\Bbb L_{\l,s}^{\dz})=\un{(\c(\Bbb L_{\l,s}^{\dz}))^{a+\nu}}((a+\nu)/2)
=\un{(\c(\cl^{\dz\sha}_{\l,s}\la|z|+\nu+\r\ra))^{a+\nu}}((a+\nu)/2)\\&
=\un{(\c(\cl^{\dz\sha}_{\l,s}))^{|z|+a+\D}}((|z|+a+\D)/2)
=\un{(\c(\cl^{\dz\sha}_{\l,s}))^{n_z}}(n_z/2)=\un{(R_{\l,s}^{\dz})^{n_z}}(n_z/2).\endalign$$
Using (a) and 6.2(d) we see that:

(d) {\it If $z\cdo\l\in\boc^s$ is such that $\un\c(\Bbb L_{\l,s}^{\dz})\ne0$ then
$z\cdo\l\underset\text{left}\to\si\ee^s(z\i)\cdo\l$.}

\subhead 6.6\endsubhead
For $z\cdo\l,z'\cdo\l'$ in $\boc^s$ we show:
$$\dim\Hom_{\cc^\boc\tG_s}(\un\c(\Bbb L_{\l,s}^{\dz}),\un\c(\Bbb L_{\l',s}^{\dz'}))
=\sum_{u\cdo\l_1\in\boc}
\tt(t_{u\i\cdo u(\l_1)}t_{z\cdo\l}t_{\ee^s(u)\cdo\ee^s(\l_1)}t_{z'{}\i\cdo z'(\l')})\tag a$$
where $\tt:\HH^\iy@>>>\ZZ$ is as in 1.9.
\nl
Let $()^\spa:\bbq@>>>\bbq$ be a field automorphism which maps any root of $1$ in $\bbq$ to its 
inverse. The field automorphism $\bbq(v)@>>>\bbq(v)$ which maps $v$ to $v$ and $x\in\bbq$ to 
$x^\spa$ is denoted again by ${}^\spa$.

Let $N_1$ (resp. $N_2$) be the left (resp. right) hand side of (a). 
Using 6.5(a) and the definitions we see that 
$$N_1=\sum_{A\in CS_{\boc,s}}(A:(R_{\l,s}^{\dz})^{n_z})(A:(R_{\l',s}^{\dz'})^{n_{z'}}).\tag b$$
Using 6.2(a) and the analogous identity for $(A:(R_{\l',s}^{\dz'})^{n_{z'}})$ in which the field 
automorphism $()^\spa:\bbq@>>>\bbq$ is applied to both sides (the left hand side is fixed by 
$()^\spa$), we deduce that 
$$\align&N_1=\\&(-1)^{|z|+|z'|}\sum_{E,E'\in\Irr_s(\HH^1_\fo)}
\sum_{A\in CS_{\boc,s}}b_{A,E}b_{A,E'}^\spa
\tr(\ee_st_{z\cdo\l},E^\iy)\tr(\ee_st_{z'\cdo\l'},E'{}^\iy)^\spa.\endalign$$
In the last sum we replace $\sum_{A\in CS_{\boc,s}}b_{A,E}b_{A,E'}^\spa$ by $1$ if 
$E'=E$ and by $0$ if $E'\ne E$. 
(In case A with $s\ne0$ we use \cite{\ORA, 3.9(i)} which assumes that the centre of $G$ is 
connected, but a similar proof applies without assumption on the centre. In case A with $s=0$
and in case B we use \cite{\CDGVII, 35.18(g)}.) 

We see that 
$$N_1=(-1)^{|z|+|z'|}\sum_{E\in\Irr_s(\HH^1_\fo)}
\tr(\ee_st_{z\cdo\l},E^\iy)\tr(\ee_st_{z'\cdo\l'},E^\iy)^\spa.$$
We now use the equality (for $E\in\Irr_s(\HH^1_\fo)$):
$$\tr(\ee_st_{z'\cdo\l'},E^\iy)^\spa=\tr(t_{z'{}\i\cdo z'(\l')}\ee_s\i,E^\iy)$$
which can be deduced from \cite{\CDGVII, 34.17}. We see that 
$$N_1=(-1)^{|z|+|z'|}\sum_{E\in\Irr_s(\HH^1_\fo)}
\tr(\ee_st_{z\cdo\l},E^\iy)\tr(t_{z'{}\i\cdo z'(\l')}\ee_s\i,E^\iy).$$
This is equal to $(-1)^{|z|+|z'|}$ times the trace of the linear map 
$\x\m t_{z\cdo\l}\ee^s(\x)t_{z'{}\i\cdo z'(\l')}$ from $\JJ_\fo$ to $\JJ_\fo$; hence it is
equal to 
$$(-1)^{|z|+|z'|}\sum_{u\cdo\l_1\in\fo}\tt(t_{u\i\cdo u(\l_1)}t_{z\cdo\l}
t_{\ee^s(u)\cdo\ee^s(\l_1)}t_{z'{}\i\cdo z'(\l')})=(-1)^{|z|+|z'|}N_2.$$
(In the last sum, the terms with $u\cdo\l_1\in\fo-\boc$ contribute $0$.) Thus, 
$N_1=(-1)^{|z|+|z'|}N_2$. Since $N_1$ and $N_2$ are natural numbers it follows that $N_1=N_2$. 
This proves (a).

\mpb

The proof above shows also that 
$\dim\Hom_{\cc^\boc\tG_s}(\un\c(\Bbb L_{\l,s}^{\dz}),\un\c(\Bbb L_{\l',s}^{\dz'}))=0$ 
whenever $(-1)^{|z|+|z'|}=-1$.

\mpb

Replacing in (a) $u\cdot\l_1$ by $\ee^{-s}(y)\cdo\ee^{-s}\l_1$ (recall that $\ee^s:\boc@>>>\boc$
is a bijection) we can rewrite (a) as follows:
$$\dim\Hom_{\cc^\boc\tG_s}(\un\c(\Bbb L_{\l,s}^{\dz}),\un\c(\Bbb L_{\l',s}^{\dz'}))
=\sum_{y\cdo\l_1\in\boc}
\tt(t_{\e^{-s}(y\i)\cdo \ee^{-s}(y(\l_1))}t_{z\cdo\l}t_{y\cdo\l_1}t_{z'{}\i\cdo z'(\l')}).$$
Since $N_1$ (in the form (b)) is symmetric in $z\cdo\l,z'\cdo\l'$, we have also
$$\dim\Hom_{\cc^\boc\tG_s}(\un\c(\Bbb L_{\l,s}^{\dz}),\un\c(\Bbb L_{\l',s}^{\dz'}))
=\sum_{y\cdo\l_1\in\boc}
\tt(t_{\e^{-s}(y\i)\cdo \ee^{-s}(y(\l_1))}t_{z'\cdo\l'}t_{y\cdo\l_1}t_{z\i\cdo z(\l)}).$$
Replacing $y\cdot\l_1$ by $y\i\cdot y(\l_1)$ (recall that $y\cdot\l_1\m y\i\cdot y(\l_1)$ is an 
involution $\boc@>>>\boc$) we can rewrite this as follows:
$$\dim\Hom_{\cc^\boc\tG_s}(\un\c(\Bbb L_{\l,s}^{\dz}),\un\c(\Bbb L_{\l',s}^{\dz'}))
=\sum_{y\cdo\l_1\in\boc}
\tt(t_{\e^{-s}(y)\cdo \ee^{-s}(\l_1)}t_{z'\cdo\l'}t_{y\i\cdo y(\l_1)}t_{z\i\cdo z(\l)}).\tag c$$

\mpb

We show:

(d) {\it There exist $z\cdo\l\in\boc^s$ such that $\un\c(\Bbb L_{\l,s}^{\dz})\ne0$.}
\nl
Let $k=u\cdo\l_1\in\boc$. Then $\ee^s(k)\in\boc$, $k^!\in\boc$ hence by 1.15(d) we have
$t_{k^!}t_jt_{\ee^s(k)}\ne0$ for some $j\in I$. From  2.5(a) we deduce that
$j\in\boc^s$. We can find $j'=z'\cdo\l'\in\boc$ such that
$t_{j'}$ appears with nonzero coefficient in $t_{k^!}t_jt_{\ee^s(k)}$. It follows that
$\tt(t_{k^!}t_jt_{\ee^s(k)}t_{j'{}^!})\ne0$. Since $\tt(\x\x')=\tt(\x'\x)$ for $\x,\x'\in\HH^\iy$
we deduce that $\tt(t_{\ee^s(k)}t_{j'{}^!}t_{k^!}t_j)\ne0$. In particular we have
$t_{\ee^s(k)}t_{j'{}^!}t_{k^!}\ne0$. Applying the antiautomorphism $t_u\m t_{u^!}$ of $\HH^\iy$
we deduce $t_kt_{j'}t_{\ee^s(k^!)}\ne0$. Using again 2.5(a) we deduce that $j'\in\boc^s$.
If $i\in\boc$, $j\in I$ satisfy $t_{i^!}t_jt_{\ee^s(i)}\ne0$ then $j\in\boc^s$. Since 
$\tt(t_{h^!}t_jt_{\ee^s(h)}t_{j'{}^!})\in\NN$ for any $h\in\boc$ and
$\tt(t_{k^!}t_jt_{\ee^s(k)}t_{j'{}^!})\ne0$, we see that
$\sum_{h\in\boc}\tt(t_{h^!}t_jt_{\ee^s(h)}t_{j'{}^!})\in\NN_{>0}$.
Using this and (a), we see that
$$\dim\Hom_{\cc^\boc\tG_s}(\un\c(\Bbb L_{\l,s}^{\dz}),\un\c(\Bbb L_{\l',s}^{\dz'}))\in\NN_{>0}.$$
This proves (d).

\mpb

The following converse to 6.2(e) is an immediate consequence of (d):

(e) {\it We have $CS_{\boc,s}\ne\emp$.}

\subhead 6.7\endsubhead
Let $L\in\cc^\boc_0Z_s$. We show that $\fD(L)\in\cc^{\wt{\boc}}_0Z_s$. (Here $\wt{\boc}$ is as in
1.14.) It is enough to note that for $w\cdo\l\in\boc^s$ and $\o\in\k\i_0(w)$ we have

(a) $\fD(\Bbb L_{\l,s}^\o)=\Bbb L_{\l\i,s}^\o$. 
\nl
We show: 

(b) {\it For $L\in\cc^\boc_0Z_s$ we have canonically $\un\c(\fD(L))=\fD(\un\c(L))$ where the 
first $\un\c$ is relative to $\wt{\boc}$ instead of $\boc$.}
\nl
Let $\p,f,\dZ_s$ be as in 6.1. By the relative hard Lefschetz theorem \cite{\BBD, 5.4.10} applied 
to the projective morphism $\p$ and to $f^*L\la\nu\ra$ (a perverse sheaf of pure weight $0$ on 
$\dZ_s$) we have canonically for any $j\in\ZZ$:
$$(\p_!f^*L\la\nu\ra)^{-j}=(\p_!f^*L\la\nu\ra)^j(j).\tag c$$
We have used the fact that $f$ is smooth with fibres of dimension $\nu$. This also shows that
$$\fD(\c(\fD(L)))=\c(L)\la2\nu\ra.\tag d$$
Using (d) we have
$$\align&\fD(\un\c(\fD(L)))=\fD((\c(\fD(L)))^{a+\nu}((a+\nu)/2)))
=(\fD(\c(\fD(L))))^{-a-\nu}((-a-\nu)/2)\\&
=(\c(L)\la2\nu\ra)^{-a-\nu}((-a-\nu)/2)=(\c(L)\la\nu\ra)^{-a}(-a/2).\endalign$$
Hence using (c) we have
$$\fD(\un\c(\fD(L)))=(\c(L)\la\nu\ra)^a(a/2)
=(\c(L))^{a+\nu}((a+\nu)/2)=\un\c(L).$$
This proves (b).

\subhead 6.8\endsubhead
We define $\z:\cd(\tG_s)@>>>\cd(Z_s)$ and $\z:\cd_m(\tG_s)@>>>\cd_m(Z_s)$ by
$\z(K)=f_!\p^*K$ where $Z_s@<f<<\dZ_s@>\p>>\tG_s$ is as in 6.1(a). We show:

(a) {\it For any $L\in\cd(Z_s)$ or $L\in\cd_m(Z_s)$ we have $\fb''(L)=\z(\c(L))$.}
\nl
We have $\z(\c(L))=f_!\p^*\p_!f^*(L)$. We have
$$\dZ_s\T_{\tG_s}\dZ_s=\{((B_0,B_1,B_2,B_3),\g)\in\cb^4\T\tG_s;
\g B_0\g\i=B_3,\tg B_1\tg\i=B_2\}.$$
We have a cartesian diagram
$$\CD
\dZ_s\T_{\tG_s}\dZ_s@>\ti\p_1>>\dZ_s\\
@V\ti\p_2VV          @V\p VV\\
\dZ_s    @>\p>>   \tG_s
\endCD$$
where $\ti\p_1((B_0,B_1,B_2,B_3),\g)=(B_0,B_3,\g)$,
$\ti\p_2((B_0,B_1,B_2,B_3),\g)=(B_1,B_2,\g)$.
It follows that $\p^*\p_!=\ti\p_{1!}\ti\p_2^*$. Thus,
$$\z(\c(L))=f_!\ti\p_{1!}\ti\p_2^*f^*(L)=(f\ti\p_1)_!(f\ti\p_2)^*(L).$$
Define $\p'_1:\dZ_s\T_{\tG_s}\dZ_s@>>>Z_s$, $\p'_2:\dZ_s\T_{\tG_s}\dZ_s@>>>Z_s$
 by 
$$\p'_1((B_0,B_1,B_2,B_3),\g)=(B_0,B_3,\g U_{B_0}),$$
$$\p'_2((B_0,B_1,B_2,B_3),\g)=(B_1,B_2,\g U_{B_1}).$$
Then $\p'_1=f\ti\p_1$, $\p'_2=f\ti\p_2$ and $\z(\c(L))=\p'_{1!}\p'_2{}^*(L)$. Let ${}^\di\cy$ be 
as in 4.14. We have an isomorphism ${}^\di\cy@>>>\dZ_s\T_{\tG_s}\dZ_s$ induced by
$$((x_0\UU,x_1\UU,x_2\UU,x_3\UU),\g)\m
((x_0\BB x_0\i,x_1\BB x_1\i,x_2\BB x_2\i,x_3\BB x_3\i),\g).$$
We use this to identify ${}^\di\cy=\dZ_s\T_{\tG_s}\dZ_s$. Then $\p'_1,\p'_2$ become
$d,{}^\di\et$ of 4.25. We see that (a) holds.

\subhead 6.9\endsubhead
{\it In the remainder of this section we assume that $s\in\ZZ_\boc$.}
\nl
Let $z\cdo\l\in\fo$. We set $\Si=\e_s^*\z(R_{\l,s}^{\dz})\la2\nu+|z|\ra\in\cd(\tcb^2)$. Let 
$j\in\ZZ$. We show:

(a) {\it If $z\cdo\l\preceq\boc$, then $\Si^j\in\cm^{\preceq}\tcb^2$.}

(b) {\it If $z\cdo\l\prec\boc$, then $\Si^j\in\cm^{\prec}\tcb^2$.}

(c) {\it If $z\cdo\l\in\boc$ and $j>\nu+2\r+2a$, then $\Si^j\in\cm^{\prec}\tcb^2$.}
\nl
If $z\cdo\l\n I^s$, then $\Si=0$ and there is nothing to prove. Now assume that $z\cdo\l\in I^s$. 
Using 4.9(a), we have
$$\Si=\e_s^*\z(\c(\cl_{\l,s}^{\dz\sha}))\la2\nu+|z|\ra=
\fb'(\cl_{\l,s}^{\dz\sha})\la2\nu+|z|\ra=\fb'(\Bbb L_{\l,s}^{\dz})\la\nu-\r\ra.$$
Now (a),(b) follow from 4.14(a),(b); (c) follows from 4.14(c). (If $j>\nu+2\r+2a$, then 
$j+\nu-r>2\nu+\r+2a$.)

\subhead 6.10\endsubhead 
We show:

(a) {\it If $K\in\cd^{\preceq}\tG_s$, then $\z(K)\in\cd^{\preceq}Z_s$.}

(b) {\it If $K\in\cd^{\prec}\tG_s$, then $\z(K)\in\cd^{\prec}Z_s$.}

(c) {\it If $K\in\cd^{\preceq}\tG_s$ and $j>\nu+a$, then $(\z(K))^j\in\cm^{\prec}Z_s$.}
\nl
We can assume in addition that $K=A\in CS_{\boc',s}$ for a two-sided cell $\boc'$ such that 
$\boc'\preceq\boc$. Assume first that $\boc'=\boc$. By 6.2(c) we can find $z\cdo\l\in\boc$ such 
that $(A:(R_{\l,s}^{\dz})^{n_z})\ne0$. Then $A[-n_z]$ (without mixed structure) is a direct 
summand of the semisimple complex $R_{\l,s}^{\dz}$. Hence $\e^*_s\z(A)[-n_z]$ is a direct summand 
of $\e^*_s\z(R_{\l,s}^{\dz})$ and $\e^*_s\z(A)[-n_z+2\nu+|z|]$ is a direct summand of $\Si$ (in 
6.9), that is, $\e^*_s\z(A)[-a-\r]$ is a direct summand of $\Si$. By 6.9, if $j\in\ZZ$ (resp. 
$j>\nu+2\r+2a$) then $\Si^j\in\cm^{\preceq}\tcb^2$ (resp. $\Si^j\in\cm^{\prec}\tcb^2$) hence 
$(\e^*_s\z(A)[-a-\r])^j\in\cm^{\preceq}\tcb^2$ (resp. 
$(\e_s^*\z(A)[-a-\r])^j\in\cm^{\prec}\tcb^2$), that is,
$(\e^*_s\z(A))^{j-a-\r}\in\cm^{\preceq}\tcb^2$ (resp. 
$(\e^*_s\z(A))^{j-a-\r}\in\cm^{\prec}\tcb^2$). We see that if $j'\in\ZZ$ (resp. 
$j'>\nu+\r+a$) then $(\e_s^*\z(A))^{j'}\in\cm^{\preceq}\tcb^2$ (resp. 
$(\e^*_s\z(A))^{j'}\in\cm^{\prec}\tcb^2$), so that $(\z(A))^{j'-\r}\in\cm^{\preceq}Z_s$ (resp. 
$(\z(A))^{j'-\r}\in\cm^{\prec}Z_s$); here we use 4.3(a). We see that if $j\in\ZZ$ (resp. 
$j>\nu+a$, so that $j+\r>\nu+\r+a$), then $(\z(A))^j\in\cm^{\preceq}Z_s$ (resp.
$(\z(A))^j\in\cm^{\prec}Z_s$). Thus the desired results hold when $\boc'=\boc$.

Assume now that $\boc'\prec\boc$. Applying the above argument with $\boc$ replaced by $\boc'$, we 
see that the desired results hold.

\subhead 6.11\endsubhead 
For $K\in\cc^\boc_0\tG_s$ we set
$$\un\z(K)=\un{(\z(K))^{\{\nu+a\}}}\in\cc^\boc_0Z_s.$$
We say that $\un\z(K)$ is the {\it truncated restriction} of $K$.

\subhead 6.12\endsubhead 
Let $L\in\cc^\boc_0Z_s$. We show:

(a) {\it We have canonically $\un\z(\un\c(L))=\un{\fb''}(L)$.}
\nl
We shall apply the method of \cite{\CONV, 1.12} with $\Ph:\cd_m(Y_1)@>>>\cd_m(Y_2)$ replaced by 
$\z:\cd_m(\tG_s)@>>>\cd_m(Z_s)$ and with $\cd^{\preceq}(Y_1)$, $\cd^{\preceq}(Y_2)$ replaced by 
$\cd^{\preceq}\tG_s$, $\cd^{\preceq}Z_s$. We shall take $\XX$ in {\it loc.cit.} equal to $\c(L)$. 
The conditions of {\it loc.cit.} are satisfied: those concerning $\XX$ are satisfied with 
$c'=a+\nu$, see 6.3. The conditions concerning $\z$ are satisfied with $c=a+\nu$, see 6.10. We see
that 
$$(\z(\c(L)))^j=0\text{ if }j>2a+2\nu\tag b$$
and
$$\un{gr_{2a+2\nu}((\z(\c(L)))^{2a+2\nu})}(a+\nu)=\un\z(\un\c(L)).\tag c$$
Since $\z(\c(L))=\fb''(L)$, we see that the left hand side of (c) equals $\un{\fb''}(L)$. Thus (a)
is proved.

\mpb

Combining (a) with 4.25(d) and 4.14(d) we see that

(b) {\it we have canonically $\ti\e_s\un\z(\un\c(L))=\un{\fb}(L)$.}

\subhead 6.13\endsubhead
Let $K\in\cd(\tG_s)$ and let $L\in\cd^\spa\tcb^2$. Let $\tL=(\ee^s)^*L$. In (a) below the
assumption $s\in\ZZ_\boc$ is not used:

(a) {\it there is a canonical isomorphism $\tL\cir\e_s^*\z(K)@>\si>>\e_s^*\z(K)\cir L$.}
\nl
Let $Y=\tcb^2\T\tG_s$. Define $j:Y@>>>\tG_s$ by $j(x_0\UU,x_1\UU,\g)=\g$. Define $j_1:Y@>>>\tcb^2$
by $j_1(x_0\UU,x_1\UU,\g)=(x_0\UU,\g\i x_1\t^s\UU)$. Define $j_2:Y@>>>\tcb^2$ by 
$j_2(x_0\UU,x_1\UU,\g)=(\g x_0\t^{-s}\UU,x_1\UU)$. From the definitions we have 
$\tL\cir\e^*_s\z(K)=j_{2!}(j_1^*(\tL)\ot j^*(K))$, $\e^*_s\z(K)\cir L=j_{2!}(j_2^*(L)\ot j^*(K))$. 
It remains to prove that $j_1^*(\tL)=j_2^*L$ that is, $j'_1{}^*L=j_2^*L$ where
$j'_1=\ee^sj_1:Y@>>>\tcb^2$ is given by 
$j'_1(x_0\UU,x_1\UU,\g)=(\t^sx_0\t^{-s}\UU,\t^s\g\i x_1\UU)$. The equality $j'_1{}^*L=j_2^*L$
follows from the $G$-equivariance of $L$. This proves (a).

\mpb

Now let $K\in\cc^\boc_0\tG_s$ and let $L\in\cc^\boc_0\tcb^2$. Since $\ee^s(\boc)=\boc$, 
we have $(\ee^s)^*L\in\cc^\boc_0\tcb^2$, see 3.11(a). We show that

(b) {\it there is a canonical isomorphism 
$(\ee^s)^*(L)\un\cir\ti\e_s\un\z(K)@>\si>>(\ti\e_s\un\z(K))\un\cir L$.}
\nl
We apply the method of \cite{\CONV, 1.12} with 
$\Ph:\cd^{\preceq}_m\tcb^2@>>>\cd^{\preceq}_m\tcb^2$, $L'\m L'\cir L$, $\XX=\ti\e_s\z(K)$ and with
$(c,c')=(a-\nu,\nu+a)$, see \cite{\MONO, 2.23(a)}, 6.10(c). We deduce that we have canonically
$$\un{\un{((\ti\e_s\z(K))^{\{a+\nu\}}}\cir L)^{\{a-\nu\}}}=\un{(\ti\e_s\z(K)\cir L)^{\{2a\}}}.
\tag c$$
We apply the method of \cite{\CONV, 1.12} with 
$\Ph:\cd^{\preceq}_m \tcb^2@>>>\cd^{\preceq}_m\tcb^2$, $L'\m(\ee^s)^*L\cir L'$, 
$\XX=\ti\e_s\z(K)$ and with $(c,c')=(a-\nu,\nu+a)$, see \cite{\MONO, 2.23(a)}, 6.10(c). We deduce 
that we have canonically
$$\un{(((\ee^s)^*L\cir\un{(\ti\e\z(K))^{\{a+\nu\}})^{\{a-\nu\}}}}
=\un{((\ee^s)^*L\cir\ti\e\z(K))^{\{2a\}}}.\tag d$$
We now combine (c),(d) with (a); we obtain (b).

\subhead 6.14\endsubhead
Let $s',s''$ be integers. Let $\mu:\tG_{s'}\T\tG_{s''}@>>>\tG_{s'+s''}$ be the multiplication map. For 
$K\in\cd(\tG_{s'}),K'\in\cd(\tG_{s''})$ (resp. $K\in\cd_m(\tG_{s'}),K'\in\cd_m(\tG_{s''})$) we set 
$K*K'=\mu_!(K\bxt K')$; this is in $\cd(\tG_{s'+s''})$ (resp. in $\cd_m(\tG_{s'+s''})$). For 
$K\in\cd(\tG_{s_1}),K'\in\cd(\tG_{s_2})$, $K''\in\cd(\tG_{s_3})$ we have canonically
$(K*K')*K''=K*(K'*K'')$ (and we denote this by $K*K'*K''$). For 
$K\in\cm(\tG_{s'}),K'\in\cm(\tG_{s''})$ we show:

(a) {\it If $K'$ is $G$-equivariant then we have canonically $K*K'=((\ee^{-s'})^*K')*K'$. If $K$ 
is $G$-equivariant then we have canonically $K*K'=K'*((\ee^{s''})^*K)$.}
\nl
The proof is immediate. It will be omitted. (Compare \cite{\CONV, 4.1}.)

\subhead 6.15\endsubhead
Let $s',s''\in\ZZ$. We show:

(a) {\it For $K\in\cd(\tG_{s'})$, $L\in\cd(Z_{s''})$ we have canonically 
$K*\c(L)=\c(L\bul\z(K))$.}
\nl
Let $Y=\tG_{s'}\T\tG_{s''}\T\cb$. Define $c:Y@>>>\tG_{s'}\T Z_{s''}$ by 
$$c(\g_1,\g_2,B)=(\g_1,(B,\g_2B\g_2\i,\g_2U_B));$$ define $d:Y@>>>\tG_{s'+s''}$ by
$d(\g_1,\g_2,B)=\g_1\g_2$. From the definitions we see that both $K*\c(L)$, $\c(L\bul\z(K))$ can 
be identified with $d_!c^*(K\bxt L)$. This proves (a).

Now let $L\in\cd(Z_{s'}),L'\in\cd(Z_{s''})$. Replacing in (a) $K,L$ by $\c(L),L'$ and using 6.8(a),
we obtain
$$\c(L)*\c(L')=\c(L'\bul\fb''(L)).\tag b$$

\subhead 6.16\endsubhead
Let $s'\in\ZZ_\boc$. Let $L\in\cd^\spa(Z_s)$, $L'\in\cd^\spa(Z_{s'})$, $j\in\ZZ$. We show:

(a) {\it If $L\in\cd^{\preceq}Z_s$ or $L'\in\cd^{\preceq}Z_{s'}$ then 
$L'\bul\fb''(L)\in\cd^{\preceq}Z_{s+s'}$.}

(b) {\it If $L\in\cd^{\prec}Z_s$ or $L'\in\cd^{\prec}Z_{s'}$ then 
$L'\bul\fb''(L)\in\cd^{\prec}Z_{s+s'}$.}

(c) {\it If $L\in\cm^{\preceq}Z_s$, $L'\in\cm^\spa Z_{s'}$ and $j>3a+\r+\nu$ then
$(L'\bul\fb''(L))^j\in\cd^{\prec}Z_{s+s'}$.}
\nl
Now (a),(b) follow from 4.25(b) and 4.23(a). To prove (c) we may assume that 
$L=\Bbb L_{\l,s}^{\dw}$, $L'=\Bbb L_{\l',s'}^{\dw'}$ with $w\cdo\l\in I^s_n$,
$w'\cdo\l'\in I^{s'}_n$ and $w\cdo\l\preceq\boc$. We apply the method of \cite{\CONV, 1.12} with
$\Ph:\cd^{\preceq}Z_s@>>>\cd^{\preceq}Z_{s+s'}$, $L_1\m L'\bul L_1$ and
$\XX=\fb''(L)$ and with $c'=2\nu+2a$ (see 4.25(c)), $c=a+\r-\nu$ (see 4.23(b)). We have 
$c+c'=\nu+\r+3a$ hence (c) holds.

\subhead 6.17\endsubhead
Let $s'\in\ZZ_\boc$. Let $L\in\cd^\spa(Z_s)$, $L'\in\cd^\spa(Z_{s'})$, $j\in\ZZ$. We show:

(a) {\it If $L\in\cd^{\preceq}Z_s$ or $L'\in\cd^{\preceq}Z_{s'}$ then 
$\c(L'\bul\fb''(L))\in\cd^{\preceq}\tG_{s+s'}$.}

(b) {\it If $L\in\cd^{\prec}Z_s$ or $L'\in\cd^{\prec}Z_{s'}$ then 
$\c(L'\bul\fb''(L))\in\cd^{\prec}\tG_{s+s'}$.}

(c) {\it If $L\in\cm^{\preceq}Z_s$, $L'\in\cm^\spa Z_{s'}$ and $j>4a+2\nu+\r$ then
$(\c(L'\bul\fb''(L)))^j\in\cm^{\prec}\tG_{s+s'}$.}
\nl
(a),(b) follow from 6.3(a) using 6.16(a),(b). To prove (c) we can assume that 
$L=\Bbb L_{\l,s}^{\dw}$, $L'=\Bbb L_{\l',s'}^{\dw'}$ with $w\cdo\l\in I^s_n$,
$w'\cdo\l'\in I^{s'}_n$ and $w\cdo\l\preceq\boc$. We apply the method of \cite{\CONV, 1.12} with 
$\Ph:\cd^{\preceq}Z_{s+s'}@>>>\cd^{\preceq}\tG_{s+s'}$, $L_1\m\c(L_1)$,
$\XX=L'\bul\fb''(L)$ and with $c'=\nu+\r+3a$ (see 6.16(c)), $c=a+\nu$ (see 6.3(b)). We have 
$c+c'=2\nu+\r+4a$ hence (c) holds. 

\subhead 6.18\endsubhead
Let $s'\in\ZZ_\boc$. Let $K\in\cd^\spa(\tG_s)$, $K'\in\cd^\spa(\tG_{s'})$. We show:

(a) {\it If $K\in\cd^{\preceq}\tG_s$ or $K'\in\cd^{\preceq}\tG_{s'}$ then 
$K*K'\in\cd^{\preceq}G_{s+s'}$.}

(b) {\it If $K\in\cd^{\prec}\tG_s$ or $K'\in\cd^{\prec}\tG_{s'}$ then 
$K*K'\in\cd^{\prec}\tG_{s+s'}$.}

(c) {\it If $K\in\cd^{\preceq}\tG_s$ or $K'\in\cd^{\preceq}\tG_{s'}$ and $j>2a+\r$ then 
$(K*K')^j\in\cd^{\prec}\tG_{s+s'}$.}
\nl
We can assume that $K=A\in CS_{\fo,s}$, $K'=A'\in CS_{\fo,s'}$. Let $A''\in\cm(\tG_{s+s'})$ be a 
composition factor of $(A*A')^j$. By 6.2(c) we can find $w\cdo\l\in\boc_A$, 
$w'\cdo\l'\in\boc_{A'}$ such that $(A:(R_{\l,s}^{\dw})^{n_w})\ne0$, 
$(A':(R_{\l',s'}^{\dw'})^{n_{w'}})\ne0$. Then $A$ is a direct summand of $R_{\l,s}^{\dw}[n_w]$ and
$A'$ is a direct summand of $R_{\l',s'}^{\dw'}[n_{w'}]$. Hence $A*A'$ is a direct summand of 
$$R_{\l,s}^{\dw}*R_{\l',s'}^{\dw'}[a(w\cdo\l)+a(w'\cdo\l')+|w|+|w'|+2\D]$$ 
and $(A*A')^j$ is a direct summand of 
$$\align&(R_{\l,s}^{\dw}*R_{\l',s'}^{\dw'}[|w|+|w'|+2\nu+2\r])^{j+a(w\cdo\l)+a(w'\cdo\l')
+2\nu}\\&=(\c(\Bbb L_{\l,s}^{\dw})*\c(\Bbb L_{\l',s'}^{\dw'}))^{j+a(w\cdo\l)+
a(w'\cdo\l')+2\nu}.\endalign$$
Using 6.15(b) we see that $(A*A')^j$ is a direct summand of 
$$(\c(\Bbb L_{\l',s'}^{\dw'}\bul\fb''(\Bbb L_{\l,s}^{\dw}))^{j+a(w\cdo\l)+a(w'\cdo\l')+2\nu}.
\tag d$$
Hence $A''$ is a composition factor of (d). Using 6.17(a) we see that 
$A''\in CS_{\fo,s+s'}$, that $\boc_{A''}\preceq w\cdo\l$ and that $\boc_{A''}\preceq w'\cdo\l'$.
In the setup of (a) we have $w\cdo\l\preceq\boc$ or $w'\cdo\l'\preceq\boc$ hence 
$\boc_{A''}\le\boc$. Thus (a) holds. Similarly, (b) holds. In the setup of (c) we have 
$w\cdo\l\preceq\boc$ and $w'\cdo\l'\preceq\boc$. Hence $a(w\cdo\l)\ge a$, $a(w'\cdo\l')\ge a$. 
(See Q3 in 1.9.) Assume that $\boc_{A''}=\boc$. Since $A''$ is a composition factor of (d), we see
from 6.17(c) that 
$$j+a(w\cdo\l)+a(w'\cdo\l')+2\nu\le 4a+2\nu+\r$$
hence $j+2a+2\nu\le4a+2\nu+\r$ and $j\le2a+\r$. This proves (c).

\subhead 6.19\endsubhead
Let $s'\in\ZZ_\boc$. 
For $K\in\cc^\boc_0\tG_s$, $K'\in\cc^\boc_0\tG_{s'}$, we set
$$K\un{*}K'=\un{(K*K')^{\{2a+\r\}}}\in\cc^\boc_0\tG_{s+s'}.$$
We say that $K\un{*}K'$ is the {\it truncated convolution} of $K,K'$. Note 
that 6.14(a) induces for $K,K'\in\cc^\boc_0G$ a canonical isomorphism
$$K\un{*}K'=K'\un{*}((\ee^{s'})^*K).\tag a$$
Let $L\in\cc^\boc_0Z_{s'}$, $K\in\cc^\boc_0\tG_s$. Using the method of \cite{\CONV, 1.12} several 
times, we see that
$$K\un{*}\un\c(L)=\un{gr_k((K*\c(L))^k)}(k/2)$$
where $k=(a+\nu)+(2a+\r)=3a+\nu+\r$ and
$$\un\c(L\un{\bul}\un\z(K))=\un{gr_{k'}((\c(L\bul\z(K))^{k'})}(k'/2)$$
where $k'=(a+\nu)+(a+\nu)+(a+\r-\nu)=3a+\nu+\r$.
Using now 6.15(a) and the equality $k=k'$ we obtain
$$K\un{*}\un\c(L)=\un\c(L\un{\bul}\un\z(K)).\tag b$$

\mpb

Let $L\in\cc^\boc_0Z_s$, $L'\in\cc^\boc_0Z_{s'}$.  Using the method of \cite{\CONV, 1.12} several 
times, we see that
$$\un\c(L)\un{*}\un\c(L')=\un{gr_k((\c(L)*\c(L'))^k)}(k/2)$$
where $k=(a+\nu)+(a+\nu)+(2a+\r)=4a+2\nu+\r$ and
$$\un\c(L'\un{\bul}\un{\fb''}(L)=\un{gr_{k'}((\c(L'\bul\fb''(L)))^{k'})}(k'/2)$$
where $k'=(2a+2\nu)+(a+\r-\nu)+(a+\nu)=4a+2\nu+\r$. Using now 6.15(b) and the equality $k=k'$ we 
obtain
$$\un\c(L)\un{*}\un\c(L')=\un\c(L'\un{\bul}(\un{\fb''}(L))).\tag c$$

\mpb

We show (assuming that $s_h\in\ZZ_\boc$ for $h=1,2,3$):

(d) {\it For $K\in\cc^\boc_0\tG_{s_1},K'\in\cc^\boc_0\tG_{s_2}, K''\in\cc^\boc_0\tG_{s_3}$, there is a 
canonical isomorphism $(K\un{*}K')\un{*}K''@>\si>>K\un{*}(K'\un{*}K'')$.}
\nl
Indeed, just as in \cite{\CONV, 4.7} we can identify, using the method of \cite{\CONV, 1.12}, both
$(K\un{*}K')\un{*}K''$ and $K\un{*}(K'\un{*}K'')$ with $\un{(K*K'*K'')^{\{4a+2\r\}}}$.

\subhead 6.20\endsubhead
Let $s',s''\in\ZZ$. For $K\in\cd(\tG_{s'})$, $K'\in\cd(\tG_{s''})$,  we show:

(a) {\it We have canonically $\z(K*K')=\z(K')\bul\z(K)$.}
\nl
Let 
$$Y=\{(B,\g U_B,\g_1,\g_2);B\in\cb,\g\in\tG_{s'+s''},\g_1\in\tG_{s'},\g_2\in\tG_{s''};
\g_1\g_2\in\g U_B\}.$$ 
Define $j_1:Y@>>>\tG_{s'}$, $j_2:Y@>>>\tG_{s''}$ by $j_1(B,\g U_B,\g_1,\g_2)=\g_1$,
$j_2(B,\g U_B,\g_1,\g_2)=\g_2$. Define $j:Y@>>>Z_{s'+s''}$ by 
$j(B,\g U_B,\g_1,\g_2)=(B,\g B\g\i,\g U_B)$. From the definitions we have 
$\z(K*K')=j_!(j_1^*(K)\ot j_2^*(K'))=\z(K')\bul\z(K)$; (a) follows.

\mpb

Let $s'\in\ZZ_\boc$. For $K\in\cd^\boc_0(G_s)$, $K'\in\cd^\boc_0(G_{s'})$, we show:

(b) {\it We have canonically $\un\z(K\un{*}K')=\un\z(K')\un{\bul}\un\z(K)$.}
\nl
Using the method of \cite{\CONV, 1.12} we see that
$$\un\z(K\un{*}K')=\un{gr_k((\z(K*K'))^k)}(k/2)$$
where $k=(a+\nu)+(2a+\r)=3a+\nu+\r$ and that
$$\un\z(K')\un{\bul}\un\z(K)=\un{gr_{k'}((\z(K)\bul\z(K'))^{k'})}(k'/2)$$
where $k'=(a+\r-\nu)+(a+\nu)+(a+\nu)=3a+\nu+\r$. It remains to use (a) and the
equality $k=k'$.

\subhead 6.21\endsubhead
Let $s'\in\ZZ$. Define $h:\tG_{s'}@>>>\tG_{-s'}$ by $\g\m\g\i$. For $K\in\cd(\tG_{-s'})$ we set 
$K^\da=h^*K\in\cd(\tG_{s'})$. We show:

(a) {\it For $L\in\cd(Z_{-s'})$ we have $(\c(L))^\da=\c(L^\da)$ with $L^\da$ as in 4.2.}
\nl
This follows from the definition of $\c$ using the commutative diagram
$$\CD
Z_{s'}@<f<<\dZ_{s'}@>\p>>\tG_{s'}\\
@V\fh VV @V\dot{\fh}VV @VhVV\\
Z_{-s'}@<f<<\dZ_{-s'}@>\p>>\tG_{-s'}
\endCD$$
where $f,\p$ are as in 6.1, $\fh$ is as in 4.2 and $\dot{\fh}:\dZ_{s'}@>>>\dZ_{-s'}$ is 
$(B,B',\g)\m(B',B,\g\i)$.

\mpb

From (a) and 4.3(e) we see that, if $w\cdo\l\in I^{-s}_n$, then
$$(\c(\Bbb L_{\l,-s}^{\dw}))^\da=\c(\Bbb L_{w(\l)\i,s}^{\dw\i}).\tag b$$
We deduce that

(c) {\it if $A\in CS_{\boc,-s}$, then $A^\da\in CS_{\wt{\boc},s}$.}
\nl
From (a),(c) we deduce:

(d) {\it For $L\in\cc^\boc_0Z_{-s}$ we have $(\un\c(L))^\da=\un\c(L^\da)$ where the second 
$\un\c$ is relative to $\ti{\boc},\fo\i$ instead of $\boc,\fo$.}

\head 7. Equivalence of $\cc^\boc\tG_s$ with the $\ee^s$-centre of $\cc^\boc\tcb^2$\endhead
\subhead 7.1\endsubhead
{\it In this section (except in 7.8) let $\boc,\fo,a,n,\Ps$ be as in 3.1(a).}
\nl
In this subsection we assume that $s\in\ZZ_\boc$.
Let $u:\tG_{-s}@>>>\pp$ be the obvious map; let $\ph:\pp@>>>G$ be the map with image $\{1\}$.
From \cite{\CSII, 7.4} we see that for $K,K'$ in $\cm_m\tG_{-s}$ we have canonically
$$(u_!(K\ot K'))^0=\Hom_{\cm(\tG_{-s})}(\fD(K),K'),\qua (u_!(K\ot K'))^j=0 \text{ if }j>0.$$
We deduce that if $K,K'$ are also pure of weight $0$ then $(u_!(K\ot K'))^0$  is pure of weight $0$
that is, $(u_!(K\ot K'))^0=gr_0(u_!(K\ot K'))^0$. From the definitions we see that we have
$u_!(K\ot K')=\ph^*(K^\da*K')$ where $K^\da\in\cm_m(\tG_s)$ is as in 6.21. Hence, for $K'$ in 
$\cc^\boc_0\tG_{-s}$ and $K$ in $\cc^{\ti{\boc}}_0\tG_{-s}$ (so that $K^\da\in\cc^\boc_0\tG_s$, 
see 6.21(c)) we have
$$\Hom_{\cm(\tG_{-s})}(\fD(K),K')=(\ph^*(K^\da*K'))^0=(\ph^*(K^\da*K'))^{\{0\}}.\tag a$$
Using \cite{\CONV, 8.2} with $\Ph:\cd^{\preceq}_m\tG_0@>>>\cd_m\pp$, $K_1\m\ph^*K_1$, $c=-2a-\r$ 
(see \cite{\MONO, 6.8(a)}), $K$ replaced by $K^\da*K'\in\cd_m(\tG_0)$ and $c'=2a+\r$, we see that 
we have canonically
$$(\ph^*(K^\da\un{*}K'))^{\{-2a-\r\}}\sub(\ph^*(K^\da*K'))^{\{0\}}.$$
In particular, if $L\in\cc^\boc_0Z_{-s}$, $L'\in\cc^\boc_0Z_s$,  then we have canonically
$$(\ph^*(\un\c(L')\un{*}\un\c(L)))^{\{-2a-\r\}}\sub(\ph^*(\un\c(L')*\un\c(L)))^{\{0\}}.$$
Using the equality
$$(\ph^*(\un\c(L')\un{*}\un\c(L)))^{\{-2a-\r\}}=\ph^*(\un\c(L\un\bul\un\z(\un\c(L')))))^{-2a-\r}$$
which comes from 6.19(b), we deduce that we have canonically
$$\ph^*(\un\c(L\un\bul\un\z(\un\c(L')))))^{-2a-\r}\sub(\ph^*(\un\c(L')*\un\c(L)))^{\{0\}},$$
or equivalently, using (a) with $K, K'$ replaced by $\un\c(L')^\da$, $\un\c(L)$,
$$\align&\ph^*(\un\c(L\un\bul\un\z(\un\c(L'))))^{-2a-\r}\sub
\Hom_{\cc^\boc\tG_{-s}}(\fD(\un\c(L')^\da),\un\c(L))\\&=
\Hom_{\cc^\boc\tG_s}(\fD(\un\c(L)^\da),\un\c(L')).\endalign$$
Using now \cite{\MONO, 6.9(d)} with $L$ replaced by $L\un\bul\un\z(\un\c(L'))\in\cc^\boc_0Z_0$,
we have canonically
$$\ph^*(\un\c(L\un\bul\un\z(\un\c(L'))))^{-2a-\r}=
\Hom_{\cc^\boc Z_0}(\bold1'_0,L\un\bul\un\z(\un\c(L'))).$$
Thus we have canonically
$$\Hom_{\cc^\boc Z_0}(\bold1'_0,L\un\bul\un\z(\un\c(L')))\sub
\Hom_{\cc^\boc\tG_s}(\fD(\un\c(L)^\da),\un\c(L'))$$
or equivalently (using 5.8(a))
$$\Hom_{\cc^\boc Z_{-s}}(\fD(\un\z(\un\c(L'))^\da),L)
\sub\Hom_{\cc^\boc\tG_s}(\fD(\un\c(L)^\da),\un\c(L')).$$
Now we have
$$\align&\Hom_{\cc^\boc Z_{-s}}(\fD(\un\z(\un\c(L'))^\da),L)
=\Hom_{\cc^{\ti{\boc}}Z_{-s}}(\fD(L),\un\z(\un\c(L'))^\da)\\&=
\Hom_{\cc^\boc Z_s}((\fD(L))^\da,\un\z(\un\c(L'))),\endalign$$
hence 
$$\Hom_{\cc^\boc Z_s}((\fD(L))^\da,\un\z(\un\c(L')))
\sub\Hom_{\cc^\boc\tG_s}(\fD(\un\c(L)^\da),\un\c(L')).$$
We set ${}^1L=\fD(L^\da)=(\fD(L))^\da\in\cc^\boc_0Z_s$ and note that 
$$\fD(\un\c(L)^\da)=\fD(\un\c(L^\da))=\un\c(\fD(L^\da))=\un\c({}^1L),$$
see 6.21(d), 6.7(b). We obtain
$$\Hom_{\cc^\boc Z_s}({}^1L,\un\z(\un\c(L')))\sub\Hom_{\cc^\boc\tG_s}(\un\c({}^1L),\un\c(L'))
\tag b$$
for any ${}^1L,L'$ in $\cc^\boc_0Z_s$. We show that (b) is an equality:
$$\Hom_{\cc^\boc Z_s}({}^1L,\un\z(\un\c(L')))=\Hom_{\cc^\boc\tG_s}(\un\c({}^1L),\un\c(L')).\tag c$$
Let $N'$ (resp. $N''$) be the dimension of the left (resp. right) hand side of (b). It is enough 
to show that $N'=N''$. We can assume that ${}^1L=\Bbb L_{\l',s}^{\dz'}$, $L'=\Bbb L_{\l,s}^{\dz}$ 
where $z\cdo\l\in\boc^s$, $z''\cdo\l'\in\boc^s$. 
By 6.12(a), $N'$ is the multiplicity of ${}^1L$ in $\un{\fb''}(L')$; by the fully faithfulness of
$\ti\e_s$ this is the same as the multiplicity of $\ti\e_s{}^1L$ in 
$\ti\e_s\un{\fb''}(L')=\un{\fb'}(L')=\un\fb(L')$ (the last two equalities use
4.25(d) and 4.14(d)). By 4.13(d), this is the same as the multiplicity of $\LL_{\l'}^{\dz'}$ in 
$$\op_{y\in W;y\cdo\l\in\boc}\LL_{\ee^{-s}(\l)}^{\ee^{-s}(\dy)}
\un{\cir}\LL_{\l}^{\dz}\un{\cir}\LL_{y(\l)}^{\dy\i}.$$
Using now \cite{\MONO, 2.22(c)} we see that $N'$ is the coefficient of $t_{z'\cdo\l'}$ in
$$\sum_{y\in W;y\cdo\l\in\boc}t_{\ee^{-s}(y)\cdo\ee^{-s}(\l)}t_{z\cdo\l}t_{y\i\cdo y(\l)}
\in\HH^\iy.$$
Hence if $\tt:\HH^\iy@>>>\ZZ$ is as in 1.9, then
$$N'=\sum_{y\in W;y\cdo\l\in\boc}\tt(t_{\ee^{-s}(y)\cdo\ee^{-s}(\l)}t_{z\cdo\l}t_{y\i\cdo y(\l)}
t_{z'{}\i\cdo z'(\l')}).$$
This can be rewritten as
$$N'=\sum_{y\cdo\l_1\in\boc}
\tt(t_{\ee^{-s}(y)\cdo\ee^{-s}(\l_1)}t_{z\cdo\l}t_{y\i\cdo y(\l_1)}t_{z'{}\i\cdo z'(\l')}).$$
(In the last sum, the terms corresponding to $y\cdo\l_1$ with $\l_1\ne\l$ are equal to zero.) By 
6.6(c) (with $z\cdo\l,z'\cdo\l'$ interchanged) we have
$$N''=\sum_{y\cdo\l_1\in\boc}
\tt(t_{\e^{-s}(y)\cdo\ee^{-s}(\l_1)}t_{z\cdo\l}t_{y\i\cdo y(\l_1)}t_{z'{}\i\cdo z'(\l)}).$$
Thus, $N'=N''$. This completes the proof of (c).

\subhead 7.2\endsubhead
Let $s,s'\in\ZZ_\boc$. We define a bifunctor $\cc^\boc\tG_s\T\cc^\boc\tG_{s'}@>>>\cc^\boc\tG_{s+s'}$
 denoted by $K,K'\m K\un{*}K'$ as follows. By replacing if necessary $\Ps$ in 7.1 by a power, we can
assume that
any $A\in CS_{\boc,s}$ and any $A\in CS_{\boc,s'}$ admits a mixed structure (defined in terms of
$\Ps$) of pure weight zero.
Let $K\in\cc^\boc\tG_s$, $K'\in\cc^\boc\tG_{s'}$; we choose mixed structures of pure weight $0$
on $K,K'$ with respect to $\Ps$ (this is possible by our choice of $\Ps$). We define $K\un{*}K'$ 
as in 6.19 in terms of these mixed structures and we then disregard the mixed structure on 
$K\un{*}K'$. The resulting object of $\cc^\boc\tG_{s+s'}$ is denoted again by $K\un{*}K'$; it is 
independent of the choice made.

In the same way the functor $\un\c:\cc^\boc_0Z_s@>>>\cc^\boc_0\tG_s$ gives rise to a functor
$\cc^\boc Z_s@>>>\cc^\boc\tG_s$ denoted again by $\un\c$; the functor 
$\un\z:\cc^\boc_0\tG_s@>>>\cc^\boc_0Z_s$ gives rise to a functor $\cc^\boc\tG_s@>>>\cc^\boc Z_s$ 
denoted again by $\un\z$.

The operation $K\un{*}K'$ is again called truncated convolution. It has a canonical associativity
isomorphism (deduced from that in 6.19(d)); this makes $\sqc_{s\in\ZZ_\boc}\cc^\boc\tG_s$ into a 
monoidal category.

From 6.20 we see that under
$\un\z:\sqc_{s\in\ZZ_\boc}\cc^\boc\tG_s@>>>\sqc_{s\in\ZZ_\boc}\cc^\boc Z_s$,
the monoidal structure on $\sqc_{s\in\ZZ_\boc}\cc^\boc\tG_s$ is compatible with the opposite of the
monoidal structure on $\sqc_{s\in\ZZ_\boc}\cc^\boc Z_s$.

If $K\in\cc^\boc\tG_s$ then the isomorphisms 6.13(b) provide an $\ee^s$-half-braiding for 
$\ti\e_s\un\z(K)\in \cc^\boc\tcb^2$ so that $\ti\e_s\un\z(K)$ can be naturally viewed as an object
of $\cz^\boc_{\ee^s}$ denoted by $\ov{\ti\e_s\un\z(K)}$. (Note that 6.13(b) is stated in the mixed
category but it implies the corresponding result in the unmixed category.) Then 
$K\m\ov{\ti\e_s\un\z(K)}$ is a functor $\cc^\boc\tG_s@>>>\cz^\boc_{\ee^s}$.

\proclaim{Theorem 7.3}Let $s\in\ZZ_\boc$. The functor $\cc^\boc\tG_s@>>>\cz^\boc_{\ee^s}$,
$K\m\ov{\ti\e_s\un\z(K)}$ is an equivalence of categories.
\endproclaim
From 6.12(a),4.14(d),4.25(d) we have canonically for any $z\cdo\l\in\boc^s$:
$$\ti\e_s\un\z(\un\c(\Bbb L_{\l,s}^{\dz}))=\un\fb(\Bbb L_{\l,s}^{\dz})\tag a$$
as objects of $\cc^\boc\tcb^2$. 
From the definitions we see that the $\ee^s$-half-braiding on the left hand side of (a) provided 
by 7.2 is the same as the $\ee^s$-half-braiding on the right hand side of (a) provided by 4.14(j).
Hence we have
$$\ov{\ti\e_s\un\z(\un\c(\Bbb L_{\l,s}^{\dz}))}=\ov{\un\fb(\Bbb L_{\l,s}^{\dz})}\tag b$$
as objects of $\cz^\boc_{\ee^s}$. Using this and 5.7(a) with
$L'=\ti\e_s\un\z(\un\c(\Bbb L_{\l',s}^{\dw}))$ 
(where $z\cdo\l,w\cdo\l'$ are in $\boc^s$), we have
$$\Hom_{\cc^\boc\tcb^2}(\LL_\l^{\dz},\ti\e_s\un\z(\un\c(\Bbb L_{\l',s}^{\dw})))=       
\Hom_{\cz^\boc_{\ee^s}}(\ov{\ti\e_s\un\z(\un\c(\Bbb L_{\l,s}^{\dz}))},
\ov{\ti\e_s\un\z(\un\c(\Bbb L_{\l',s}^{\dw}))}).$$
Combining this with the equalities
$$\align&\Hom_{\cc^\boc\tG_s}(\un\c(\Bbb L_{\l,s}^{\dz}),\un\c(\Bbb L_{\l',s}^{\dw}))=
\Hom_{\cc^\boc Z_s}(\Bbb L_{l,s}^{\dz},\un\z(\un\c(\Bbb L_{\l',s}^{\dw})))\\&=
\Hom_{\cc^\boc\tcb^2}(\LL_l^{\dz},\ti\e_s\un\z(\un\c(\Bbb L_{\l',s}^{\dw}))),\endalign$$
of which the first comes from 6.10(c) and the second comes from the fully faithfulness of 
$\ti\e_s$, we obtain
$$\Hom_{\cc^\boc\tG_s}(\un\c(\Bbb L_{\l,s}^{\dz}),\un\c(\Bbb L_{\l',s}^{\dw}))=
\Hom_{\cz^\boc_{\ee^s}}(\ov{\ti\e_s\un\z(\un\c(\Bbb L_{\l,s}^{\dz}))},
\ov{\ti\e_s\un\z(\un\c(\Bbb L_{\l',s}^{\dw}))}).$$
In other words, setting 
$$\AA_{z\cdo\l,w\cdo\l'}=\Hom_{\cc^\boc\tG_s}(\un\c(\Bbb L_{\l,s}^{\dz}),
\un\c(\Bbb L_{\l',s}^{\dw})),$$
$$\AA'_{z\cdo\l,w\cdo\l'}=\Hom_{\cz^\boc_{\ee^s}}(\ov{\ti\e_s\un\z(\un\c(\Bbb L_{\l,s}^{\dz}))},
\ov{\ti\e_s\un\z(\un\c(\Bbb L_{\l',s}^{\dw}))}),$$
we have 
$$\AA_{z\cdo\l,w\cdo\l'}=\AA'_{z\cdo\l,w\cdo\l'}.\tag c$$
Note that the identification (c) is induced by the functor $K\m\ov{\ti\e_s\un\z(K)}$.
Let $\AA=\op\AA_{z\cdo\l,w\cdo\l'}$, $\AA'=\op\AA_{z\cdo\l,w\cdo\l'}$ (both direct sums are taken 
over all $z\cdo\l,w\cdo\l'$ in $\boc^s$). Then from (c) we have $\AA=\AA'$. Note that this
identification is compatible with the obvious algebra structures of $\AA,\AA'$.

For any $A\in CS_{\boc,s}$ we denote by $\AA_A$ the set of all $f\in\AA$ such that for any
$z\cdo\l,w\cdo\l'$, the $(z\cdo\l,w\cdo\l')$-component of $f$ maps the $A$-isotypic component of 
$\un\c(\Bbb L_{\l,s}^{\dz})$ to the $A$-isotypic component of $\un\c(\Bbb L_{\l',s}^{\dw})$ and 
any other isotypic component of $\un\c(\Bbb L_{\l,s}^{\dz})$ to $0$. Thus, 
$\AA=\op_{A\in CS_{\boc,s}}\AA_A$ is the decomposition of $\AA$ into a sum of simple algebras. 
(Each $\AA_A$ is nonzero since, by 6.2(c) and 6.5(a), any
$A$ is a summand of some $\un\c(\Bbb L_{\l,s}^{\dz})$.)

Let $\fS$ be a set of representatives for the isomorphism classes of simple objects of 
$\cz^\boc_{\ee^s}$. For any $\s\in\fS$ we denote by $\AA'_\s$ the set of all $f'\in\AA'$ such that
for any $z\cdo\l,w\cdo\l'$, the $(z\cdo\l,w\cdo\l')$-component of $f'$ maps the $\s$-isotypic 
component of $\ov{\ti\e_s\un\z(\un\c(\Bbb L_{\l,s}^{\dz}))}$ to the $\s$-isotypic component of 
$\ov{\ti\e_s\un\z(\un\c(\Bbb L_{\l',s}^{\dw}))})$ and all other isotypic components of 
$\ov{\ti\e_s\un\z(\un\c(\Bbb L_{\l,s}^{\dz}))}$ to zero. Then $\AA'=\op_{\s\in\fS}\AA'_\s$ is the 
decomposition of $\AA'$ into a sum of simple algebras. (Each $\AA'_\s$ is nonzero since any $\s$ 
is a summand of some $\ov{\ti\e_s\un\z(\un\c(\Bbb L_{\l,s}^{\dz}))}$ with $z\cdo\l\in\boc^s$.
Indeed, we can find $z\cdo\l\in\boc$ such that $\LL_\l^{\dz}$ is a direct summand of $\s$, viewed 
as an object of $\cc^\boc\tcb^2$; then, by 5.5(a), $\s$ is a summand of 
$\ov{\ci_s(\LL_\l^{\dz})}$. If in addition, $z\cdo\l\in\boc^s$ then, by 5.6(a),(b), we have 
$\ov{\ci_s(\LL_\l^{\dz})}=\ov{\un\fb(\Bbb L_{\l,s}^{\dz})}$ hence $\s$ is a summand of 
$\ov{\un\fb(\Bbb L_{\l,s}^{\dz})}$ hence, by (a), $\s$ is a summand of 
$\ov{\ti\e_s\un\z(\un\c(\Bbb L_{\l,s}^{\dz}))}$, as required. If $z\cdo\l\n\boc^s$ then, by 
5.5(b), we have $\ci_s(\LL_\l^{\dz})=0$ which is a contradiction.)
Since $\AA=\AA'$, from the uniqueness of decomposition of a semisimple algebra as a direct sum of
simple algebras, we see that there is a unique bijection $CS_{\boc,s}\lra\fS$, $A\lra\s_A$ such 
that $\AA_A=\AA'_{\s_A}$ for any $A\in CS_{\boc,s}$. From the definitions we now see that for any 
$A\in CS_{\boc,s}$ we have $\ov{\ti\e_s\un\z(K)}\cong\s_A$. Therefore, Theorem 7.3 holds.

\proclaim{Theorem 7.4}We preserve the setup of Theorem 7.3.
Let $L\in\cc^\boc Z_s$, $K\in\cc^\boc\tG_s$. We have canonically
$$\Hom_{\cc^\boc Z_s}(L,\un\z(K))=\Hom_{\cc^\boc\tG_s}(\un\c(L),K).\tag a$$
\endproclaim
We can assume that $L=\Bbb L_{\l,s}^{\dz}$ where $z\cdo\l\in\boc^s$. From 7.3 and its proof we 
see that
$$\Hom_{\cc^\boc\tG_s}(\un\c(L),K)
=\Hom_{\cz^\boc_{\ee^s}}(\ov{\ti\e_s\un\z(\un\c(L))},\ov{\ti\e_s\un\z(K)})=
\Hom_{\cz^\boc_{\ee^s}}(\ov{\ci_s(\LL_\l^{\dz})},\ov{\ti\e_s\un\z(K)}).$$
Using 5.5(a) we see that 
$$\Hom_{\cz^\boc_{\ee^s}}(\ov{\ci_s(\LL_\l^{\dz})},\ov{\ti\e_s\un\z(K)})
\Hom_{\cc^\boc\tcb^2}(\LL_\l^{\dz},\ti\e_s\un\z(K))=\Hom_{\cc^\boc Z_s}(L,\un\z(K)).$$
This proves the theorem.

\subhead 7.5\endsubhead
We preserve the setup of Theorem 7.3. We show that for $K\in\cc^\boc\tG_s$ we have canonically
$$\fD(\un\z(\fD(K))))=\un\z(K).\tag a$$
Here the first $\un\z$ is relative to $\wt\boc$. It is enough to show that for any 
$L\in\cc^\boc Z_s$ we have canonically
$$\Hom_{\cc^\boc Z_s}(L,\fD(\un\z(\fD(K)))))=\Hom_{\cc^\boc Z_s}(L,\un\z(K)).$$
Here the left side equals
$$\align&\Hom_{\cc^{\wt{\boc}}Z_s}(\un\z(\fD(K)),\fD(L))=\Hom_{\cc^\boc\tG_s}(\fD(K),
\un\c(\fD(L)))\\&=\Hom_{\cc^\boc\tG_s}(\fD(K),\fD(\un\c(L))).\endalign$$
(We have used 7.4(a) for $\ti\boc$ and 6.7(b).) The right hand side equals
$$\Hom_{\cc^\boc\tG_s}(\un\c(L),K)=\Hom_{\cc^\boc\tG_s}(\fD(K),\fD(\un\c(L))).$$
(We have again used 7.4(a).) This proves (a).

\proclaim{Theorem 7.6} Let $s\in\ZZ_\boc$. Let $K\in\cc^\boc\tG_s$. In $\cc^\boc\tcb^2$ we
have 
$$\ti\e_s\un\z(K)\cong\op_{z\cdo\l\in\boc^s;z\cdo\l
\underset\text{left}\to\si\ee^s(z\i)\cdo\l}(\LL_\l^{\dz})^{\op N(z,\l)}$$
where $N(z,\l)\in\NN$.
\endproclaim
In $\cc^\boc Z_s$ we have
$$\un\z(K)\cong\op_{z\cdo\l\in\boc^s}(\Bbb L_{\l,s}^{\dz})^{\op N(z,\l)}\tag a$$
where $N(z,\l)\in\NN$. If $N(z,\l)>0$ then 
$$\Hom_{\cc^\boc Z_s}(\Bbb L_{\l,s}^{\dz},\un\z(K))\ne0$$
hence by 7.4 we have $\Hom_{\cc^\boc\tG_s}(\un\c(\Bbb L_{\l,s}^{\dz}),K)\ne0$
and in particular $\un\c(\Bbb L_{\l,s}^{\dz})\ne0$. Using 6.5(d) we deduce that

(b) $z\cdo\l\underset\text{left}\to\si\ee^s(z\i)\cdo\l$. 
\nl
Thus the direct sum in (a) can be restricted to $z\cdo\l$ satisfying (b). We now apply
$\ti\e_s$ to both sides of (a) and use that $\ti\e_s\Bbb L_{\l,s}^{\dz}=\LL_\l^{\dz}$. The 
theorem follows.

\subhead 7.7\endsubhead
Let $s\in\ZZ_\boc$. From 7.3 and 7.6 we see that any object of $\cz^\boc_{\ee^s}$, when 
viewed as an object of $\cc^\boc\tcb^2$, is a direct sum of objects of the form 
$\LL_\l^{\dz}$ with $z\cdo\l\in\boc^s$ such that
$z\cdo\l\underset\text{left}\to\si\ee^s(z\i)\cdo\l$.

In the remainder of this subsection we assume that $\tG$ is as in case A with $G$
simple of type $A_2$ (resp. $B_2$ or $G_2$). In this case $W$ is generated by 
$\s_1,\s_2$ in $S$ with relation $(\s_1\s_2)^m=1$ where $m=3$ (resp. $m=4$ or $m=6$). We 
assume that $\boc$ is the two-sided cell of $I$ consisting of all $w\cdo1$ where $w\in W$, 
$1\le|w|\le m-1$. We shall write $\LL^{iji\do}$ instead of 
$\LL^{\dot\s_i\dot\s_j\dot\s_i\do}_1$ where $iji\do$ is $121\do$ or $212\do$.
The objects of $\cc^\boc\tcb^2$ of the form $\ti\e_s\un\z(K)$ with $K$ a simple object of 
$\cc^\boc\tG_s$ are (up to isomorphism) the following ones:
$$\LL^1\op\LL^2\text{ for type }A_2;$$
$$\LL^1\op\LL^2, \LL^1\op\LL^{212}, \LL^2\op\LL^{121}, \LL^{121}\op\LL^{212}
\text{ for type }B_2;$$
$$\align&\LL^1\op\LL^2, \LL^1\op\LL^2\op\LL^{121}\op\LL^{212},
\LL^2\op\LL^{121}\op\LL^{21212},\\& \LL^1\op\LL^{212}\op\LL^{12121}, 
\LL^{121}\op\LL^{212}\op\LL^{12121}\op\LL^{21212}, \LL^{121}\op\LL^{212}\text{ for type }
G_2.\endalign$$
Note that in type $G_2$, $\LL^{121}\op\LL^{212}$ comes from two nonisomorphic objects
$K$ of $\cc^\boc\tG_s$.

\subhead 7.8\endsubhead
In this subsection we assume that $\tG$ is as in case A with $G=SL_2(\kk)$ and $p\ne2$. 
In this case we may identify $\TT=\kk^*$ and $W=\{1,\s\}$ with $\s(t)=t\i$ for $t\in\TT$.
We take $\t\in\tG_1$ such that $\ee:G@>>>G$ in 2.3 satisfies $\ee(t)=t^q$ for any $t\in T$.
Then for $\l\in\fs_\iy\cong\kk^*$ we have $\ee(\l)=\l^{q\i}$, $\s(\l)=\l\i$.
Let $\l_0$ be the unique element of $\fs_\iy$ such that $\l_0^2=1,\l_0\ne1$.
In $\HH$ we have $c_{1\cdo\l}=T_11_\l$ for all $\l$, $c_{\s\cdo\l}=v\i T_\s1_\l$ if $\l\ne1$, 
$c_{\s\cdo 1}=v\i T_\s1_1+v\i T_11_1$. It follows that the two-sided cells in 
$I=\{w\cdo\l;w\in W,\l\in\fs_\iy\}$ are the following subsets of $I$:

$\boc_\l=\boc_{\l\i}=\{1\cdo\l,1\cdo\l\i,\s\cdo\l,\s\cdo\l\i\}$ with $\l\in\fs_\iy;\l^2\ne1$;

$\boc_{\l_0}=\{1\cdo\l_0,\s\cdo\l_0\}$;

$\boc'_1=\{\s\cdo1\}$;

$\boc_1=\{1\cdo1\}$.

Let $s\in\ZZ$. The two-sided cells of $I$ which are stable under $\ee^s$ are:

(i) $\boc_\l=\boc_{\l\i}$ where $\l\in\fs_\iy$, $\l^2\ne1$, $\l^{q^{-s}}=\l$ (note that $\ee^s$ 
acts as $1$ on this two-sided cell);

(ii) $\boc_\l=\boc_{\l\i}$ where $\l\in\fs_\iy$, $\l^2\ne1$, $\l^{q^{-s}}=\l\i$ (note that $\ee^s$ 
acts as a fixed point free involution on this two-sided cell and that we have necessarily $s\ne0$);

(iii) $\boc_{\l_0}$ (note that $\ee^s$ acts as $1$ on this two-sided cell);

(iv) $\boc'_1$ (note that $\ee^s$ acts as $1$ on this two-sided cell);

(v) $\boc_1$ (note that $\ee^s$ acts as $1$ on this two-sided cell).

For $\boc$ in (i)-(v), the $\ee^s$-centre of $\cc^\boc\tcb^2$ has exactly $N$ simple objects (up to
isomorphism) where $N=1$ in the cases (i),(ii),(iv),(v) and $N=4$ in the case (iii).


\widestnumber\key{ENO}
\Refs
\ref\key\BBD \by A.Beilinson, J.Bernstein and P.Deligne\paper Faisceaux pervers\jour Ast\'erisque
\vol100\yr1982\endref
\ref\key\BFO\by R.Bezrukavnikov, M.Finkelberg and V.Ostrik\paper Character D-modules via Drinfeld 
center of Harish-Chandra bimodules\jour Invent. Math.\vol188\yr2012\pages589-620\endref
\ref\key\DL\by P.Deligne and G.Lusztig\paper Representations of reductive groups over finite fields
\jour Ann. Math\vol103\yr1976\pages103-161\endref
\ref\key\ENO\by P.Etingof, D.Nikshych and V.Ostrik\paper On fusion categories\jour Ann. Math.
\vol162\yr2005\pages581-642\endref
\ref\key\KL\by D.Kazhdan and G.Lusztig\paper Representations of Coxeter groups and Hecke algebras
\jour Invent.Math.\vol53\yr1979\pages165-184\endref
\ref\key\ORA\by G.Lusztig\book Characters of reductive groups over a finite field\bookinfo 
Ann.Math. Studies 107\publ Princeton Univ.Press\yr1984\endref
\ref\key\ICM\by G.Lusztig\paper Characters of reductive groups over finite fields\inbook
Proc.Int.Congr.Math. Warsaw 1983\publ North Holland\yr1984\pages877-880\endref
\ref\key\CELLSII\by G.Lusztig\paper Cells in affine Weyl groups II\jour J.Alg.\vol109\yr1987
\pages536-548\endref
\ref\key\CSI\by G.Lusztig\paper Character sheaves I\jour Adv. Math.\vol56\yr1985\pages193-237
\endref
\ref\key\CSII\by G.Lusztig\paper Character sheaves II\jour Adv. Math.\vol57\yr1985\pages226-265
\endref
\ref\key\CSIII\by G.Lusztig\paper Character sheaves III\jour Adv. Math.\vol57\yr1985\pages266-315
\endref
\ref\key\TENS\by G.Lusztig\paper Cells in affine Weyl groups and tensor categories\jour Adv.Math.
\vol129\yr1997\pages85-98\endref
\ref\key\HEC\by G.Lusztig\book Hecke algebras with unequal parameters\bookinfo CRM Monograph Ser.18
\publ Amer. Math. Soc.\yr2003\endref
\ref\key\CDGVI\by G.Lusztig\paper Character sheaves on disconnected groups VI\jour Represent. Th.
\vol8\yr2004\pages377-413\endref
\ref\key\CDGVII\by G.Lusztig\paper Character sheaves on disconnected groups VII\jour Represent. Th.
\vol9\yr2005\pages209-266\endref
\ref\key\CDGVIII\by G.Lusztig\paper Character sheaves on disconnected groups VIII\jour 
Represent. Th. \vol10\yr2006\pages314-352\endref
\ref\key\CDGIX\by G.Lusztig\paper Character sheaves on disconnected groups IX\jour Represent. Th.
\vol10\yr2006\pages353-379\endref
\ref\key\CDGX\by G.Lusztig\paper Character sheaves on disconnected groups X\jour Represent. Th.
\vol13\yr2009\pages82-140\endref
\ref\key\CONV\by G.Lusztig\paper Truncated convolution of character sheaves\jour Bull. Inst. Math.
Acad. Sin. (N.S.)\vol10\yr2005\pages1-72\endref
\ref\key\URE\by G.Lusztig Unipotent representations as a categorical centre\jour Represent. Th.\vol19\yr2015\pages211-235
\endref
\ref\key\MONO\by G.Lusztig Non-unipotent character sheaves as a categorical centre\jour 
Bull. Inst. Math.Acad. Sin. (N.S.)\vol11\yr2016\pages603-731\endref
\ref\key\MUG\by M.M\"uger\paper From subfactors to categories and topology II. The quantum double 
of tensor categories and subfactors\jour J. Pure Appl. Alg.\vol180\yr2003\page159-219\endref
\ref\key\YO\by T.Yokonuma\paper Sur la structure des anneaux de Hecke d'un groupe de Chevalley fini
\jour C.R. Acad.Sci. Paris Ser.A\vol264\yr1967\pages A334-A347\endref
\endRefs
\enddocument